\documentclass{article}

\usepackage{arxiv}
\usepackage{lipsum}
\usepackage{amsmath}
\usepackage{xr-hyper}
\usepackage{hyperref}
\usepackage{graphics,graphicx}
\usepackage{epstopdf}
\usepackage[amsmath,thmmarks,hyperref]{ntheorem}
\usepackage{amsfonts}
\usepackage{graphicx}
\usepackage{epstopdf}
\usepackage{algorithm}
\usepackage{algorithmic}
\usepackage[algo2e,linesnumbered,lined, ruled]{algorithm2e}
\usepackage[capitalize,nameinlink]{cleveref}
\usepackage{amssymb}
\usepackage{nicefrac}
\usepackage{bm} 
\usepackage{mathtools}
\usepackage{pgfplotstable} 
\usepackage{booktabs}
\usepackage{multicol}
\usepackage{colortbl}
\usepackage{makecell}
\usepackage{siunitx}
\usepackage{ifpdf}

\usepackage[caption=false]{subfig} 

\Crefname{figure}{Figure}{Figures}

\crefformat{equation}{\textup{#2(#1)#3}}
\crefrangeformat{equation}{\textup{#3(#1)#4--#5(#2)#6}}
\crefmultiformat{equation}{\textup{#2(#1)#3}}{ and \textup{#2(#1)#3}}
{, \textup{#2(#1)#3}}{, and \textup{#2(#1)#3}}
\crefrangemultiformat{equation}{\textup{#3(#1)#4--#5(#2)#6}}%
{ and \textup{#3(#1)#4--#5(#2)#6}}{, \textup{#3(#1)#4--#5(#2)#6}}{, and \textup{#3(#1)#4--#5(#2)#6}}

\Crefformat{equation}{#2Equation~\textup{(#1)}#3}
\Crefrangeformat{equation}{Equations~\textup{#3(#1)#4--#5(#2)#6}}
\Crefmultiformat{equation}{Equations~\textup{#2(#1)#3}}{ and \textup{#2(#1)#3}}
{, \textup{#2(#1)#3}}{, and \textup{#2(#1)#3}}
\Crefrangemultiformat{equation}{Equations~\textup{#3(#1)#4--#5(#2)#6}}%
{ and \textup{#3(#1)#4--#5(#2)#6}}{, \textup{#3(#1)#4--#5(#2)#6}}{, and \textup{#3(#1)#4--#5(#2)#6}}

\theoremstyle{plain}

\newtheorem{theorem}{Theorem}[section]

\newtheorem{definition}{Definition}[section]

\ifpdf
  \DeclareGraphicsExtensions{.eps,.pdf,.png,.jpg}
\else
  \DeclareGraphicsExtensions{.eps}
\fi

\DeclareSIUnit\Dalton{Da}

\usepackage{enumitem}
\setlist[enumerate]{leftmargin=.5in}
\setlist[itemize]{leftmargin=.5in}

\pgfplotsset{compat = 1.16}

\definecolor{myGreen}{RGB}{25,142,33}

\newcommand{\R}{\mathbb{R}}

\newcommand{\ra}[1]{\renewcommand{\arraystretch}{#1}}

\DeclareMathOperator{\dom}{dom}
\DeclareMathOperator{\divergence}{div}
\DeclareMathOperator{\TV}{TV}
\DeclareMathOperator{\TVeps}{\TV_{\varepsilon_{\TV}}}
\DeclareMathOperator{\prox}{prox}

\DeclareMathOperator*{\argmin}{arg\,min}

\newcommand{\email}[1]{\protect\href{mailto:#1}{#1}}

\newcommand{\diff}{\mathop{} \! \mathrm{d}}

\crefname{hypothesis}{Hypothesis}{Hypotheses}

\Crefname{ALC@unique}{Line}{Lines} 

\crefname{algocf}{alg.}{algs.}
\Crefname{algocf}{Algorithm}{Algorithms}

\title{Spatially Coherent Clustering Based on Orthogonal Nonnegative Matrix Factorization}

\author{ \href{https://orcid.org/0000-0003-3424-8031}{\includegraphics[scale=0.06]{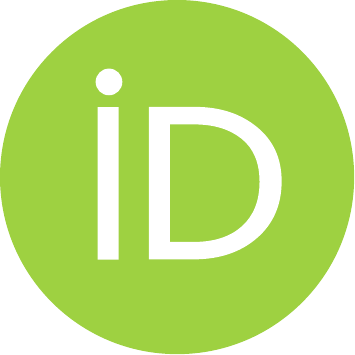}\hspace{1mm}Pascal~Fernsel}\\
        Center for Industrial Mathematics\\
        University of Bremen\\
        Bremen, Germany \\
	\email{p.fernsel@uni-bremen.de} \\
}

\usepackage{amsopn}
\DeclareMathOperator{\dist}{dist}

\renewcommand{\shorttitle}{Spatially Coherent Clustering Based on ONMF}

\ifpdf
\hypersetup{
  pdftitle={Spatially Coherent Clustering Based on Orthogonal Nonnegative Matrix Factorization},
  pdfauthor={P. Fernsel},
  pdfkeywords={Orthogonal nonnegative matrix factorization, Clustering, Spatial coherence, Hyperspectral data, MALDI imaging}
}
\fi

\begin{document}

    \maketitle
    
    \begin{abstract}
        Classical approaches in cluster analysis are typically based on a feature space analysis. However, many applications lead to datasets with additional spatial information and a ground truth with spatially coherent classes, which will not necessarily be reconstructed well by standard clustering methods. Motivated by applications in hyperspectral imaging, we introduce in this work clustering models based on orthogonal nonnegative matrix factorization, which include an additional total variation (TV) regularization procedure on the cluster membership matrix to enforce the needed spatial coherence in the clusters. We propose several approaches with different optimization techniques, where the TV regularization is either performed as a subsequent postprocessing step or included into the clustering algorithm. Finally, we provide a numerical evaluation of all proposed methods on a hyperspectral dataset obtained from a matrix-assisted laser desorption/ionisation imaging measurement, which leads to significantly better clustering results compared to classical clustering models.
    \end{abstract}


    \AMSClass{15A23 \and 65F22 \and 65F55}

    \section{Introduction} \label{sec: Introduction}
        Cluster analysis has been studied over fifty years in the machine learning community and is one of the central topics in unsupervised learning with a wide range of possible research directions and application fields including image segmentation, document clustering and bioinformatics \cite{AggarwalReddy:2013:BookClustering}. The general clustering problem is to partition a given set of objects $\mathcal{O}$ into different groups, such that each object within a group are more similar to each other compared to the objects in other groups. One typical approach is based on a feature space analysis. The basic concept is to assign to each object $\sigma\in \mathcal{O}$ a feature vector $x_\sigma\in \mathcal{X}$ containing the characteristics of $\sigma,$ where $X$ is a suitably defined feature space. Furthermore, a similarity measure and a suitable minimization problem is defined to introduce the notion of similarity between the feature vectors and to formulate the clustering problem.
        
        However, many types of datasets contain additional spatial information, which is typically not used in cluster analysis. Characteristic examples are images or, more generally, hyperspectral datasets, where each measured spectrum is associated to a point in a two- or three-dimensional space. Furthermore, many application fields like mass spectrometry imaging or Earth remote sensing naturally lead to datasets with spatially coherent regions. Hence, a classical cluster analysis, which is entirely based on a feature space analysis, does not lead necessarily to spatially coherent clusters and is therefore not sufficient to reconstruct the spatial coherent regions in these kind of data.
        
        Hence, this work focuses on a combined clustering analysis, which takes into account both the feature space and the spatial coherence of the clusters. We introduce numerous clustering methods based on orthogonal nonnegative matrix factorization (ONMF) for general nonnegative datasets with spatial information and include a TV regularization procedure to regularize the cluster membership matrix to induce the needed spatial coherence. Furthermore, we discuss different optimization techniques for the ONMF models and derive the corresponding minimization algorithms. Finally, we perform a numerical evaluation on a mass spectrometry imaging dataset acquired from a matrix-assisted laser desorption/ionization imaging measurement of a human colon tissue sample and compare the proposed clustering methods to classical ONMF approaches.
        
        This paper is organized as follows. After a short description of the related work and the used notation in \cref{subsec: Related Work} and \cref{subsec: Notation}, we give a short outline of the basics of ONMF approaches, its relations to K-means clustering and details on possible solution algorithms in \cref{sec: Background}. In \cref{sec: Orthogonal NMF with Spatial Coherence}, we introduce the proposed methods in this work which are divided into so-called separated methods and combined methods. \cref{sec: Numerical Experiments} is entirely devoted to the numerical experiments and the evaluation of the discussed methods. Finally, \cref{sec: Conclusion} concludes the findings and gives an outlook for future possible research directions.
    
        \subsection{Related Work} \label{subsec: Related Work}
            The natural relation between ONMF and clustering models is well studied.
            
            One of the first theoretical analysis was provided by Ding et al.\ in \cite{Ding:2005:NMFKMeans}. By comparing the cost functions of different nonnegative matrix factorization (NMF) and $K$-means models, the authors could show in their work for example the strong relationship between $K$-means clustering and ONMF with an orthogonality constraint on one of the factorization matrices as well as kernel $K$-means and symmetric ONMF. Furthermore, the connections to spectral clustering and tri-factorizations were studied. Several works followed with a similar theoretical emphasis on tri-factorizations \cite{Ding:2006:triFac} and multiple other NMF models \cite{LiDing:2006:Relationships}. 
            
            The previous mentioned works focus on the theoretical side but also give some first update rules to solve the corresponding ONMF problems. However, more work has been done for the algorithm development. Many classical approaches are based on multiplicative update rules \cite{Ding:2006:triFac,Choi:2008:CompMethod,YangOja:2010:CompMethod,LiWuPeng:2010:CompMethod}. More recent works are, for instance, based on nuclear norm optimization \cite{PanNg:2018:CompMethod}, further multiplicative update schemes \cite{Mirzal:2014:CompMethod,ZhangEtal:2016:CompMethodNMFLPP}, hierarchical alternating least squares \cite{Kimura:2015:CompMethod,LiEtal:2015:CompMethod,LiEtal:2020:CompMethod}, proximal alternating linearized Minimization \cite{WangEtal:2019:CompMethod}, EM like algorithms and augmented Lagrangian approaches \cite{Pompili:2014:CompMethod}, deep neural networks \cite{QiuEtal:2017:DeepNMFCompMethod} and other techniques \cite{AsterisEtal:2015:CompMethod,ZhangEtal:2016:CompMethodStiefelManifold}. Finally, we would like to refer the interested reader at this point to two review articles on NMF methods for clustering \cite{LiDing:2014:NMFClusteringSurvey,Turkmen:2015:NMFClusteringSurvey} and a book on NMF and nonnegative tensor factorizations \cite{Cichocki:2009:BookNMF}. 
            
            While analyzing and developing optimization algorithms for NMF clustering methods is a major topic throughout the literature, studying clustering techniques with spatial coherence by incorporating the local information of the considered datapoints is far less common. Spatially coherent clusterings are primarily analyzed in the context of image segmentation problems \cite{DespotovicEtal:2013:ImgSegment,ZabihKolmogorov:2004:ImgSegment,HuangEtal:2011:ImgSegment,Mignotte:2011:ImgSegment}. 
            The subject of spatial coherence can also be found in the literature of hyperspectral image analysis. Several NMF models with total variation regularization, which include the local neighbourhood information of each datapoint, have been analyzed for the critical processing step of hyperspectral unmixing \cite{HeEtal:2017:HypUnmix,FengEtal:2018:HypUnmix,FengEtal:2019:HypUnmix}. These articles can be considered as closest to our approach, since we also focus on the application of generalized NMF models to hyperspectral images. Further works also consider NMF models with different TV penalty terms for hyperspectral image denoising \cite{ZhangEtal:2008:NMFTV,YinLiu:2010:NMFTV,LengEtal:2017:PropMethod}, some of which will be used in the later course of this work for the derivation of optimization algorithms of the respective clustering models.
            
            However, all the aforementioned works only include either orthogonality constraints or TV regularization into their NMF models, whereas we focus on combining both of these properties to obtain a spatially coherent clustering method for hyperspectral datasets.
            
            The only work which includes TV regularization as well as penalty terms to enforce an orthogonality constraint on one of the matrices is, to the best of our knowledge, the survey article \cite{Fernsel:2018:Survey} with a rather general focus on the development of algorithms based on surrogate functions leading to multiplicative update schemes.
            
            
            
        
        \subsection{Notation} \label{subsec: Notation}
            Matrices will play a major role throughout this work and are denoted, unless otherwise stated, by capital Latin or Greek letters (e.g.\ $X, U, \Psi, \dots$). The entry of a matrix $U$ in the $i$-th row and $j$-th column is indicated as $U_{ij}.$ The same holds true for a matrix product, where its $ij$-th entry is given by $(UV)_{ij}.$ Furthermore, we use a dot to indicate rows and columns of matrices. The $i$-th row and $j$-th column of $U$ are written as $U_{i, \bullet}$ and $U_{\bullet,j}$ respectively. Moreover, we denote the $i$-th iteration of a matrix $U$ in an algorithm by $U^{[i]}.$
            
            What is more, we write $\Vert U\Vert_F$ and $\Vert U_{\bullet,j} \Vert_2$ for the Frobenius norm of a matrix $U$ and the usual Euclidean norm of a vector $U_{\bullet,j}.$ We also use the notion of nonnegative matrices and write $U\geq 0$ or $U\in \R_{\geq 0}^{m\times n} $ with $\R_{\geq 0} \coloneqq \{ x\in \R \ \vert \ x \geq 0\}$ for an $m\times n$ matrix $U,$ which has only nonnegative entries. The notation for the dimension of the matrices in the NMF problems are reused throughout the article and will be introduced in the following \cref{sec: Background}.
        
    
    
    \section{Background} \label{sec: Background}
        \subsection{Orthogonal NMF and K-means} \label{subsec:Orthogonal NMF and K-Means}
            Nonnegative matrix factorization (NMF), originally introduced by Paatero and Tapper in 1994 \cite{PaateroTapper:1994:PositiveMF} as positive matrix factorization, is a specific matrix factorization method designed to obtain a low-rank approximation of a given and typically large nonnegative data matrix. Different from the widely used principal component analysis (PCA), which is based on the singular value decomposition and allows to compute a best rank $K$ approximation of a given arbitrary matrix, the NMF constraints the matrix factors to be nonnegative. This property makes the NMF the method of choice where the considered data naturally fulfills a nonnegativity constraint so that the interpretability of the factor matrices is ensured. NMF has been widely used for data compression, source separation, feature extraction, clustering or even for solving inverse problems. Possible application fields are hyperspectal unmixing \cite{HeEtal:2017:HypUnmix,FengEtal:2018:HypUnmix,FengEtal:2019:HypUnmix}, document clustering \cite{KimPark:2008:SparseNMF,PanNg:2018:CompMethod}, music analysis \cite{Fevotte:2009:NMFforMusic} but also medical imaging problems such as dynamic computed tomography to perform a joint reconstruction and low-rank decomposition of the corresponding dynamic inverse problem \cite{Arridge:2020:NMFForIP} or matrix-assisted laser desorption/ionisation (MALDI) imaging, where it can be used for tumor typing in the field of bioinformatics as a supervised classification method \cite{Leuschner:2018:NMFSupervised}.
            
            Mathematically, the standard NMF problem can be formulated as follows: For a given nonnegative matrix $X\in \R_{\geq 0}^{M\times N},$ the task is to find two nonnegative matrices $U\in \R_{\geq 0}^{M\times K}$ and $V\in \R_{\geq 0}^{K\times N}$ with $K \ll \min\{M, N\},$ such that
            \begin{equation*}
                X \approx UV = \sum_{k=1}^K U_{\bullet, k} V_{k, \bullet}.
            \end{equation*}
            This allows the approximation of the columns $ X_{\bullet,n } $ and rows $ X_{m,\bullet } $ via a superposition of just a few basis vectors $\{ U_{\bullet, k} \}_k$ and $\{ V_{k, \bullet} \}_k, $ such that $X_{\bullet,n } \approx \sum_k V_{kn} U_{\bullet, k} $ and  $X_{m,\bullet } \approx \sum_k U_{mk} V_{k, \bullet}.$ In this way, the NMF can be seen as a basis learning tool with additional nonnegativity constraints.
    
            The typical variational approach to tackle the NMF problem is to reformulate it as a minimization problem by defining a suitable discrepancy term $\mathcal{D}$ according to the noise assumption of the underlying problem. The default case of Gaussian noise corresponds to the Frobenius norm on which we will focus on in this work. Further possible choices include the Kullback-Leibler divergence or more generalized divergence functions \cite{Cichocki:2009:BookNMF}. Moreover, NMF problems are usually ill-posed since their solution is non-unique \cite{Klingenberg:2009:NMFIllposedness}. Thus, they require stabilization techniques, which is typically done by adding regularization terms $\mathcal{R}_j$ into the NMF cost function $\mathcal{F}.$ However, besides the use case of regularization, the penalty terms $\mathcal{R}_j$ can be also used to enforce additional properties to the factorization matrices $U$ and $V.$ The general NMF minimization problem can therefore be written as
            \begin{equation} \label{eq:Background:NMFProblemGeneral}
                \min_{U, V\geq 0} \mathcal{D}(X, UV) + \sum_{j=1}^J \alpha_j \mathcal{R}_j(U, V) \eqqcolon \min_{U, V\geq 0} \mathcal{F}(U, V),
            \end{equation}
            where $\alpha_j$ are regularization parameters. Common choices for $\mathcal{R}_j$ are $\ell_1$ and $ \ell_2$ penalty terms. Further possible options are total variation regularization or other penalty terms which enforce orthogonality or even allow a supervised classification scheme in case the NMF is used as a prior feature extraction step \cite{Fernsel:2018:Survey,Leuschner:2018:NMFSupervised}. In this work, we will focus on the combination of orthogonality constraints and total variation penalty terms to construct an NMF model for spatially coherent clustering methods.
            
            On the other hand, $K$-means clustering is one of the most commonly used prototype-based, partitional clustering technique. As for any other clustering method, the main task is to partition a given set of objects into groups such that objects in the same group are more similar compared to the ones in other groups. These groups are usually referred to as \textit{clusters}. In mathematical terms, the aim is to partition the index set $\{ 1, 2,\dots,M\}$ of a corresponding given data set $\{ x_m\in \R^{N} \ \vert \ m = 1,\dots,M \}$ into disjoint sets $\mathcal{I}_k \subset \{1,\dots,M\},$ such that $\cup_{k=1,\dots,K} \mathcal{I}_k = \{1,\dots,M\}.$
            
            Many different variations and generalizations of $K$-means have been proposed and analyzed (see for instance \cite{Chen:2009:K-MeansCoresets,AggarwalReddy:2013:BookClustering} and the references therein), but we will focus in this section on the most common case.
            The method is based on two main ingredients. On the one hand, a similarity measure $\dist(\cdot, \cdot)$ is needed to specify the similarity between data points. The default choice is the squared Euclidean distance $\dist(x_i, x_j) \coloneqq \Vert x_i - x_j \Vert_2^2.$ On the other hand, so-called representative \textit{centroids} $c_k\in \R^N$ are computed for each cluster $\mathcal{I}_k.$ The computation of the clusters and centroids is based on the minimization of the within-cluster variances given by
            \begin{equation*}
                \mathcal{J} = \sum_{k=1}^K \sum_{m\in \mathcal{I}_k} \dist(x_m, c_k).
            \end{equation*}
            Due to the NP-hardness of the minimization problem \cite{MahajanEtal:2012:K-MeansNPHard}, heuristic approaches are commonly used to find an approximate solution. The $K$-means algorithm is the most common optimization technique which is based on an alternating minimization. After a suitable initialization, the first step is to assign each data point $x_m$ to the cluster with the closest centroid with respect to the distance measure $\dist(\cdot, \cdot).$ In the case of the squared Euclidean distance, the centroids are recalculated in a second step based on the mean of its newly assigned datapoints to minimize the sum of the within-cluster variances. Both steps are repeated until the assignments do not change anymore. 
            
            The relationship between the NMF and $K$-means clustering can be easily seen by adding proper further constraints to both problems. First of all, from the point of view of $K$-means, we assume nonnegativity of the given data and write the vectors $x_m$ row-wise to a data matrix $X$ such that $X = [x_1, \dots, x_M]^\intercal \in \R_{\geq 0}^{M\times N}.$ Furthermore, we define the so-called cluster membership matrix $\Tilde{U}\in \{0,1\}^{M\times K},$ such that
            \begin{equation*}
                \Tilde{U}_{mk} \coloneqq    \begin{cases}
                                                0 & \text{if} \ \ m\not\in \mathcal{I}_k \\
                                                1 & \text{if} \ \ m\in \mathcal{I}_k
                                            \end{cases}
            \end{equation*}
            and the centroid matrix $\Tilde{V} \coloneqq [c_1,\dots,c_K]^\intercal \in \R_{\geq 0}^{K\times N}.$ With this and by choosing the squared Euclidean distance function, it can be easily shown that the objective function $\mathcal{J}$ of $K$-means can be rewritten as $\mathcal{J}=\Vert X - \Tilde{U} \Tilde{V}\Vert_F^2,$ which has the same structure of the usual cost function of an NMF problem. However, the usual NMF does not constrain one of the matrices to have binary entries or, more importantly, to be row-orthogonal as it is the case for $\Tilde{U}.$ This ensures that each row of $\Tilde{U}$ contains only one nonzero element which gives the needed clustering interpretability.
            
            This gives rise to the problem of orthogonal nonnegative matrix factorization (ONMF), which is given by
            \begin{equation*}
                \min_{U, V\geq 0} \mathcal{F}(U, V), \ \ \text{s.t.} \ \ U^\intercal U = I,
            \end{equation*}
            where $I$ is the identity matrix. The matrices $U$ and $V$ of an ONMF problem will be henceforth also referred to as cluster membership matrix and centroid matrix respectively. In the case of the Frobenius norm as discrepancy term and without any regularization terms $\mathcal{R}_j,$ it can be shown that this problem is equivalent to weighted variant of the spherical $K$-means problem \cite{Pompili:2014:CompMethod}. For further variants of the relations between ONMF and $K$-means, we refer to the works \cite{Ding:2005:NMFKMeans,Ding:2006:triFac,LiDing:2006:Relationships} and the review articles \cite{LiDing:2014:NMFClusteringSurvey,Turkmen:2015:NMFClusteringSurvey}.
            
        \subsection{Algorithms for Orthogonal NMF} \label{subsec: Algorithms for Orthogonal NMF}
            Due to the ill-posedness of NMF problems and possible constraints on the matrices, tailored minimization approaches are needed. In this section, we review shortly common optimization techniques of NMF and ONMF problems which will also be used in this work for the derivation of algorithms for ONMF models including spatial coherence.
            
            For usual choices of $\mathcal{D}$ and $\mathcal{R}_j$ in the NMF problem \cref{eq:Background:NMFProblemGeneral}, the corresponding cost function $\mathcal{F}$ is convex in each of the variables $U$ and $V$ but non-convex in $(U,V).$ Therefore, the majority of optimization algorithms for NMF and ONMF problems are based on alternating minimization schemes
            \begin{align}
                U^{[i+1]} &= \arg\min_{U\geq 0} \mathcal{F}(U,V^{[i]}), \label{eq:AlgorithmsForONMF:AlternatingMin:U} \\
                V^{[i+1]} &= \arg\min_{V\geq 0} \mathcal{F}(U^{[i+1]},V). \label{eq:AlgorithmsForONMF:AlternatingMin:V}
            \end{align}
            One classical technique to tackle these minimization problems are alternating multiplicative algorithms, which only consist of summations and multiplications of matrices and therefore ensure the nonnegativity of $U$ and $V$ without any additional projection step provided that they are initialized appropriately. This approach was mainly popularized by the works of Lee and Seung \cite{Lee:1999:NMFMultiplicative,Lee:2001:NMFMultiplicative}, which also brought much attention to the NMF in general. The update rules are usually derived by analyzing the Karush-Kuhn-Tucker (KKT) first-order optimality conditions for each of the minimization problems in \cref{eq:AlgorithmsForONMF:AlternatingMin:U} and \cref{eq:AlgorithmsForONMF:AlternatingMin:V} or via the so-called majorize-minimization (MM) principle. The basic idea of the latter technique is to replace the NMF cost function $\mathcal{F}$ by a majorizing surrogate function $\mathcal{Q}_\mathcal{F}: \dom(\mathcal{F}) \times \dom(\mathcal{F}) \to \R, $ which is easier to minimize and whose tailored construction leads to the desired multiplicative updates rules defined by
            \begin{equation*}
                x^{[i+1]} \coloneqq \arg\min_{x\in \dom(\mathcal{F})} \mathcal{Q}_{\mathcal{F}} (x, x^{[i]}).
            \end{equation*}
            With the defining properties of a surrogate function that $\mathcal{Q}_\mathcal{F}$ majorizes $\mathcal{F}$ and $\mathcal{Q}_\mathcal{F}(x,x) = \mathcal{F}(x)$ for all $x\in \dom(\mathcal{F}),$ it can be easily shown that the above update rule leads to a monotone decrease of the cost function $\mathcal{F}.$ However, the whole method is based on an appropriate construction of the surrogate functions which is in general non-trivial. Possible techniques for common choices of $\mathcal{D}$ and $\mathcal{R}_j$ in the NMF cost function are based on the quadratic upper bound principle and Jensen's inequality \cite{Fernsel:2018:Survey}. Overall, multiplicative algorithms offer a flexible approach to various choices of NMF cost functions and will also be used in this work for some of the proposed and comparative methods.
            
            Another classical method are alternating (nonnegative) least squares (ALS) algorithms. They are based on the estimation of the stationary points of the cost function with a corresponding fixed point approach and a subsequent projection step to ensure the nonnegativity of the matrices. In case of the standard NMF problem with $\mathcal{D}$ being the Frobenius norm and without additional regularization terms, this leads for example to the update rules
            \begin{align}
                U^{[i+1]} &= [X{V^{[i]}}^\intercal(V^{[i]}{V^{[i]}}^\intercal)^{-1}]_+,\label{eq:AlgorithmsForONMF:ALS:U} \\
                V^{[i+1]} &= [({U^{[i+1]}}^\intercal U^{[i+1]})^{-1} {U^{[i+1]}}^\intercal X]_+,\label{eq:AlgorithmsForONMF:ALS:V}
            \end{align}
            where $[\cdot]_+$ denotes a suitable projection step to enforce the nonnegativity of all matrix entries.
            
            An extension to this approach is given by hierarchical alternating (nonnegative) least squares (HALS) algorithms, which solve nonnegative ALS problems column-wise for both matrices $U$ and $V$ \cite{Kimura:2015:CompMethod,LiEtal:2015:CompMethod} and will also be used as a comparative methodology.
            
            An optimization approach which was recently used for NMF problems is the widely known proximal alternating linearized minimization (PALM) \cite{WangEtal:2019:CompMethod} together with its extensions including stochastic gradients \cite{DriggsEtal:2020:SPRING}. As a first step, the cost function is splitted up into a differentiable part $\mathcal{F}_1$ and a non-differentiable part $\mathcal{F}_2.$ In its basic form, the PALM update rules concist of alternating gradient descent steps for $U$ and $V$ with learning rates based on the Lipschitz constants of the gradient of $\mathcal{F}_1$ in combination with a subsequent computation of a proximal operator of the function $\mathcal{F}_2.$ Some of these techniques will be used for the proposed methods in this work and will be discussed in more detail in \cref{sec: Orthogonal NMF with Spatial Coherence}.
            
            Further well-known techniques are for example projected gradient algorithms consisting of additive update rules, quasi-Newton approaches based on second-order derivatives of the cost function or algorithms based on an augmented Lagrangian concept \cite{Cichocki:2009:BookNMF,Pompili:2014:CompMethod}.
            
            All these methods can be also used to derive suitable algorithms for ONMF problems. Common approaches to include the orthogonality constraints are the use of Lagrangian multipliers \cite{Ding:2006:triFac,Choi:2008:CompMethod,YangOja:2010:CompMethod,Kimura:2015:CompMethod,LiEtal:2015:CompMethod} or the replacement of the hard constraint $U^\intercal U = I$ by adding a suitable penalty term into the NMF cost function to enforce approximate row-orthogonality for $U$ controlled by a regularization parameter \cite{LiWuPeng:2010:CompMethod,Mirzal:2014:CompMethod,ZhangEtal:2016:CompMethodNMFLPP,Fernsel:2018:Survey}. Other methods include optimization algorithms on the Stiefel manifold \cite{ZhangEtal:2016:CompMethodStiefelManifold}, the use of sparsity and nuclear norm minimization \cite{KimPark:2008:SparseNMF,PanNg:2018:CompMethod} or other techniques \cite{Pompili:2014:CompMethod,AsterisEtal:2015:CompMethod,WangEtal:2019:CompMethod}.
            
            In the next section, we will introduce the considered NMF models in this work and derive the corresponding optimization algorithms.
    
    \section{Orthogonal NMF with Spatial Coherence} \label{sec: Orthogonal NMF with Spatial Coherence}
        In this section, we present the spatially coherent clustering methods based on ONMF models together with the derivation of their optimization algorithms. Different from classical clustering approaches via ONMF, our proposed technique include the local information of a measured datapoint into the clustering process. This is done by including a TV regularization procedure on the cluster membership matrix $U,$ which naturally leads to spatially coherent regions while preserving their edges.
        
        This is, for instance, especially helpful for clusterings of hyperspectral data. If the neighbourhood of a measured spectrum $X_{m, \bullet}$ is associated to one cluster $\mathcal{I}_k,$ the inclusion of spatial coherence in the clustering model makes it more likely that $X_{m, \bullet}$ will be also classified to $\mathcal{I}_k.$ This spatial coherence can be also observed in many spectral imaging applications like Earth remote sensing or MALDI imaging, where locally close spectra have a higher probability to belong to the same class. 
        
        In the following, we divide our proposed techniques into so-called separated and combined methods.
        \subsection{Separated Methods} \label{subsec: Separated Methods}
            One straightforward approach to design a spatially coherent clustering method is to compute a clustering based on a classical ONMF model and subsequently perform a postprocessing on the obtained cluster memberhsip matrix based on total variation denoising. The general workflow is provided in \cref{alg:ONMF-TV} and will be henceforth referred to as \texttt{ONMF-TV}.
            \begin{algorithm}
                \setcounter{ALC@unique}{0}
                \caption{\texttt{ONMF-TV}}
                \label{alg:ONMF-TV}
                \begin{algorithmic}[1] 
                    \STATE\label{alg:ONMF-TV_Line1}{\textbf{Input:} $\displaystyle \ X\in \R_{\geq 0}^{M\times N},\ \tau>0, \ \theta\in \R^J$}
                    \STATE\label{alg:ONMF-TV_Line2}{\textbf{Initialize:} $\displaystyle \ U^{[0]}\in \R_{\geq 0}^{M\times K},\ V^{[0]}\in \R_{\geq 0}^{K\times N}$}
                    \STATE\label{alg:ONMF-TV_Line3}{$\displaystyle (U, V) \gets \textsc{ONMF}_{\theta}(X, U^{[0]}, V^{[0]})$}
                    \STATE\label{alg:ONMF-TV_Line4}{$\displaystyle U \gets \textsc{TVDenoiser}_{\tau}(U)$}
                \end{algorithmic}
            \end{algorithm}
            
            After a suitable initialization of the cluster membership matrix and the centroid matrix in \cref{alg:ONMF-TV_Line2}, a classical ONMF model is used in \cref{alg:ONMF-TV_Line3} to perform a clustering of the given data $X,$ where $\theta$ are possible hyperparameters which typically have to be chosen a priori. These ONMF models will be chosen based on some of the works described in \cref{subsec: Algorithms for Orthogonal NMF} and will be specified in more detail in \cref{subsec: Choice of Separated Methods}. Afterwards, a TV denoising algorithm is applied on the cluster membership matrix $U$ obtained by $\textsc{ONMF}_{\theta}$ to induce the spatial coherence in the clustering. In this work and for all considered separated methods, the denoising step is based on evaluating the proximal mapping column-wise on $U,$ which is defined by the minimization problem
             \begin{equation} \label{eq:proximalMapping}
                \prox_{\tau \Vert \cdot \Vert_{\TV}}(x) \coloneqq \underset{y\in \R^M}{\argmin} \left\{ \frac{1}{2} \Vert y - x\Vert^2 + \tau \Vert y\Vert_{\TV} \right\}.
            \end{equation}
            Here, $\tau>0$ is a regularization parameter and $\Vert \cdot \Vert_{\TV}$ denotes the classical isotropic TV \cite{ChambolleEtal:2010:TVIntro,BeckTeboulle:2009:TVFastGradient} and corresponds to the definition in \cref{eq:TV:DiscreteEpsilon:Definition} with $\varepsilon_{\TV}=0$ if $\Vert \cdot \Vert_{\TV}$ is applied on a matrix $U.$ Afterwards, every negative entry of $U$ is projected to zero to ensure the corresponding nonnegativity constraint.
            
            We consider this workflow as baseline for our comparison with the combined methods presented in the following \cref{subsec: Combined Methods}.
            
        
        \subsection{Combined Methods} \label{subsec: Combined Methods} 
            In this section, we present the so-called combined methods together with different optimization algorithms for their numerical solution. Different from the separated methods in \cref{subsec: Separated Methods}, this coupled approach includes a total variation penalty term into the ONMF model to induce the spatial coherence in the clustering. The combination of orthogonality constraints and a total variation penalty term in the NMF problem leads in general to a higher effort to derive suitable solution algorithms. However, the main motivation is that this joint workflow allows the clustering process to take advantage of the TV regularization leading to better local minima, which could therefore enhance the quality of the clustering compared to classical approaches or the previous described separated methods, where the spatial coherence is just enforced in an independent subsequent TV denoising step.
            
            In the following, we will present different multiplicative update rules and algorithms based on proximal gradient descent approaches. Furthermore, we analyzed a combined approach based on an alternating least squares procedure, which can be found in \cite[Section 4.3]{Cichocki:2009:BookNMF}. However, first numerical tests showed that this approach does not deliver satisfactory clustering results. Hence, this method will be omitted in the later numerical evaluation of this work.
            
            The choice of the stopping criteria and of the hyperparameters for all considered approaches will be described in \cref{sec: Numerical Experiments} of the numerical experiments of this work.

            
            
            \subsubsection{Multiplicative Update Rules} \label{subsubsec: Multiplicative Update Rules}
                \paragraph{\texttt{ONMFTV-MUL1}}\label{par:ONMFTV-MUL1}
                    Our first multiplicative algorithm is taken from the work \cite{Fernsel:2018:Survey} without any modification and is based on the ONMF model
                    \begin{equation}\label{eq:CombinedMethods:MUL:ONMFTV-MUL1:NMFProblem1}
                        \min_{U, V\geq 0} \frac{1}{2} \Vert X - UV\Vert_F^2 + \frac{\sigma}{2} \Vert I - U^\intercal U\Vert_F^2 + \frac{\tau}{2} \Vert U \Vert_{\TVeps},
                    \end{equation}
                    where $\sigma, \tau\geq 0$ are regularization parameters and $\Vert \cdot \Vert_{\TVeps}$ is the smoothed, isotropic total variation \cite{Defrise:2011:TV,Fernsel:2018:Survey} defined by
                    \begin{equation} \label{eq:TV:DiscreteEpsilon:Definition}
                	    \Vert U \Vert_{\TVeps} \coloneqq \sum_{k=1}^K \sum_{m=1}^M \vert \nabla_{mk}U\vert \coloneqq \sum_{k=1}^K \sum_{m=1}^M \sqrt{\varepsilon_{\TV}^2 + \sum_{\tilde{m} \in \mathcal{N}_m} (U_{mk}-U_{\tilde{m} k})^2}.
                	\end{equation}
                    Furthermore, $\varepsilon_{\TV}>0$ is a positive, small, predefined constant to ensure the differentiability of the TV penalty term, which is needed due to the MM principle for the optimization approach. Furthermore, $\mathcal{N}_m$ are index sets referring to spatially neighboring pixels. The default case in two dimensions for the neighbourhood of a non-boundary pixel in $(i,j)$ is $\mathcal{N}_{(i,j)} = \{ (i+1, j), (i, j+1)\}$ to obtain an estimate of the gradient components in both directions.
                    
                    The derivation of the solution algorithm is described in \cite{Fernsel:2018:Survey} and is based on the MM principle mentioned in \cref{subsec: Algorithms for Orthogonal NMF}. The surrogate function $\mathcal{Q}_\mathcal{F}(x,a)$ of such approaches majorizes $\mathcal{F}$ and is typically quadratic in $x.$ In order to avoid complications in constructing a suitable surrogate function for the cost function in \cref{eq:CombinedMethods:MUL:ONMFTV-MUL1:NMFProblem1}, the fourth order terms from $\Vert I - U^\intercal U\Vert_F^2$ have to be avoided. In \cite{Fernsel:2018:Survey}, this problem is solved by introducing an auxiliary variable $W\in \R_{\geq 0}^{M\times K}$ and by reformulating the minimization problem in \cref{eq:CombinedMethods:MUL:ONMFTV-MUL1:NMFProblem1} as
                    \begin{equation}\label{eq:CombinedMethods:MUL:ONMFTV-MUL1:NMFProblem2}
                        \min_{U, V, W\geq 0} \underbrace{\frac{1}{2} \Vert X - UV\Vert_F^2 + \frac{\sigma_1}{2} \Vert I - W^\intercal U\Vert_F^2 + \frac{\sigma_2}{2} \Vert W - U\Vert_F^2 + \frac{\tau}{2} \Vert U \Vert_{\TVeps}}_{\eqqcolon \mathcal{F}(U,V,W)}.
                    \end{equation}
                    For this problem, a suitable surrogate function and a multiplicative algorithm can be derived. The details of the derivation will not be discussed here and can be found in \cite{Fernsel:2018:Survey} in all details. We also provide a short outline of the derivation in \cref{app:sec:Derivation of the Algorithm for ONMFTV-MUL1}. However, we describe the final update rules obtained by the MM principle in \cref{alg:ONMFTV-MUL1} and define for this purpose the matrices $ P(U), Z(U)\in \mathbb{R}_{\geq 0}^{M\times K} $ as
                    \begin{align}
                    	P(U)_{m k} &\coloneqq  \dfrac{\sum_{\tilde{m} \in \mathcal{N}_m} 1}{\vert \nabla_{mk} U \vert} + \sum_{\tilde{m} \in \bar{\mathcal{N}}_m} \dfrac{1}{\vert \nabla_{\tilde{m} k} U \vert}, \label{eq:ONMFTV-MUL1-P}\\
                		Z(U)_{m k} &\coloneqq  \dfrac{1}{P(U)_{m k}} \left ( \dfrac{1}{\vert \nabla_{mk} U \vert} \sum_{\tilde{m} \in {\mathcal{N}}_m} \dfrac{U_{m k} + U_{\tilde{m}  k}}{2} + \sum_{\tilde{m} \in \bar{\mathcal{N}}_m} \dfrac{U_{m k} + U_{\tilde{m}  k}}{2 \vert \nabla_{\tilde{m} k} U \vert} \right ), \label{eq:ONMFTV-MUL1-Z}
                    \end{align}
                    where $\bar{\mathcal{N}}_m$ is the so-called adjoint neighborhood given by $\tilde{m}\in \bar{\mathcal{N}}_m \Leftrightarrow m\in \mathcal{N}_{\tilde{m}}.$

    		        \begin{algorithm}
    		            \setcounter{ALC@unique}{0}
                        \caption{\texttt{ONMFTV-MUL1}}
                        \label{alg:ONMFTV-MUL1}
                        \begin{algorithmic}[1] 
                            \STATE{\textbf{Input} $\ X\in \R_{\geq 0}^{M\times N},\ \ \sigma_1, \sigma_2, \tau>0,\ \ i=0 $}
                            \STATE{\textbf{Initialize} $\ U^{[0]}, W^{[0]}\in \R_{> 0}^{M\times K},\ \ V^{[0]}\in \R_{> 0}^{K\times N}$}
                            \REPEAT
                            \STATE\label{alg:ONMFTV-MUL1-Line4}{$ \displaystyle U^{[i+1]} = \left[U^{[i]} \circ \left( \frac{X{V^{[i]}}^\intercal + \tau P(U^{[i]}) \circ Z(U^{[i]}) + (\sigma_1 + \sigma_2)W^{[i]}}{\tau \circ P(U^{[i]}) \circ U^{[i]} + \sigma_2U^{[i]} + U^{[i]}V^{[i]}{V^{[i]}}^\intercal + \sigma_1 W^{[i]}{W^{[i]}}^\intercal U^{[i]} } \right) \right]_{>0}$}
                            \STATE{$ \displaystyle V^{[i+1]} = \left[V^{[i]} \circ \left( \frac{{U^{[i+1]}}^\intercal X }{{U^{[i+1]}}^\intercal U^{[i+1]} V^{[i]} } \right)\right]_{>0}$}
                            \STATE{$ \displaystyle W^{[i+1]} = \left[W^{[i]} \circ \left( \frac{(\sigma_1 + \sigma_2) U^{[i+1]}}{\sigma_1 U^{[i+1]} {U^{[i+1]}}^\intercal W^{[i]} + \sigma_2 W^{[i]} } \right)\right]_{>0}$}
                            \STATE{$ \displaystyle i \leftarrow i + 1 $ }
                            \UNTIL{\textit{Stopping criterion satisfied}}
                        \end{algorithmic}
                    \end{algorithm}
                    We denote by $\circ$ as well as the fraction line the element-wise (Hadamard) multiplication and division respectively. Due to the multiplicative structure of the update rules, the nonnegativity of the iterates is preserved. However, a strict positive initialization of $U, V$ and $W$ is needed to avoid numerical instabilities and the zero locking phenomenon caused by zero entries in the matrices, which is characteristic for all algorithms based on multiplicative update rules (see e.g.\ \cite{Mirzal:2014:CompMethod}). For the same reason, we perform a subsequent element-wise projection step for every matrix defined by $[\lambda]_{>0} \coloneqq \max\{\lambda, \varepsilon_{P_1}\}$ with $\varepsilon_{P_1} = 1\cdot 10^{-16}.$ Analogously, a projection step is applied for too large entries with $\varepsilon_{P_2} = 1\cdot 10^{35}$ being the corresponding parameter. Finally, the above algorithm ensures a monotone decrease of the cost function in \cref{eq:CombinedMethods:MUL:ONMFTV-MUL1:NMFProblem2} due to its construction based on the MM principle \cite{Fernsel:2018:Survey}.
                    \begin{theorem}[\texttt{ONMFTV-MUL1}] \label{thm:ONMFTV-MUL1}
                        The \cref{alg:ONMFTV-MUL1} ensures a monotone decrease of the cost function defined by the NMF model in \cref{eq:CombinedMethods:MUL:ONMFTV-MUL1:NMFProblem2}.
                    \end{theorem}
                    
                \paragraph{\texttt{ONMFTV-MUL2}}\label{par:ONMFTV-MUL2}
                    In this section, we derive another multiplicative algorithm following the ideas in the work \cite{LengEtal:2017:PropMethod} based on a continuous formulation of an isotropic and differentiable version of the TV penalty term given by
                    \begin{equation}\label{eq:ONMFTV-MUL2-TV Penalty Term}
                        \TVeps(u) \coloneqq \int_\Omega \Vert \nabla u\Vert_{\varepsilon_{\TV}} \diff (x_1, x_2)
                    \end{equation}
                    for $\Omega\subset \R^2,$ a sufficiently smooth $u: \Omega \to \R$ with bounded variation (see e.g.\ \cite{ChambolleEtal:2010:TVIntro}) and for
                    \begin{equation}\label{eq:ONMFTV-MUL2-norm of nabla u}
                        \Vert \nabla u\Vert_{\varepsilon_{\TV}} \coloneqq \sqrt{ \left( \frac{\partial u}{\partial x_1} \right)^2 + \left( \frac{\partial u}{\partial x_2} \right)^2 + \varepsilon_{\TV}^2 }
                    \end{equation}
                    with a small $\varepsilon_{\TV} > 0.$ The application of the TV regularization on the discrete matrix $U$ is done via a subsequent discretization step, which is specified in more detail in \cref{app:sec:Details on the Algorithm ONMFTV-MUL2}. Thus, we consider in this section the orthogonal NMF model
                    \begin{equation}\label{eq:CombinedMethods:MUL:ONMFTV-MUL2:NMFProblem}
                        \min_{U, V\geq 0} \underbrace{\frac{1}{2} \Vert X - UV\Vert_F^2 + \frac{\sigma_1}{4} \Vert I - U^\intercal U\Vert_F^2 + \tau \TVeps(U)}_{\eqqcolon \mathcal{F}(U,V)},
                    \end{equation}
                    with the sloppy notation of $\TVeps(U)$ for the matrix $U,$ where the discretization step is implicitly included. Different from the model \texttt{ONMFTV-MUL1} in the previous section, we do not include any auxiliary variable $W.$
                    
                    The update rule for $U$ is based on a classical gradient descent approach
                    \begin{equation}\label{eq:ONMFTV-MUL2-GradientDescentStep}
                        U^{[i+1]} \coloneqq U^{[i]} - \Gamma^{[i]} \circ \nabla_U \mathcal{F}(U^{[i]},V^{[i]}),
                    \end{equation}
                    with a step size $\Gamma^{[i]} \in \R_{\geq 0}^{M\times K}.$ By computing the gradient $\nabla_U \mathcal{F}(U^{[i]},V^{[i]})$ and choosing an appropriate step size $\Gamma^{[i]},$ this leads to the multiplicative update rule shown in \cref{alg:ONMFTV-MUL2}. For the minimization with respect to $V,$ we simply choose the multiplicative update rule given in \cref{alg:ONMFTV-MUL1}.
                    \begin{algorithm}
                        \setcounter{ALC@unique}{0}
                        \caption{\texttt{ONMFTV-MUL2}}
                        \label{alg:ONMFTV-MUL2}
                        \begin{algorithmic}[1] 
                            \STATE{\textbf{Input} $\ X\in \R_{\geq 0}^{M\times N},\ \ \sigma_1, \tau>0,\ \ i=0 $}
                            \STATE{\textbf{Initialize} $\ U^{[0]} \in \R_{> 0}^{M\times K},\ \ V^{[0]}\in \R_{> 0}^{K\times N}$}
                            \REPEAT
                            \STATE\label{alg:ONMFTV-MUL2_Line4}{$ \displaystyle U^{[i+1]} = \left[ U^{[i]} \circ \frac{X{V^{[i]}}^\intercal + \tau \divergence\left( \frac{\nabla U^{[i]}}{ \Vert \nabla U^{[i]} \Vert_{\varepsilon_{\TV}} } \right) + \sigma_1 U^{[i]} } {U^{[i]} V^{[i]} {V^{[i]}}^\intercal + \sigma_1 U^{[i]} {U^{[i]}}^\intercal U^{[i]} }\right]_{>0} $}
                            \STATE{$ \displaystyle V^{[i+1]} = \left[V^{[i]} \circ \left( \frac{{U^{[i+1]}}^\intercal X }{{U^{[i+1]}}^\intercal U^{[i+1]} V^{[i]} } \right)\right]_{>0}$}
                            \STATE{$ \displaystyle i \leftarrow i + 1 $ }
                            \UNTIL{\textit{Stopping criterion satisfied}}
                        \end{algorithmic}
                    \end{algorithm}
                    The term $\divergence\left( \nabla U^{[i]} / \Vert \nabla U^{[i]} \Vert_{\varepsilon_{\TV}} \right)$ is again an abuse of notation and implicitly includes a discretization step (see \cref{app:subsec:Discretization of the TV Gradient}). It can be seen as the gradient of the TV penalty term and is obtained by analyzing the corresponding Euler-Lagrange equation (see \cref{app:subsec:Derivation of the Update Rules}). Note that in discretized form, it is a matrix of size $M\times K$ and can also contain negative entries. Hence, this update step is, strictly speaking, not a multiplicative update, since it cannot enforce the nonnegativity of the matrix $U$ by itself. However, as in \cref{alg:ONMFTV-MUL1}, a projection step given by $[\cdot]_{>0}$ is applied subsequently for both matrices to avoid numerical instabilities and to ensure the nonnegativity of $U.$ Different from the approach in \texttt{ONMFTV-MUL1}, a monotone decrease of the cost function cannot be guaranteed based on a classical surrogate functional since $\mathcal{F}$ also contains fourth-order terms due to the penalty term $\Vert I - U^\intercal U\Vert_F^2.$
                    
                    Furthermore, note the similarity of the update rules given in \cite{Mirzal:2014:CompMethod}, where no total variation penalty term is considered. For more details on the derivation of this algorithm and the discretization of the divergence term, we refer the reader to \cref{app:sec:Details on the Algorithm ONMFTV-MUL2}.
            \subsubsection{Proximal Alternating Linearized Minimization} \label{subsubsec: Promximal Alternating Linearized Minimization}
                Following the optimization procedure of the very recent work \cite{DriggsEtal:2020:SPRING}, we consider in this section several proximal alternating linearized minimization (PALM) schemes. 
                For all presented methods in this section, we consider the NMF model
                \begin{equation}\label{eq:CombinedMethods:PALM:NMFModel}
                    \min_{U, V, W\geq 0} \underbrace{\frac{1}{2} \Vert X - UV\Vert_F^2 + \frac{\sigma_1}{2} \Vert I - W^\intercal U\Vert_F^2 + \frac{\sigma_2}{2} \Vert W - U\Vert_F^2}_{\eqqcolon \mathcal{F}(U,V,W)} + \tau \underbrace{\Vert U \Vert_{\TV}}_{\eqqcolon \mathcal{J}(U)},
                \end{equation}
                where $\Vert \cdot \Vert_{\TV}$ is defined according to \cref{eq:proximalMapping}. The additional auxiliary variable $W$ is needed to ensure the Lipschitz continuity of the partial gradients of $\mathcal{F}$ and hence the convergence rates of the respective algorithms showed in \cite{DriggsEtal:2020:SPRING}.
                
                The general scheme of these algorithms are based on an alternating minimization procedure with a gradient descent step with respect to the differentiable function $\mathcal{F}(U,V,W)$ via the computation of full gradients or gradient estimates and a subsequent column-wise application of the proximal operator for the non-differentiable part $\mathcal{J}$ given by the update step
                \begin{equation}\label{eq:CombinedMethods:PALM:GeneralUpdateForUWithProx}
                    U^{[i+1]}_{\bullet, k} \coloneqq \prox_{\tau \eta \mathcal{J}}\left( U^{[i]}_{\bullet, k} - \eta (\tilde{\nabla}_U \mathcal{F}(U^{[i]},V,W))_{\bullet, k} \right)
                \end{equation}
                for a suitable step size $\eta>0.$ Furthermore, $\tilde{\nabla}_U \mathcal{F}$ is either the partial derivative ${\nabla}_U \mathcal{F}$ of $\mathcal{F}$ with respect to $U$ or some random gradient estimate. The nonnegativity constraint of the matrices is ensured by a final projection of all negative values to zero. 
                
                As usual for PALM algorithms, adaptive step sizes based on the local Lipschitz constants of the partial gradients of $\mathcal{F}$ are used, which will be approximated via the power iteration for all considered approaches. More information on the computation of these estimates and the choice of the step sizes are given in the subsequent descriptions of the specific algorithms and in \cref{sec: Numerical Experiments} of the numerical experiments as well as in \cref{app:sec:Algorithmic Details on the Proximal Gradient Descent Approach}. Details on the derivation of the gradients and the Lipschitz constants for the specific NMF model in \cref{eq:CombinedMethods:PALM:NMFModel} are given in \cref{app:sec:Algorithmic Details on the Proximal Gradient Descent Approach}.
                
                In the following, we will write in short according to the update rule in \cref{eq:CombinedMethods:PALM:GeneralUpdateForUWithProx}
                \begin{equation*}
                    U^{[i+1]} \coloneqq \prox_{\tau \eta \mathcal{J}}\left( U^{[i]} - \eta \tilde{\nabla}_U \mathcal{F}(U^{[i]},V,W) \right)
                \end{equation*}
                 for the whole matrix $U.$ For convergence results on the PALM algorithms, we refer the reader to \cite{DriggsEtal:2020:SPRING}.
                \paragraph{\texttt{ONMFTV-PALM}} \label{par:ONMFTV-PALM}
                    The proximal alternating linearized minimization is based on the general algorithm scheme dscribed above. The main steps are illustrated in \cref{alg:ONMFTV-PALM}.
                    \begin{algorithm}
                        \setcounter{ALC@unique}{0}
                        \caption{\texttt{ONMFTV-PALM}}
                        \label{alg:ONMFTV-PALM}
                        \begin{algorithmic}[1] 
                            \STATE{\textbf{Input} $\ X\in \R_{\geq 0}^{M\times N},\ \ \sigma_1, \sigma_2, \tau>0,\ \ i=0 $}
                            \STATE{\textbf{Initialize} $\ U^{[0]}, W^{[0]}\in \R_{\geq 0}^{M\times K},\ \ V^{[0]}\in \R_{\geq 0}^{K\times N}$}
                            \REPEAT
                            \STATE{$ \displaystyle \eta_{U^{[i]}} = \textsc{PowerIt}_U(V^{[i]}, W^{[i]})$}
                            \STATE{$ \displaystyle U^{[i+1]} = \left[ \prox_{\tau \eta_{U^{[i]}} \mathcal{J}}\left( U^{[i]} - \eta_{U^{[i]}} \nabla_U \mathcal{F}(U^{[i]},V^{[i]},W^{[i]}) \right) \right]_{\geq 0} $ }
                            \STATE{$ \displaystyle \eta_{V^{[i]}} = \textsc{PowerIt}_V(U^{[i+1]})$}
                            \STATE{$ \displaystyle V^{[i+1]} = \left[ V^{[i]} - \eta_{V^{[i]}} \nabla_V \mathcal{F}(U^{[i+1]},V^{[i]},W^{[i]})  \right]_{\geq 0} $ }
                            \STATE{$ \displaystyle \eta_{W^{[i]}} = \textsc{PowerIt}_W(U^{[i+1]})$}
                            \STATE{$ \displaystyle W^{[i+1]} = \left[ W^{[i]} - \eta_{W^{[i]}} \nabla_W \mathcal{F}(U^{[i+1]},V^{[i+1]},W^{[i]})  \right]_{\geq 0} $ }
                            \STATE{$ \displaystyle i \leftarrow i + 1 $ }
                            \UNTIL{\textit{Stopping criterion satisfied}}
                        \end{algorithmic}
                    \end{algorithm}
                    For the gradient descent step, the full classical partial derivatives $\nabla_U \mathcal{F}, \nabla_V \mathcal{F}$ and $\nabla_W \mathcal{F}$ are computed without the consideration of any batch sizes. Furthermore, the function \textsc{PowerIt} denotes the power method to compute the needed step sizes and $[\cdot ]_{\geq 0}$ is the componentwise projection of negative matrix entries to zero. The algorithm terminates until a suitable stopping criterion is satisfied, which will be further specified in \cref{sec: Numerical Experiments} and \cref{app:sec:Parameter Choice}. This approach will be henceforth referred to as \texttt{ONMFTV-PALM}.
                \paragraph{\texttt{ONMFTV-iPALM}} \label{par:ONMFTV-iPALM}
                    A slightly extended version of \texttt{ONMFTV-PALM} is the so-called inertial proximal alternating linearized minimization (iPALM) algorithm, which introduces an additional momentum term into \cref{alg:ONMFTV-PALM} and hence improves the convergence rate. Since the iPALM algorithm still follows the rough outline of \texttt{ONMFTV-PALM} and uses the classical partial gradients of $\mathcal{F}$, we will not present the whole algorithm at this point and refer the reader to the corresponding work in \cite{PockSabach:2016:iPALM}.
                    
                    This method will be referred to as \texttt{ONMFTV-iPALM}.
                \paragraph{\texttt{ONMFTV-SPRING}} \label{par:ONMFTV-SPRING}
                    The so-called stochastic proximal alternating linearized minimization (SPRING) is an extended version of the PALM method, where the full gradients are replaced by random estimates.
                    
                    One basic assumption of this approach is that $\mathcal{F}$ is separable depending on the considered variable for which the cost function is minimized. In the case of the minimization with respect to $U,$ the function $\mathcal{F}$ can be expressed as
                    \begin{equation*}
                        \mathcal{F}(U, V, W) = \sum_{n=1}^N \frac{1}{2} \Vert X_{\bullet, n} - UV_{\bullet, n} \Vert_2^2 + \underbrace{\sum_{k=1}^K \frac{\sigma_1}{2} \Vert I_{\bullet, k} - U^\intercal W_{\bullet, k}\Vert_2^2 + \frac{\sigma_2}{2} \Vert U_{\bullet, k} - W_{\bullet, k}\Vert_2^2}_{\eqqcolon \tilde{F}(U, W)}.
                    \end{equation*}
                    However, since $K$ is usually chosen to be small, we use instead the formulation
                    \begin{equation*}
                        \mathcal{F}(U, V, W) = \sum_{n=1}^N \left[ \frac{1}{2} \Vert X_{\bullet, n} - UV_{\bullet, n} \Vert_2^2 + \frac{1}{N} \tilde{\mathcal{F}}(U, W) \right] \eqqcolon  \sum_{n=1}^N \mathcal{F}_n(U, V, W)
                    \end{equation*}
                    to compute the estimates of the gradients based on the functions $\mathcal{F}_n.$ These random estimates are formed by using just a few indices of $\mathcal{F}_n,$ which are the elements of the so-called mini-batch $n\in \mathcal{B}_{i,j}^U \subset \{1,\dots, N\},$ where $i$ denotes the iteration number of the SPRING algorithm and $j\in \{ 1,\dots, s_r \}$ specifying the currently used subsample of indices with $1/s_r$ being the subsample ratio. One classical example of a gradient estimator, which is also used in this work, is the stochastic gradient descent (SGD) estimator given by
                    \begin{equation}\label{eq:ONMFTV-SPRING-SGDEstimator}
                        \tilde{\nabla}_U^{i,j} \mathcal{F}(U,V,W) \coloneqq \sum_{n\in \mathcal{B}_{i,j}^U} \nabla_U \mathcal{F}_n(U,V,W).
                    \end{equation}
                    For other possible gradient estimators like the SAGA or SARAH estimator, we refer the reader to the work \cite{DriggsEtal:2020:SPRING} and the references therein.
                    
                    The case of the minimization with respect to $V$ and $W$ is more straightforward. For $V,$ the penalty terms for the orthogonality can be dropped and we use the separability of the function with respect to the rows of $X$ and $U.$ The computation of the SGD estimator is done analogously. For minimizing with respect to $W,$ the expression based on $\tilde{\mathcal{F}}(U, W)$ can be used by omitting the data fidelity term. However, since $K\ll \min\{M,N\},$ we use in this work the full gradient of $\mathcal{F}$ for the minimization with respect to $W.$
                    
                    The main steps of the algorithm are presented in \cref{alg:ONMFTV-SPRING}, which will be referred to as \texttt{ONMFTV-SPRING}. For details on the choice of the hyperparameters and the computation of the step sizes, we refer the reader to \cref{sec: Numerical Experiments} and \cref{app:sec:Algorithmic Details on the Proximal Gradient Descent Approach}.
                    \begin{algorithm}
                        \setcounter{ALC@unique}{0}
                        \caption{\texttt{ONMFTV-SPRING}}
                        \label{alg:ONMFTV-SPRING}
                        \begin{algorithmic}[1] 
                            \STATE{\textbf{Input} $\ X\in \R_{\geq 0}^{M\times N},\ \ \sigma_1, \sigma_2, \tau>0,\ \ s_r\in \mathbb{N},\ \ i=0 $}
                            \STATE{\textbf{Initialize} $\ U^{[0,1]}, W^{[0,1]}\in \R_{\geq 0}^{M\times K},\ \ V^{[0,1]}\in \R_{\geq 0}^{K\times N}$}
                            \REPEAT
                                \FOR{$j=1,\dots, s_r$}
                                    \STATE{$ \displaystyle \eta_{U^{[i,j]}} = \textsc{PowerIt}_U(V^{[i,j]}, W^{[i,j]})$}
                                    \STATE\label{alg:ONMFTV-SPRING_Line6}{$ \displaystyle U^{[i,j+1]} = \left[ \prox_{\tau \eta_{U^{[i,j]}} \mathcal{J}}\left( U^{[i,j]} - \eta_{U^{[i,j]}} \tilde{\nabla}_U^{i,j} \mathcal{F}(U^{[i,j]},V^{[i,j]},W^{[i,j]}) \right) \right]_{\geq 0} $ }
                                    \STATE{$ \displaystyle \eta_{V^{[i,j]}} = \textsc{PowerIt}_V(U^{[i,j+1]})$}
                                    \STATE{$ \displaystyle V^{[i,j+1]} = \left[ V^{[i,j]} - \eta_{V^{[i,j]}} \tilde{\nabla}_V^{i,j} \mathcal{F}(U^{[i,j+1]},V^{[i.j]},W^{[i.j]})  \right]_{\geq 0} $ }
                                    \STATE{$ \displaystyle \eta_{W^{[i.j]}} = \textsc{PowerIt}_W(U^{[i,j+1]})$}
                                    \STATE{$ \displaystyle W^{[i,j+1]} = \left[ W^{[i,j]} - \eta_{W^{[i,j]}} \nabla_W \mathcal{F}(U^{[i,j+1]},V^{[i,j+1]},W^{[i,j]})  \right]_{\geq 0} $ }
                                \ENDFOR
                                \STATE{$ \displaystyle U^{[i+1,1]} = U^{[i,s_r+1]},\ \ \ V^{[i+1,1]} = V^{[i,s_r+1]},\ \ \ W^{[i+1,1]} = W^{[i,s_r+1]}$}
                                \STATE{$ \displaystyle i \leftarrow i + 1 $ }
                            \UNTIL{\textit{Stopping criterion satisfied}}
                        \end{algorithmic}
                    \end{algorithm}

    \section{Numerical Experiments} \label{sec: Numerical Experiments}
        
        \subsection{Dataset} \label{subsec: Dataset}
            Concerning the numerical experiments of this work, we consider a hyperspectral dataset obtained from a matrix-assisted laser desorption/ionization (MALDI) imaging measurement of a human colon tissue sample \cite{AlexandrovEtal:2011:MALDI}. This technique is a mass spectrometry imaging (MALDI-MSI) method, which is able to provide a spatial molecular profile of a given analyte. Together with the technological advancements in acquisition speed and robustness over the last decade, MALDI-MSI has become a standard tool in proteomics and applications in medical research, such as characterization of tumor subtypes and extraction of characteristic spectra \cite{Leuschner:2018:NMFSupervised}, have become feasible.
            
            In general, a measurement with a mass spectrometer can be subdivided into three main steps: the ionization of the sample, followed by the separation and detection of the ions.
            Considering MALDI imaging, the ionization part is characterized by the term MALDI. Different from other ionization techniques like the so-called secondary-ion mass spectrometry (SIMS) or desorption electrospray ionization (DESI), MALDI allows to measure a wider mass range of the extracted ions. This is possible due to an application of a matrix onto the tissue sample to optimize the transfer of the ionization energy and to extract the molecules out of the analyte. The needed ionization energy is provided by a laser, which is shot on the tissue sample by following a grid pattern.
            The separation of the molecules is done by a mass analyzer. One typical method, which was also used for the dataset in this work, is to accelerate the ions by an electric field and to measure the time-of-flight (TOF) of the particles. The final step of the whole measurement is the detection of the ions with the help of an ion detector.
            
            Another major aspect of the whole MALDI-MSI workflow is the preparation of the analytes before the actual measurement with a MALDI-TOF device, which follows standardized protocols. For more information, we refer the reader to \cite{AlexandrovEtal:2011:MALDI,AichlerWalch:2015:MALDIInfo}.
            \begin{figure}[tbhp]
            	\centering
            	\includegraphics[width=0.7\textwidth]{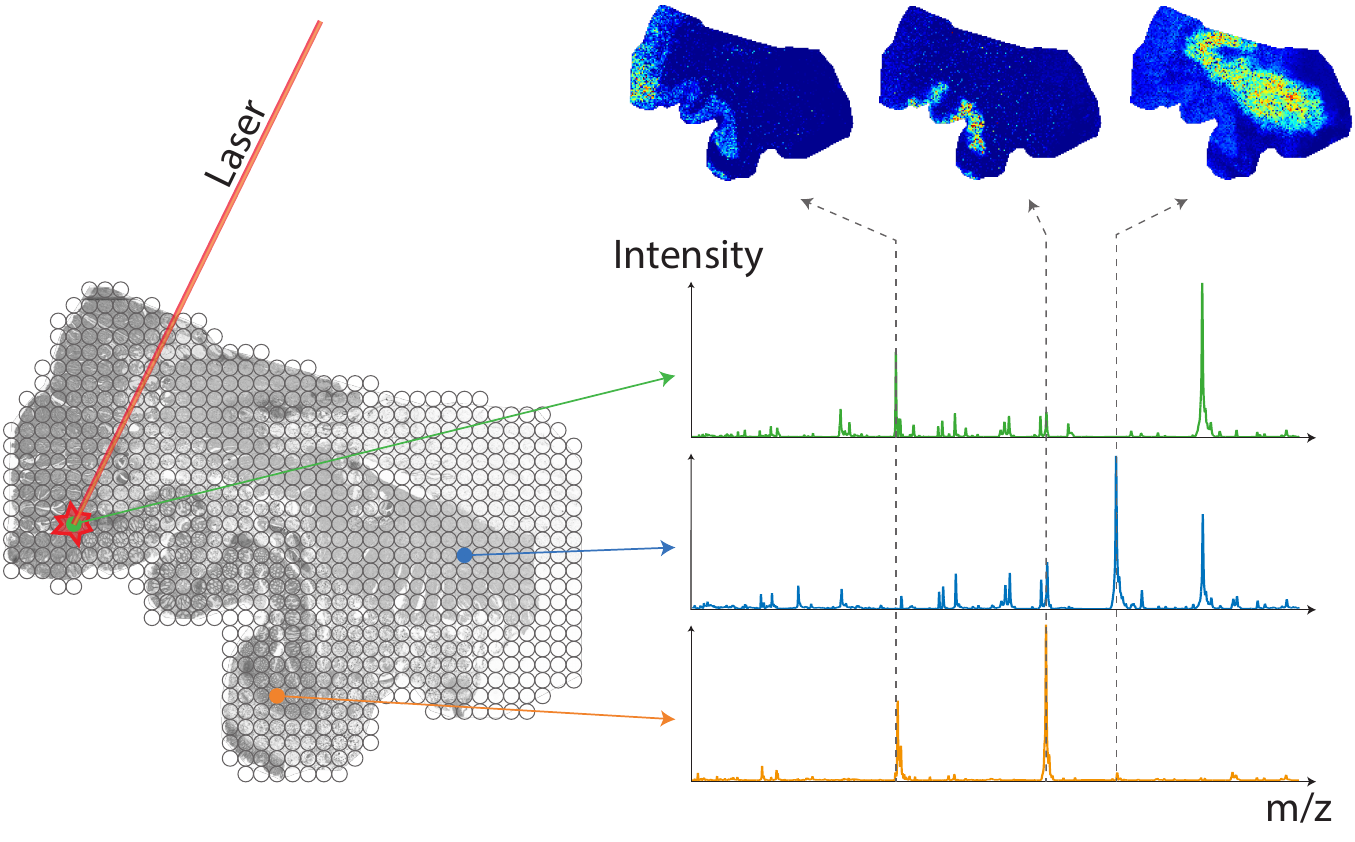}
            	\caption{Structure of a typical MALDI Dataset (slightly changed from the original image in \cite[p.1010]{Fernsel:2018:Survey}).}
            	\label{fig:dataStructure:MALDIColon}
            \end{figure}
            
            The structure of a typical MALDI dataset is shown in \cref{fig:dataStructure:MALDIColon}. For every point on the tissue slide, which is shot by the laser, a whole mass spectrum is acquired leading to a hyperspectral dataset. The whole data is then written into the matrix $X\in \R^{M\times N}_{\geq 0},$ where each entry denotes the intensity of the detected particles of a specific mass-to-charge ratio (m/z-value). Typically, the measured spectra are ordered row-wise to the matrix, such that every column corresponds to an intensity plot of the whole tissue slide for a specific m/z-value (m/z-image).
            
            Typical data sizes range from $1\cdot 10^4$ to one million spectra and m/z-images respectively. Furthermore, MALDI datasets are naturally nonnegative. These properties make the NMF an ideal analyzing tool for MALDI imaging datasets due to the nonnegativity constraints leading to a meaningful physical interpretation of the acquired matrices $U$ and $V.$ For more details on the application of NMF models to MALDI imaging problems and the interpretation of the factorization matrices, we refer the reader to the works \cite{Leuschner:2018:NMFSupervised,Fernsel:2018:Survey}.
            \begin{figure}[tbhp]
                \centering
                \subfloat[Histological image after H\&E staining]{\label{fig:dataSetAndGroundTruth:dataSet}\includegraphics[width=0.4\textwidth]{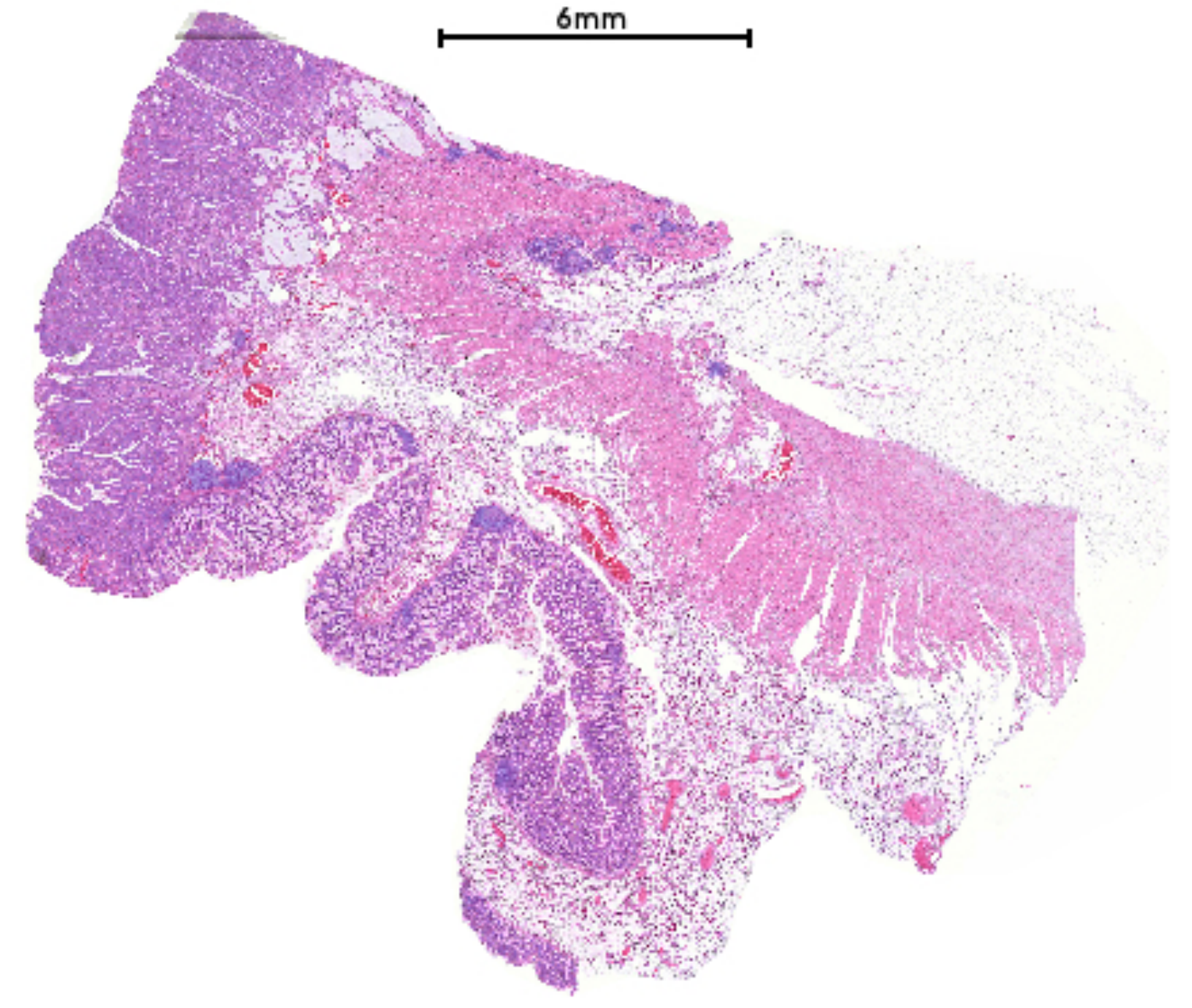}}\hfill
                \subfloat[Ground truth]{\label{fig:dataSetAndGroundTruth:groundTruth}\includegraphics[trim=0.3cm 2.4cm 4.4cm 1.6cm, clip, width=0.5\textwidth]{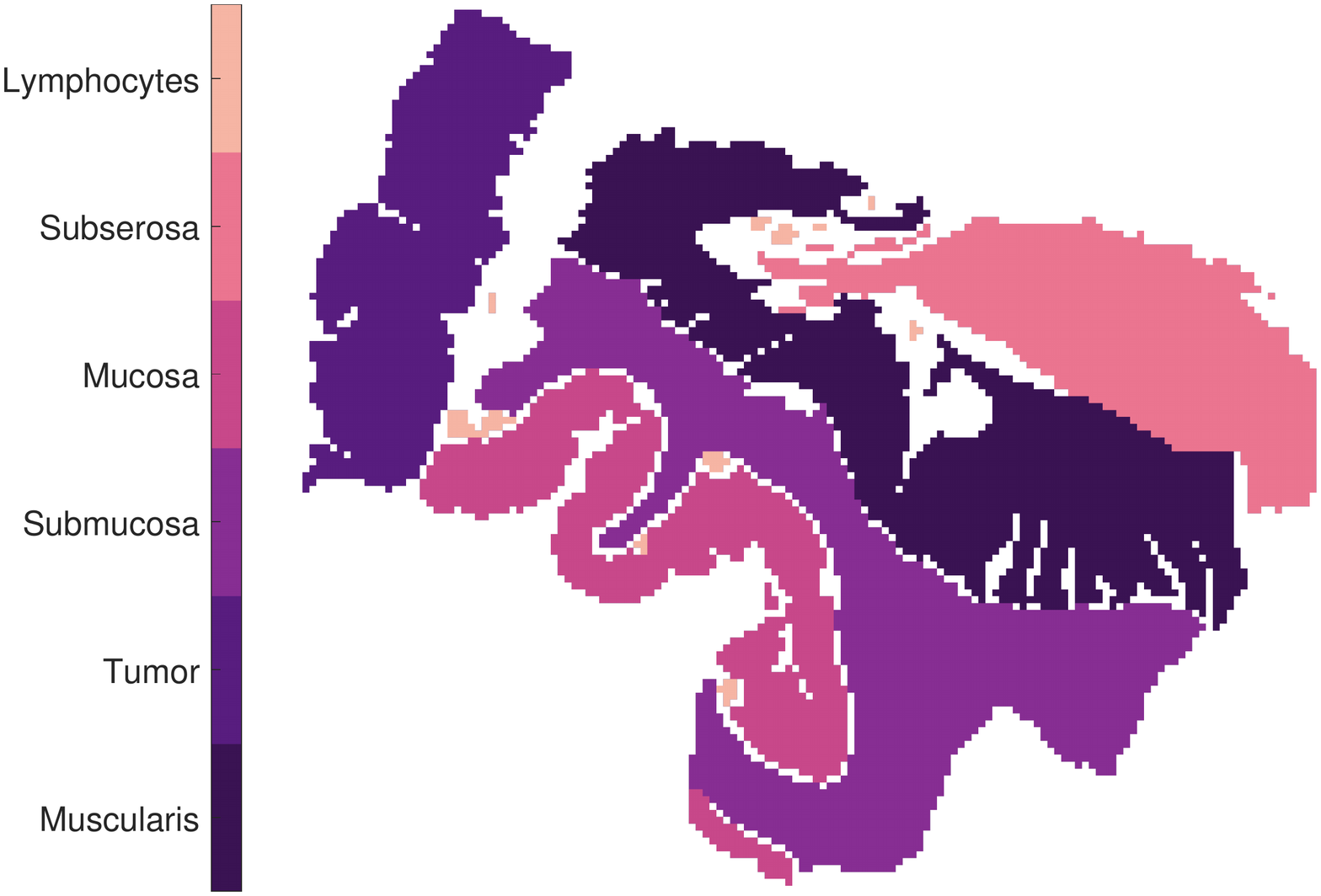}}
                \caption{Histological Image of the considered MALDI dataset after H\&E staining (\cref{fig:dataSetAndGroundTruth:dataSet}) and the histological annotation (\cref{fig:dataSetAndGroundTruth:groundTruth}).}
                \label{fig:dataSetAndGroundTruth}
            \end{figure}
            
            \cref{fig:dataSetAndGroundTruth} shows the human colon tissue dataset used for the subsequent numerical evaluation of the methods. \cref{fig:dataSetAndGroundTruth:dataSet} is the histological image after the application of a so-called hematoxylin and eosin (H\&E) staining, which allows a distinction between different tissue types. \cref{fig:dataSetAndGroundTruth:groundTruth} shows the histological annotation, which divides the dataset into six different classes and constitutes the ground truth for the subsequent numerical evaluation. The human colon dataset consists of 12049 acquired spectra, each containing 20000 measured m/z values covering a mass range of \SIrange{600}{4000}{\Dalton}. However, we only consider the actual annotated spectra to ensure that each considered spectra can be reasonably classified in one of the 6 classes shown in \cref{fig:dataSetAndGroundTruth:groundTruth}. Hence, we restrict ourselves to 8725 spectra leaving out the zero and the non-annotated ones, which leads to the nonnegative data matrix $X\in \R^{8725 \times 20000}_{\geq 0}$ for the numerical experiments.

        \subsection{Choice of Separated Methods} \label{subsec: Choice of Separated Methods}
            As discussed in \cref{subsec: Separated Methods}, we consider the separated approaches based on the workflow given in \cref{alg:ONMF-TV} as baseline for the comparison with the combined methods presented in \cref{subsec: Combined Methods}. In this section, we specify shortly seven ONMF algorithms to compute the clustering in \cref{alg:ONMF-TV_Line3} of \cref{alg:ONMF-TV}, which are based on different works throughout the literature and will be used for the numerical comparison in \cref{subsec: Results and Discussion} (see \cref{tab:listOfSepMethods}). Numerical tests for the optimization approach in \cite{PanNg:2018:CompMethod}, which is based on a sparsity and nuclear norm minimization, did not lead to satisfactory results and hence will not be presented in this work.
            
            The choice of the hyperparameters and details on the initialization of the matrices of all separated methods are discussed in \cref{subsec: Setup} and \cref{app:subsec:Parameter Choice of Separated Methods}.
            
            \begin{table}[tbhp]
                \caption{Designations of the considered separated methods (left column) and short explanation of the corresponding ONMF algorithm (right column).}
                \ra{1.1}
                \resizebox{\columnwidth}{!}{%
                \pgfplotstabletypeset[
                    col sep=&, row sep=\\,
                    string type,
                    every head row/.style={
                        before row={
                            \toprule
                        },
                        after row=\midrule,
                    },
                    every last row/.style={
                        after row=\bottomrule},
                    every even row/.style={
                        before row={\rowcolor[gray]{0.9}}
                    },
                    columns/alg/.style     ={column name=\textbf{Separated Method},column type=l},
                    columns/descr/.style    ={column name=\textbf{Description},column type=l},
                    ]{
                        alg & descr \\
                        \texttt{K-means-TV} & \makecell[{{p{12cm}}}]{Classical K-means clustering algorithm.} \\
            			\texttt{ONMF-TV-Choi} & \makecell[{{p{12cm}}}]{Alternating multiplicative update rules based on \cite{Choi:2008:CompMethod}.} \\
            			\texttt{ONMF-TV-Ding} & \makecell[{{p{12cm}}}]{Alternating multiplicative update rules based on \cite{Ding:2006:triFac}.} \\
            			\texttt{ONMF-TV-Pompili1} & \makecell[{{p{12cm}}}]{Alternating expectation-maximization algorithm similar to the default spherical K-means algorithm \cite{Pompili:2014:CompMethod,Banerjee:2003:SphericalK-means}.} \\
            			\texttt{ONMF-TV-Pompili2} & \makecell[{{p{12cm}}}]{Alternating algorithm based on an augmented Lagrangian approach with strict orthogonality constraints \cite{Pompili:2014:CompMethod}. Nonnegativity is obtained asymptotically by using a quadratic penalty.}\\
            			\texttt{ONMF-TV-Kimura} & \makecell[{{p{12cm}}}]{Hierarchical alternating least squares algorithm, which is applied column-wise on $U$ and row-wise on $V$ \cite{Kimura:2015:CompMethod}.} \\
            			\texttt{ONMF-TV-Li} & \makecell[{{p{12cm}}}]{Hierarchical alternating least squares algorithm with approximate orthogonality constraints and subsequent projection steps to ensure nonnegativity \cite{LiEtal:2015:CompMethod}.}\\
                    }
                }
        		\label{tab:listOfSepMethods}
            \end{table}
            
        \subsection{Setup} \label{subsec: Setup}
            In this section, we describe the initalization methods and stopping criteria, the calculation of the final hard clustering, the choice of the various hyperparameters, the used cluster validation measures and we give some further details on the general numerical setup.
            
            For every considered separated and combined method in \cref{subsec: Separated Methods} and \cref{subsec: Combined Methods}, we perform 30 replicates of the experiment to get an impression of the performance stability, since the used initialization methods are partially based on randomly selected data points or matrices. For each method, we use the same random seed in the beginning of the experiment. Furthermore, we measure for each replicate the computational time including the time for the calculation of the initialization of the factorization matrices.
            
            For the evaluation of the clusterings, we use several different clustering validation measures discussed in \cite{XiongLi:2014:ClusteringValidation}. Due to the known ground truth of the data, we restrict ourselves to external clustering validation measures. Based on the results in this work, we primarily consider the so-called normalized van Dongen criterion (VD\textsubscript{n}) and the normalized variation of information (VI\textsubscript{n}), since they give the most representative, quantitative measures in most of the general cases. Furthermore, they are normalized into the $[0, 1]$ range and give reasonable results in cases there are clusters without any corresponding datapoints, which is different from classical measures like the purity. In addition, we consider as secondary measures the entropy (E) and the non-normalized versions of the VD\textsubscript{n} and VI\textsubscript{n}. which we will call VD and VI respectively. For all considered measures, it holds that a lower value indicates a better clustering performance.
            
            To provide the definition of the clustering validation measures, we denote $n_{k\tilde{k}}$ as the number of data points in cluster $\mathcal{I}_k$ from class $\mathcal{C}_{\tilde{k}}$ for $k,\tilde{k}\in \{ 1, \dots, K\}$ and define
            \begin{align*}
                n&\coloneqq \sum_{k,\tilde{k}=1}^K n_{k\tilde{k}}, & n_{k,\bullet}&\coloneqq \sum_{\tilde{k}=1}^K n_{k\tilde{k}}, & n_{\bullet, \tilde{k}}&\coloneqq \sum_{k=1}^K n_{k\tilde{k}}, \\
                p_{k\tilde{k}} &\coloneqq \frac{n_{k\tilde{k}}}{n}, & p_k &\coloneqq \frac{n_{k,\bullet}}{n}, & \tilde{p}_{\tilde{k}}&\coloneqq \frac{n_{\bullet,\tilde{k}}}{n}.
            \end{align*}
            Using this notation, the definition of all considered clustering validation measures are given in \cref{tab:listOfClustValidMeasures}.
            \begin{table}[tbhp]
                \caption{Definitions of all considered clustering validation measures.}
                \ra{1.1}
                \resizebox{\columnwidth}{!}{%
                \pgfplotstabletypeset[
                    col sep=&, row sep=\\,
                    string type,
                    every head row/.style={
                        before row={
                            \toprule
                        },
                        after row=\midrule,
                    },
                    every last row/.style={
                        after row=\bottomrule},
                    every even row/.style={
                        before row={\rowcolor[gray]{0.9}}
                    },
                    columns/measure/.style      ={column name=\textbf{Measure},column type=l},
                    columns/def/.style          ={column name=\textbf{Definition},column type=l},
                    columns/range/.style        ={column name=\textbf{Range},column type=l},
                    ]{
                        measure & def & range \\
                        \makecell[{{p{5cm}}}]{Entropy (E)} & $\displaystyle -\sum_{k=1}^K p_k \left( \sum_{\tilde{k}=1}^K \frac{p_{k\tilde{k}}}{p_k} \log\left( \frac{p_{k\tilde{k}}}{p_k} \right) \right)$ & $[0, \log(K)]$ \\
            			\makecell[{{p{5cm}}}]{Variation of Information (VI)} & $\displaystyle -\sum_{k=1}^K p_k \log(p_k) - \sum_{\tilde{k}=1}^K \tilde{p}_{\tilde{k}} \log(\tilde{p}_{\tilde{k}}) - 2\sum_{k,\tilde{k}=1}^K p_{k\tilde{k}} \log \left( \frac{p_{k\tilde{k}}}{p_k \tilde{p}_{\tilde{k}}} \right) $ & $[0, 2\cdot\log(K)]$ \\
            			\makecell[{{p{5cm}}}]{Van Dongen criterion (VD)} & $\displaystyle \frac{2n - \sum_{k=1}^K \max_{\tilde{k}}\{ n_{k\tilde{k}} \} - \sum_{\tilde{k}=1}^K \max_{k}\{ n_{k\tilde{k}} \}}{2n} $ & $[0, 1)$ \\
            			\makecell[{{p{4cm}}}]{Normalized Variation of Information (VI\textsubscript{n})} & $\displaystyle 1 + 2\cdot \frac{\sum_{k,\tilde{k}=1}^K p_{k\tilde{k}} \log \left( p_{k\tilde{k}}/ (p_k \tilde{p}_{\tilde{k}}) \right)}{\sum_{k=1}^K p_k \log(p_k) + \sum_{\tilde{k}=1}^K \tilde{p}_{\tilde{k}} \log(\tilde{p}_{\tilde{k}})}$ & $[0, 1]$ \\
            			\makecell[{{p{4cm}}}]{Normalized van Dongen criterion (VD\textsubscript{n})} & $\displaystyle \frac{2n - \sum_{k=1}^K \max_{\tilde{k}}\{ n_{k\tilde{k}} \} - \sum_{\tilde{k}=1}^K \max_{k}\{ n_{k\tilde{k}} \}}{2n - \max_{k}\{ n_{k,\bullet} \} - \max_{\tilde{k}}\{ n_{\bullet,\tilde{k}} \}}$ & $[0, 1]$ \\
                    }
                }
        		\label{tab:listOfClustValidMeasures}
            \end{table}
            
            Concerning the initialization approaches, we consider either the classical K-means++ method or an initialization based on the singular value decomposition of the data matrix $X$ by following the works \cite{HalkoEtal:2011:NMFInit1,HalkoEtal:2011:NMFInit2,BoutsidisGallopoulos:2008:NMFInit}. The initialization method, which leads to a better clustering stability and performance in terms of the VD\textsubscript{n}, is chosen. The concrete choices for every ONMF model is described in \cref{app:sec:Parameter Choice}.
            Regarding the stopping criteria, we simply set for the considered ONMF models a maximal iteration number until a sufficient convergence is reached except for \texttt{K-means-TV}, \texttt{ONMF-TV-Pompili1} and \texttt{ONMF-TV-Pompili2}, where we use the internal stopping criteria of the respective algorithms. For more information, we refer the reader to \cref{app:sec:Parameter Choice} and the work \cite{Pompili:2014:CompMethod}.
            
            Another aspect is the computation of the final clustering based on the cluster membership matrix $U.$ Most of the considered ONMF models yield a cluster membership matrix $U$ having multiple positive entries in each row, which is related to a soft clustering. To obtain the final hard clustering, we assign every data point to the cluster, which corresponds to the column in $U,$ where the maximal value in the row is attained. In the case that there are two or more equal entries in one row, we choose the cluster by a random choice.
            
            A main part of the whole workflow of the numerical evaluation is the choice of the various hyperparameters of the considered ONMF models. In particular, these include the regularization parameters $\sigma_1, \sigma_2$ and $\tau$ of the combined methods and the parameter $\tau$ of the TV postprocessing of the separated methods. For all considered methods, we perform various experiments for a wide range of all parameters and choose the configuration leading to the best stability and VD\textsubscript{n}. More details and the specific selection of the hyperparameters are given in \cref{app:sec:Parameter Choice}.
            
            Furthermore, we perform several projection steps to enhance the numerical stability. For all considered methods, we project the data matrix $X$ by applying $[\cdot ]_{>0}$ defined as in \cref{subsubsec: Multiplicative Update Rules}. Furthermore, we perform specifically for the multiplicative update rules and for \texttt{ONMF-TV-Kimura} the same projection step after the initialization of the matrices. For the combined method \texttt{ONMFTV-SPRING}, we perform an additional projection step of the parameter $\tau \eta_{U^{[i,j]}}$ used for the application of the proximal operator described in \cref{alg:ONMFTV-SPRING} to avoid too large parameters for the TV regularization. A similar projection is done for the step size in each gradient descent step. More details are given in \cref{app:subsec:ONMFTV-SPRING}.
    
        \subsection{Results and Discussion} \label{subsec: Results and Discussion}
\begin{figure}[tbhp]
                \centering
                \subfloat[\texttt{K-means-TV}]{\label{fig:bestClusters-UTV:k-Means-TV}\includegraphics[trim=5.31cm 2cm 4.27cm 1.1cm, clip, width=0.29\textwidth]{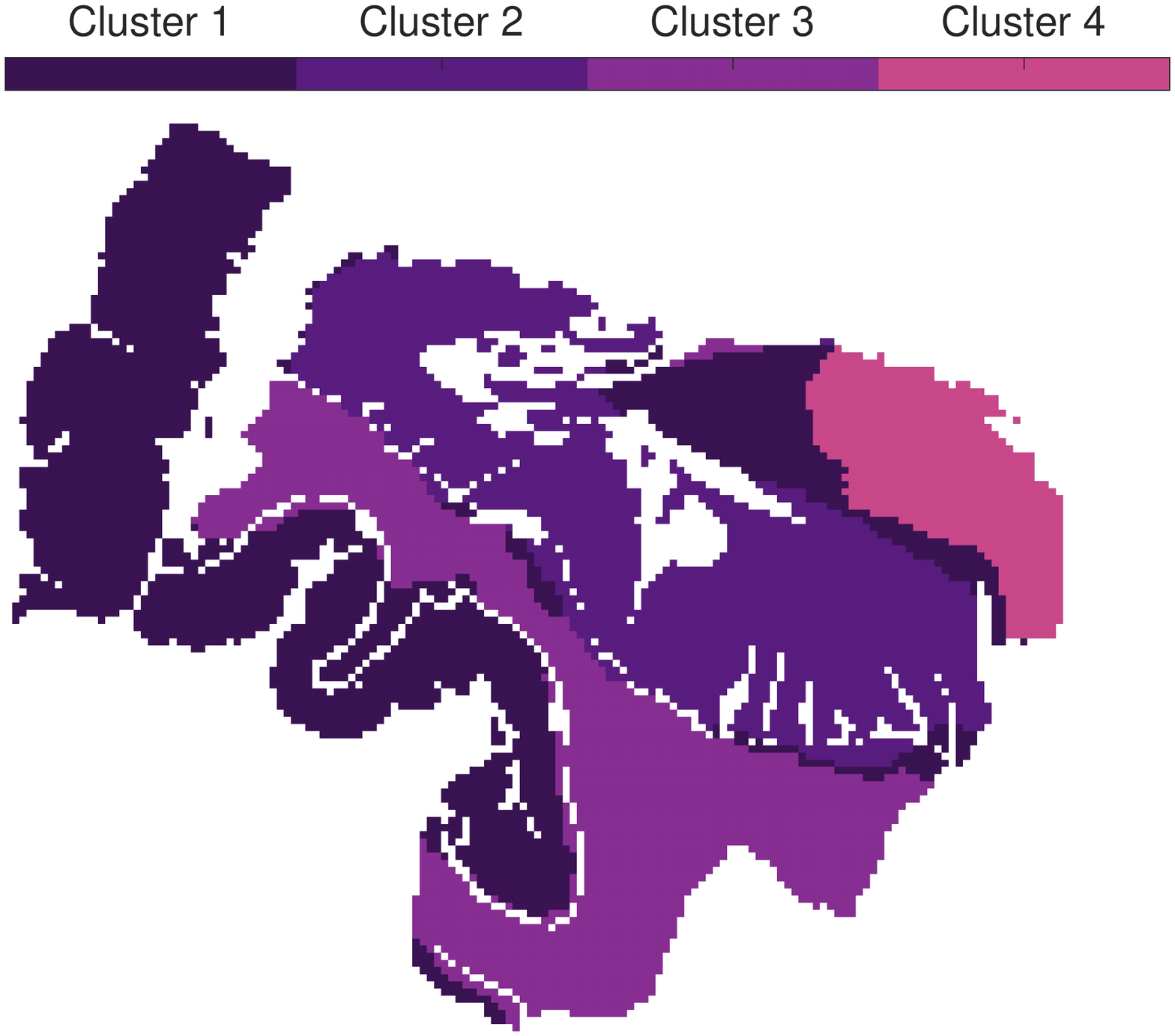}}\hfill
                \subfloat[\texttt{ONMF-TV-Choi}]{\label{fig:bestClusters-UTV:ONMF-TV-Choi}\includegraphics[trim=5.31cm 2cm 4.27cm 1.1cm, clip, width=0.29\textwidth]{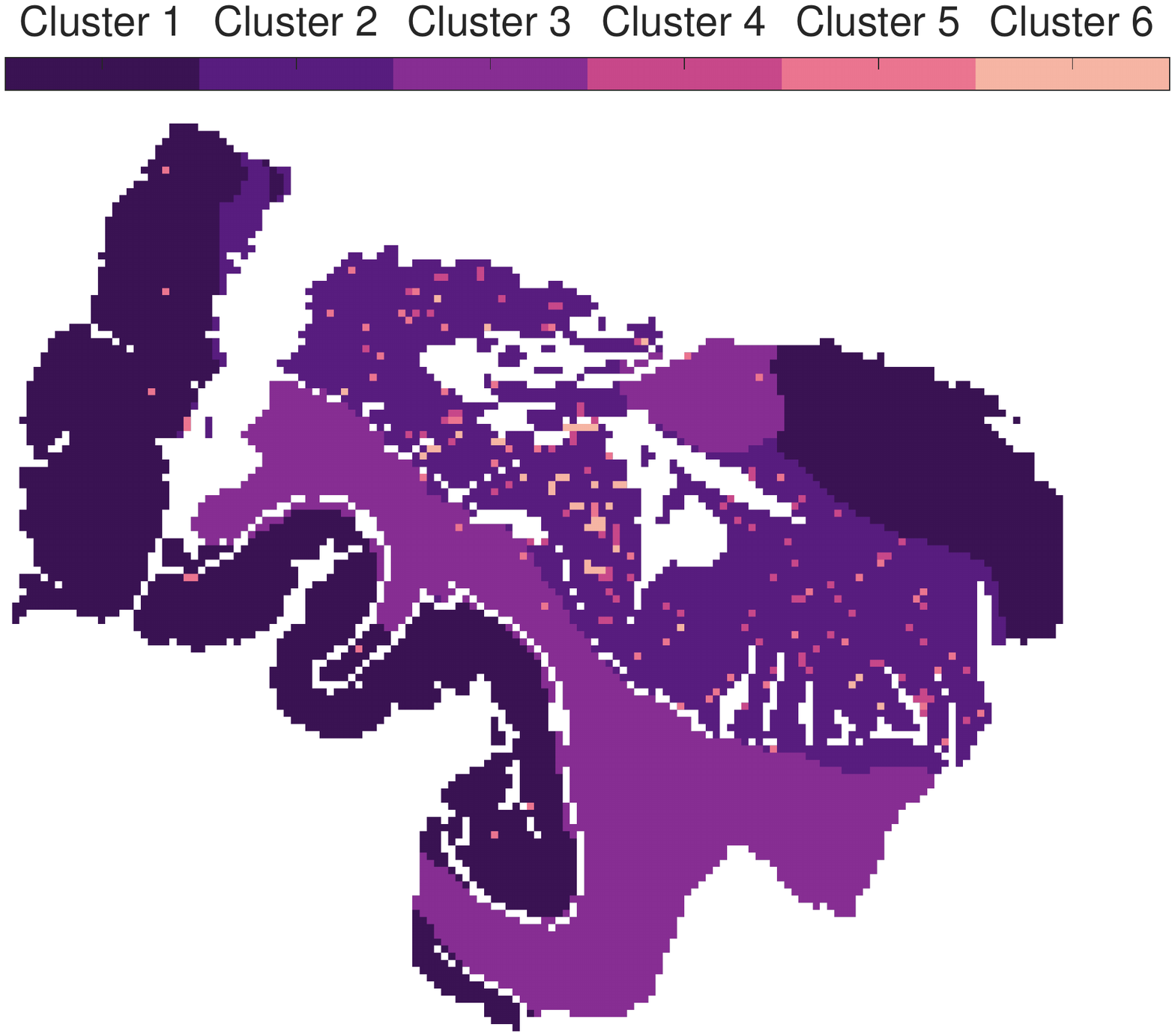}}\hfill
                \subfloat[\texttt{ONMF-TV-Ding}]{\label{fig:bestClusters-UTV:ONMF-TV-Ding}\includegraphics[trim=5.31cm 2cm 4.27cm 1.1cm, clip, width=0.29\textwidth]{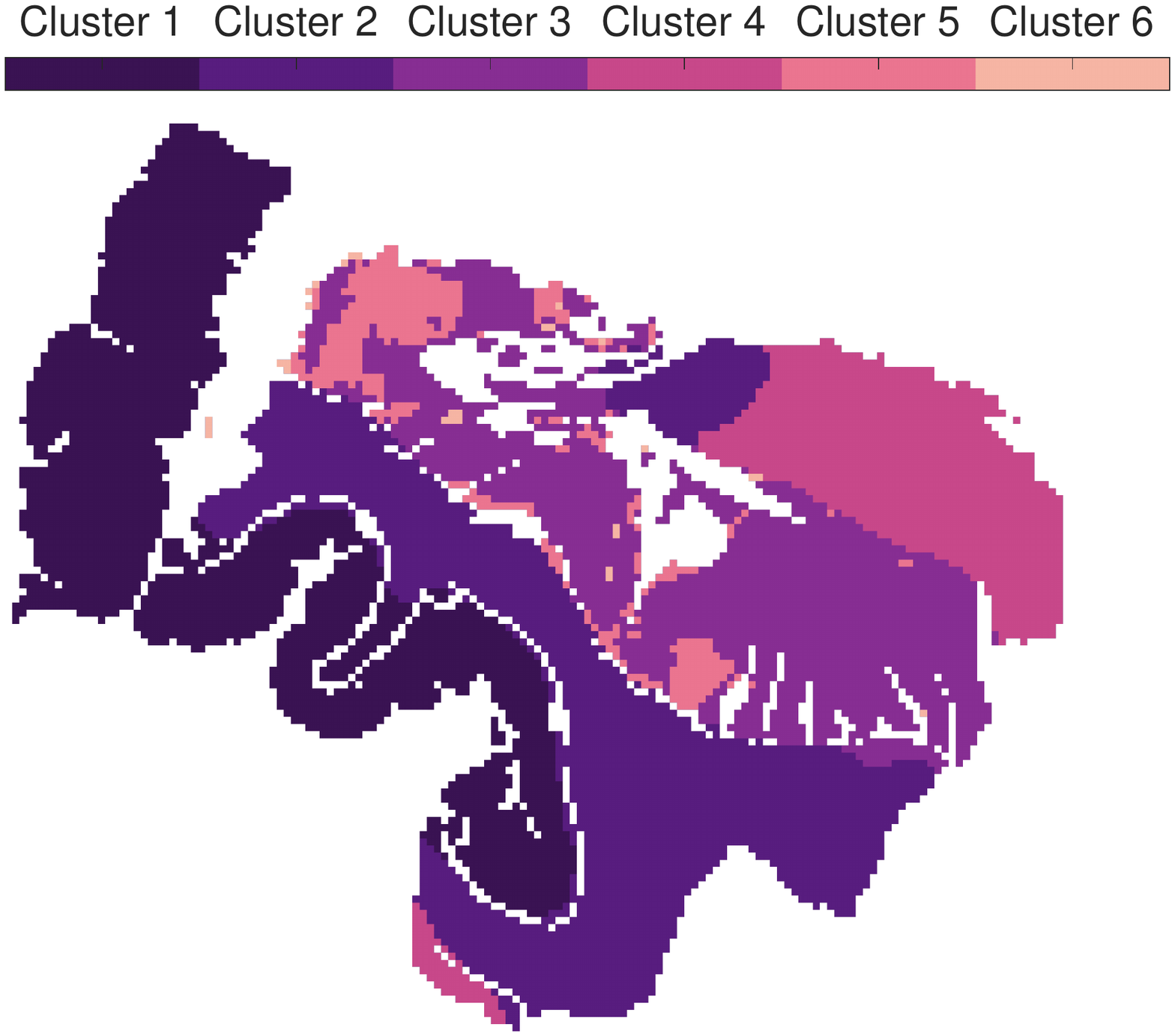}}\\
                \subfloat[\texttt{ONMF-TV-Pompili1}]{\label{fig:bestClusters-UTV:ONMF-TV-Pompili1}\includegraphics[trim=5.31cm 2cm 4.27cm 1.1cm, clip, width=0.29\textwidth]{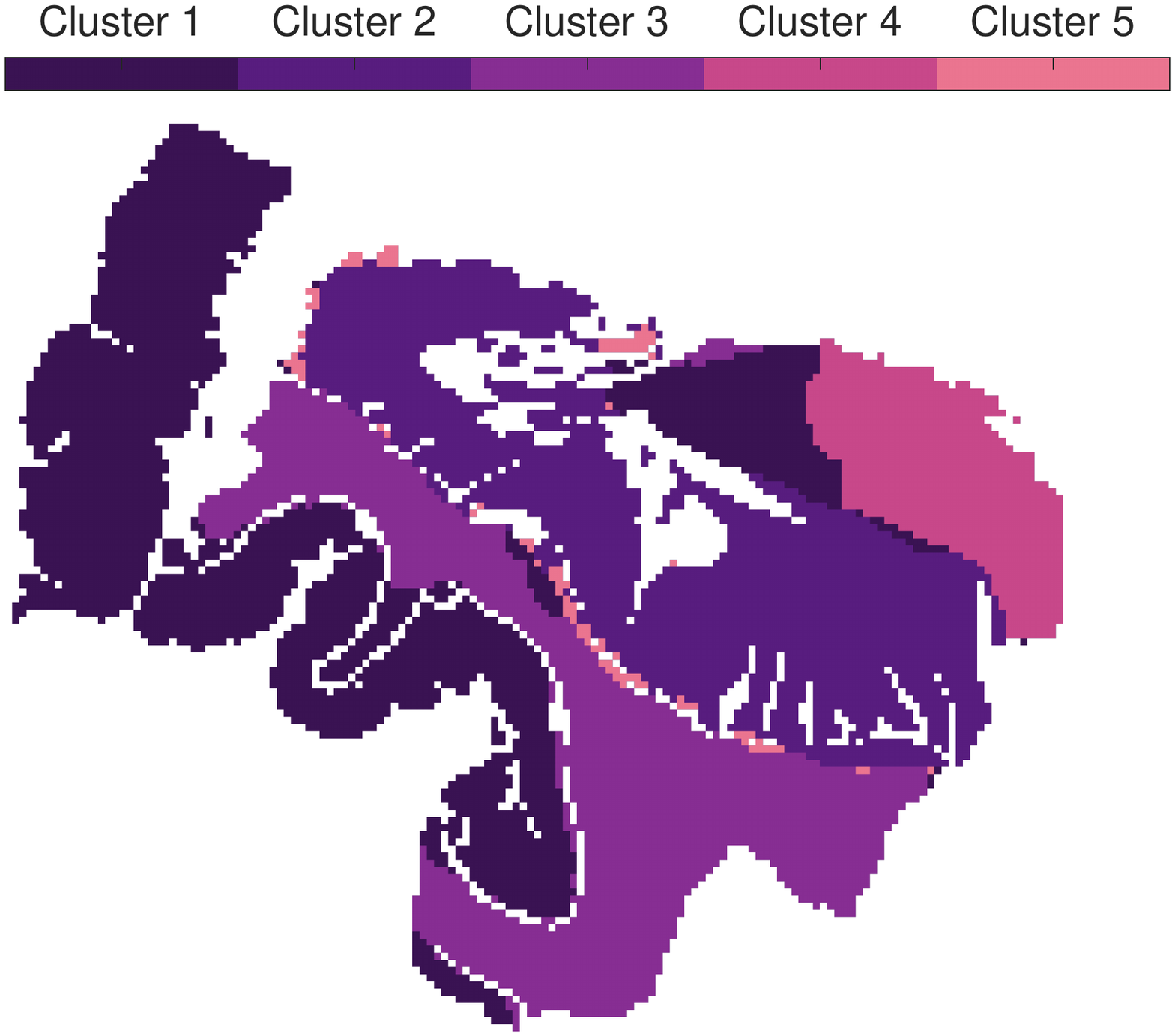}}\hfill
                \subfloat[\texttt{ONMF-TV-Pompili2}]{\label{fig:bestClusters-UTV:ONMF-TV-Pompili2}\includegraphics[trim=5.31cm 2cm 4.27cm 1.1cm, clip, width=0.29\textwidth]{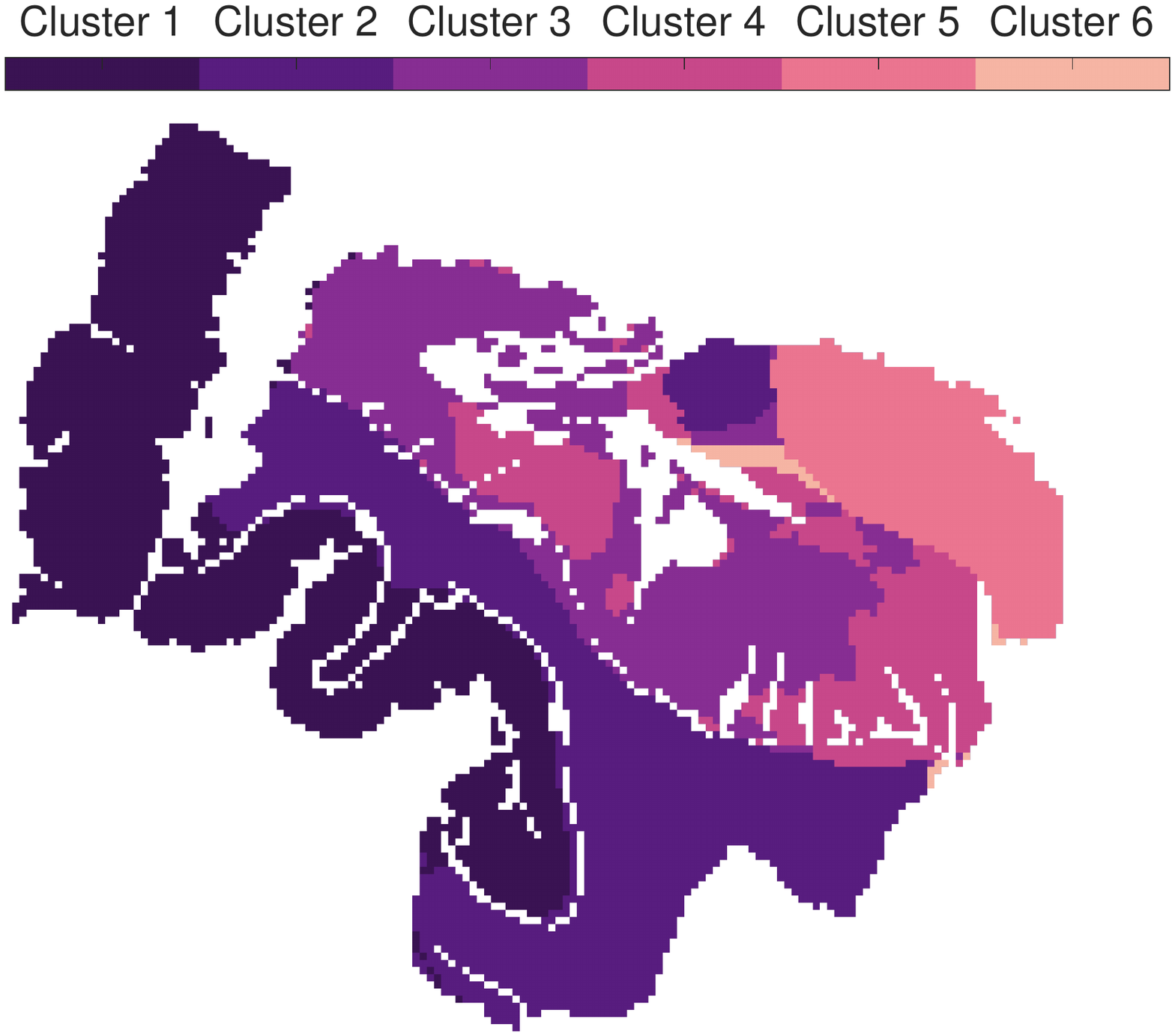}}\hfill
                \subfloat[\texttt{ONMF-TV-Kimura}]{\label{fig:bestClusters-UTV:ONMF-TV-Kimura}\includegraphics[trim=5.31cm 2cm 4.27cm 1.1cm, clip, width=0.29\textwidth]{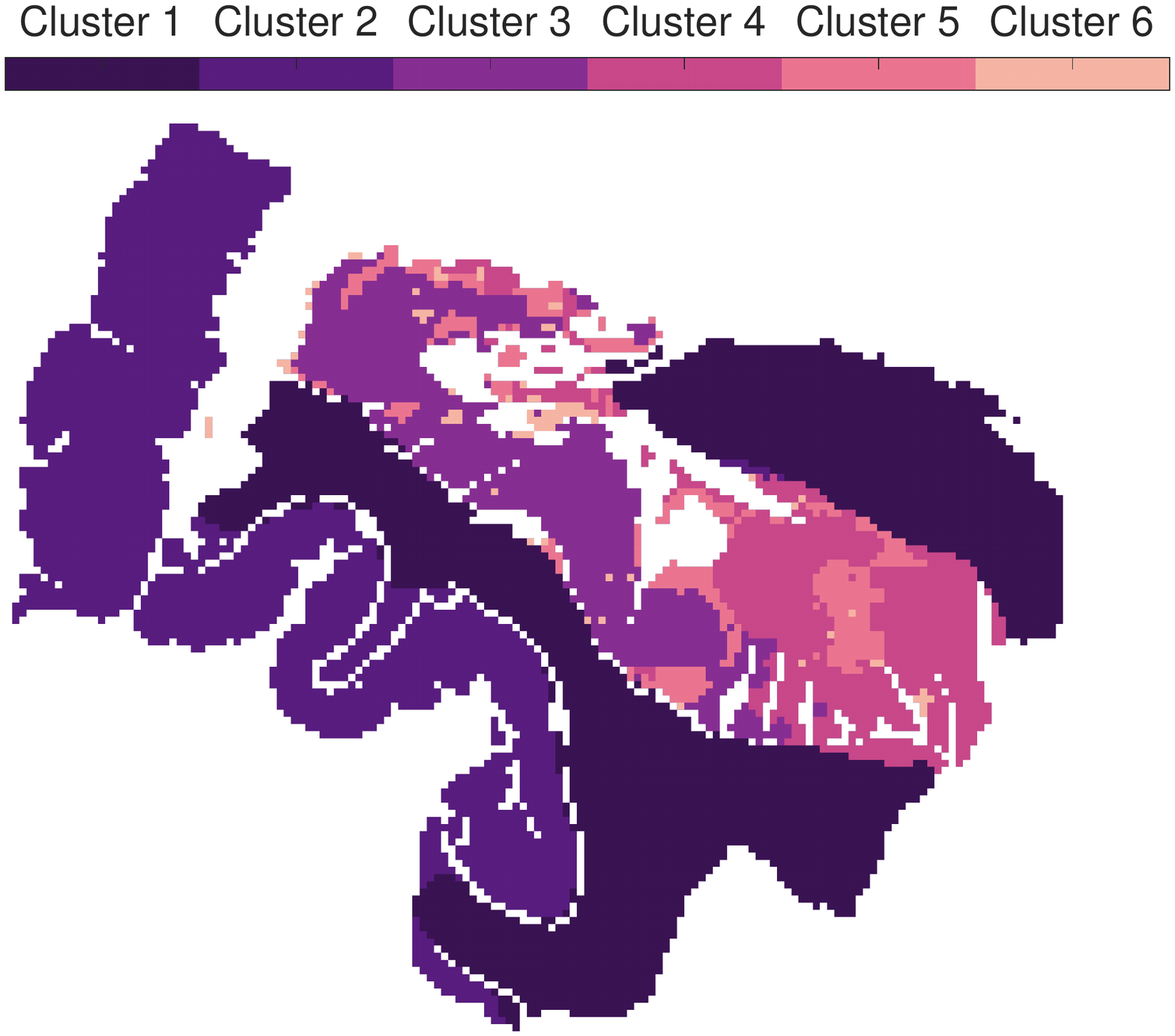}}\\
                \subfloat[\texttt{ONMF-TV-Li}]{\label{fig:bestClusters-UTV:ONMF-TV-Li}\includegraphics[trim=5.31cm 2cm 4.27cm 1.1cm, clip, width=0.29\textwidth]{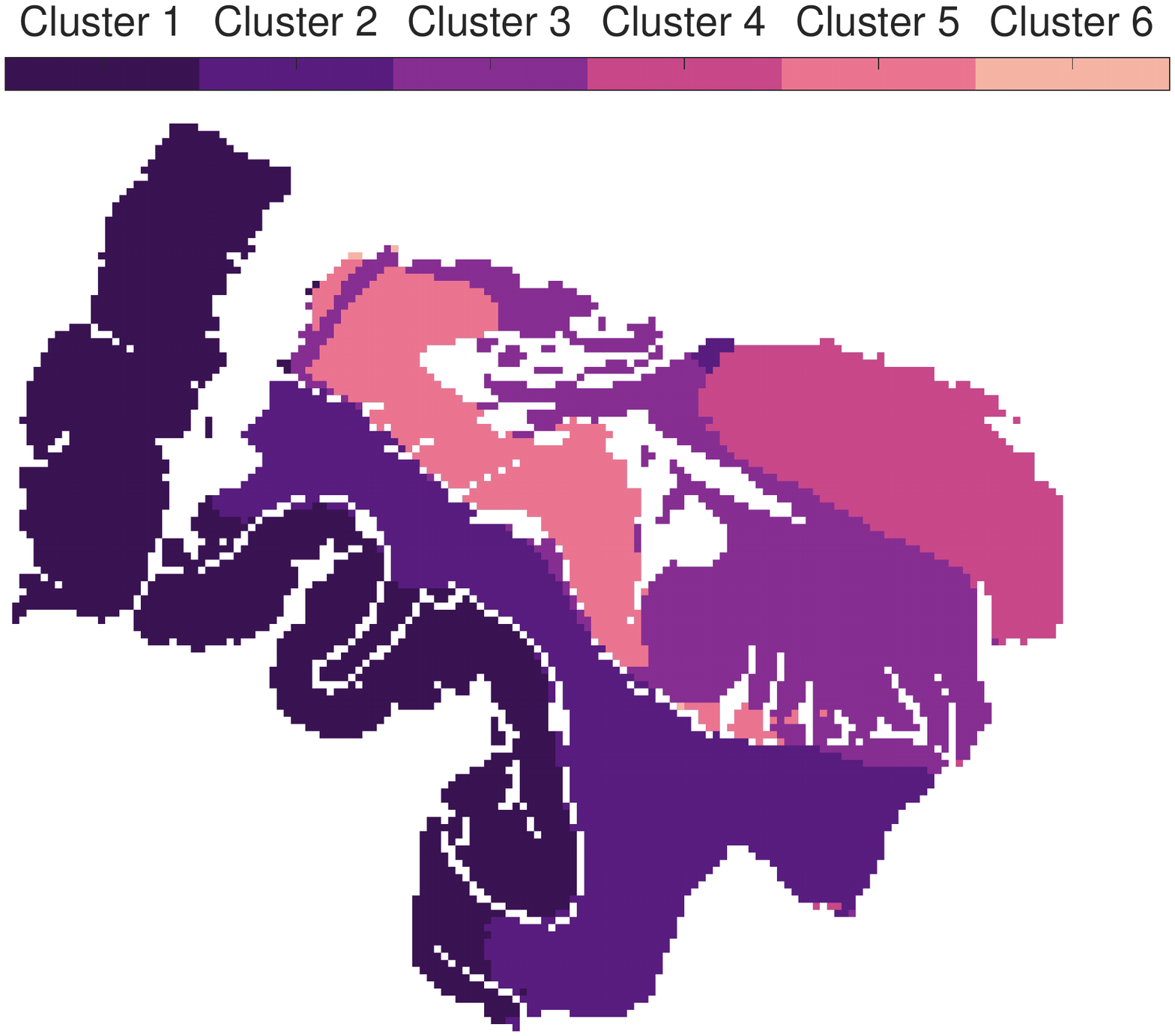}}\hfill
                \subfloat[\texttt{ONMFTV-MUL1}]{\label{fig:bestClusters-UTV:ONMFTV-MUL1}\includegraphics[trim=5.31cm 2cm 4.27cm 1.1cm, clip, width=0.29\textwidth]{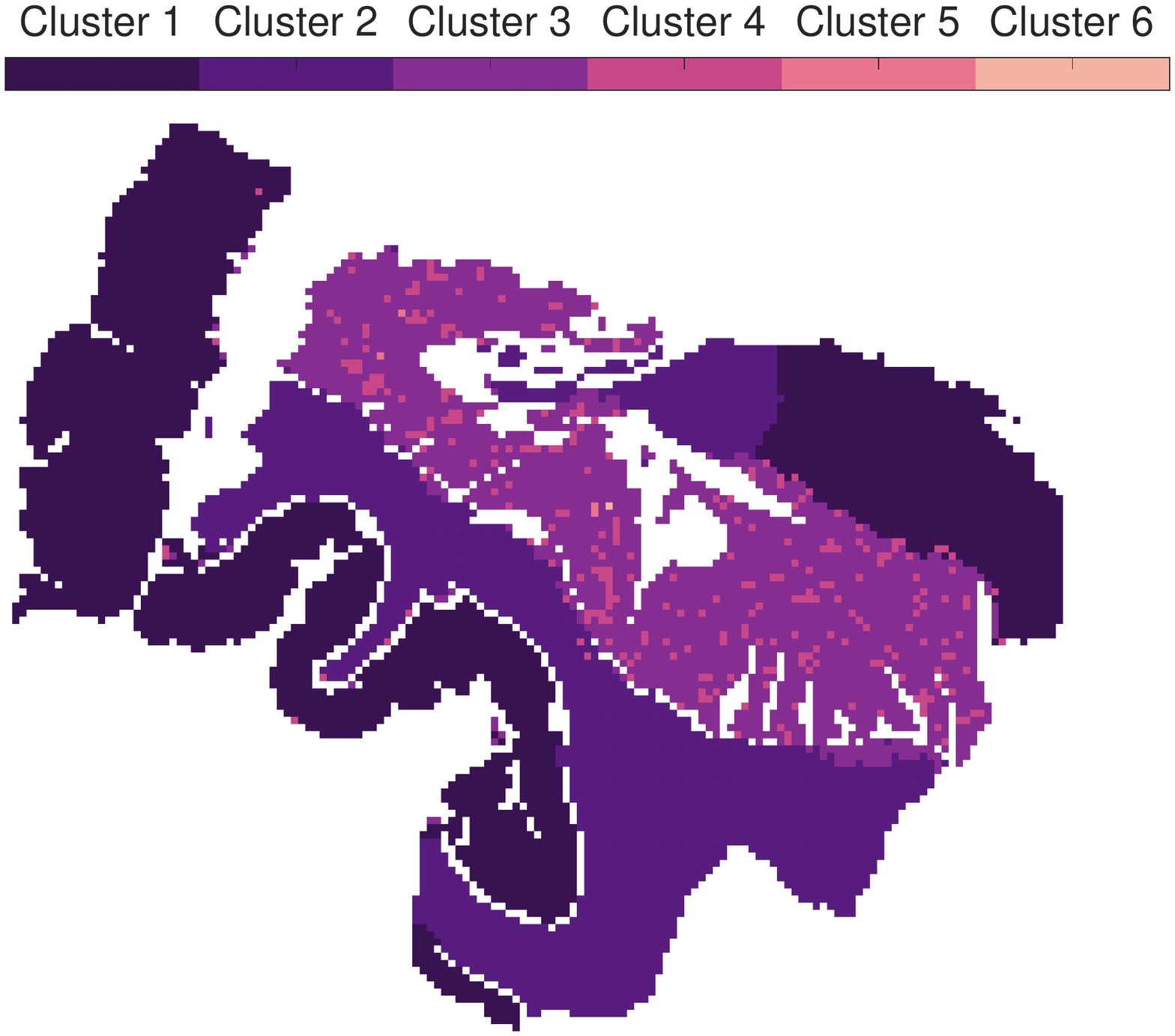}}\hfill
                \subfloat[\texttt{ONMFTV-MUL2}]{\label{fig:bestClusters-UTV:ONMFTV-MUL2}\includegraphics[trim=5.31cm 2cm 4.27cm 1.1cm, clip, width=0.29\textwidth]{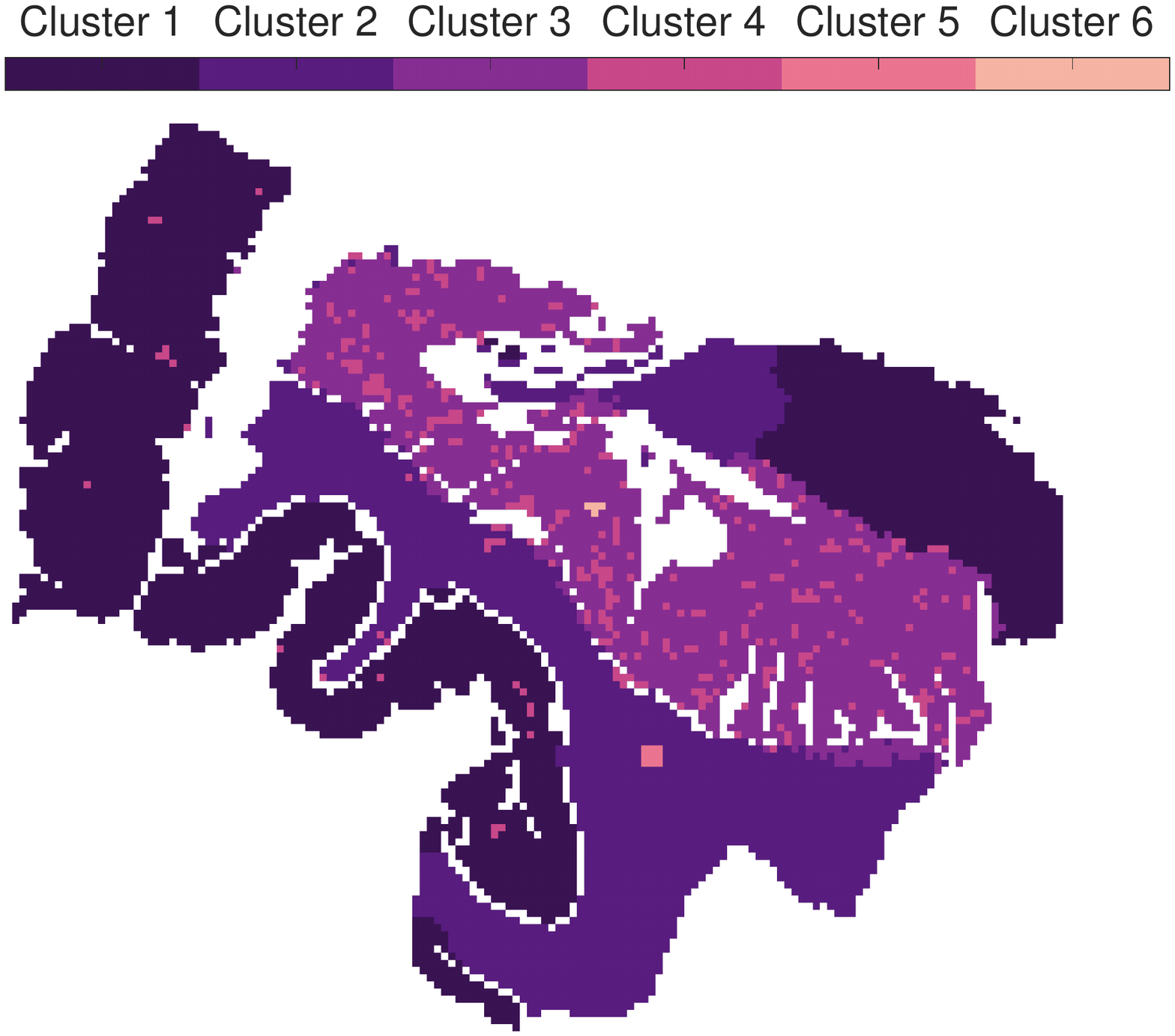}}\\
                \subfloat[\texttt{ONMFTV-PALM}]{\label{fig:bestClusters-UTV:ONMFTV-PALM}\includegraphics[trim=5.31cm 2cm 4.27cm 1.1cm, clip, width=0.29\textwidth]{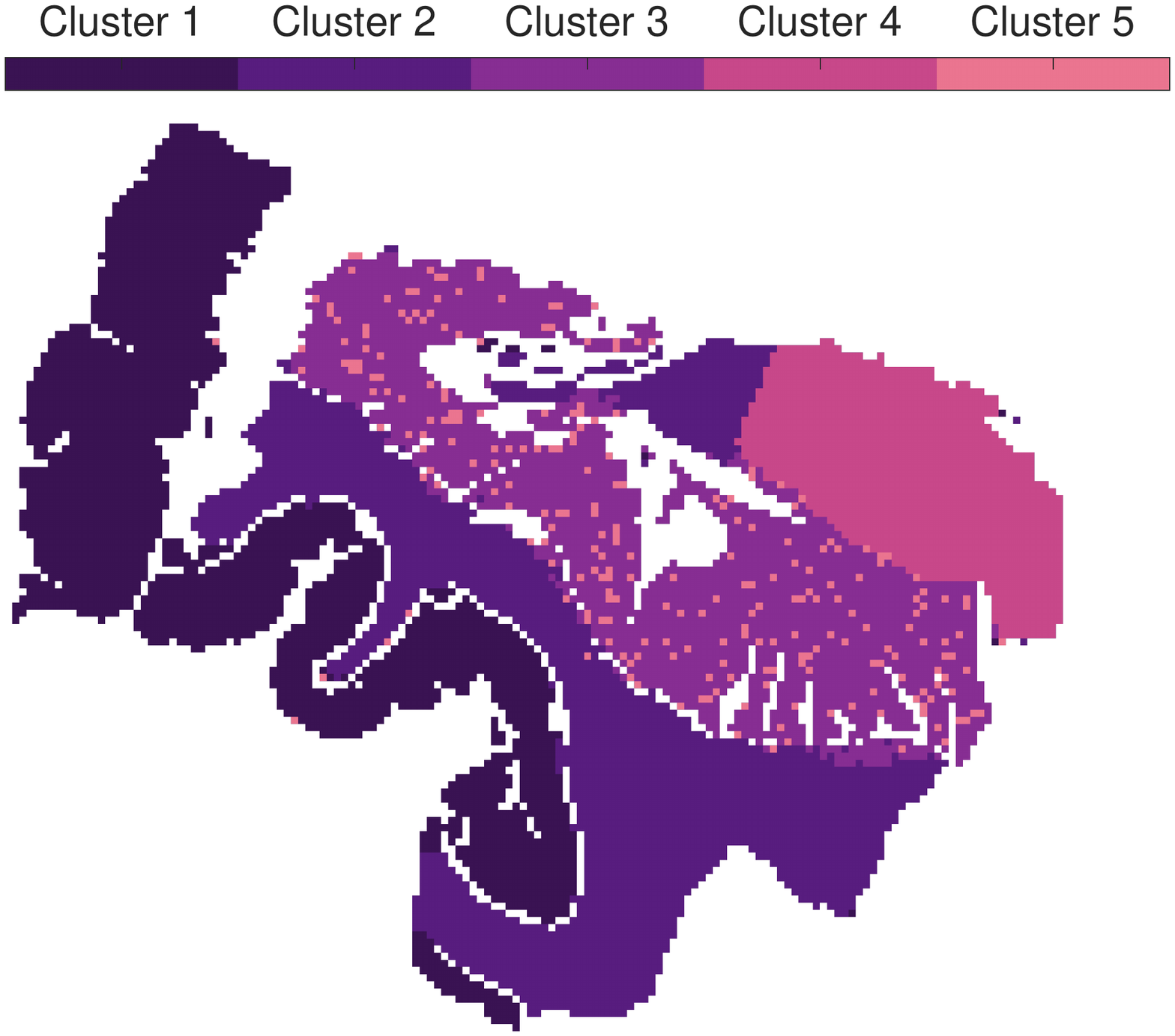}}\hfill
                \subfloat[\texttt{ONMFTV-iPALM}]{\label{fig:bestClusters-UTV:ONMFTV-iPALM}\includegraphics[trim=5.31cm 2cm 4.27cm 1.1cm, clip, width=0.29\textwidth]{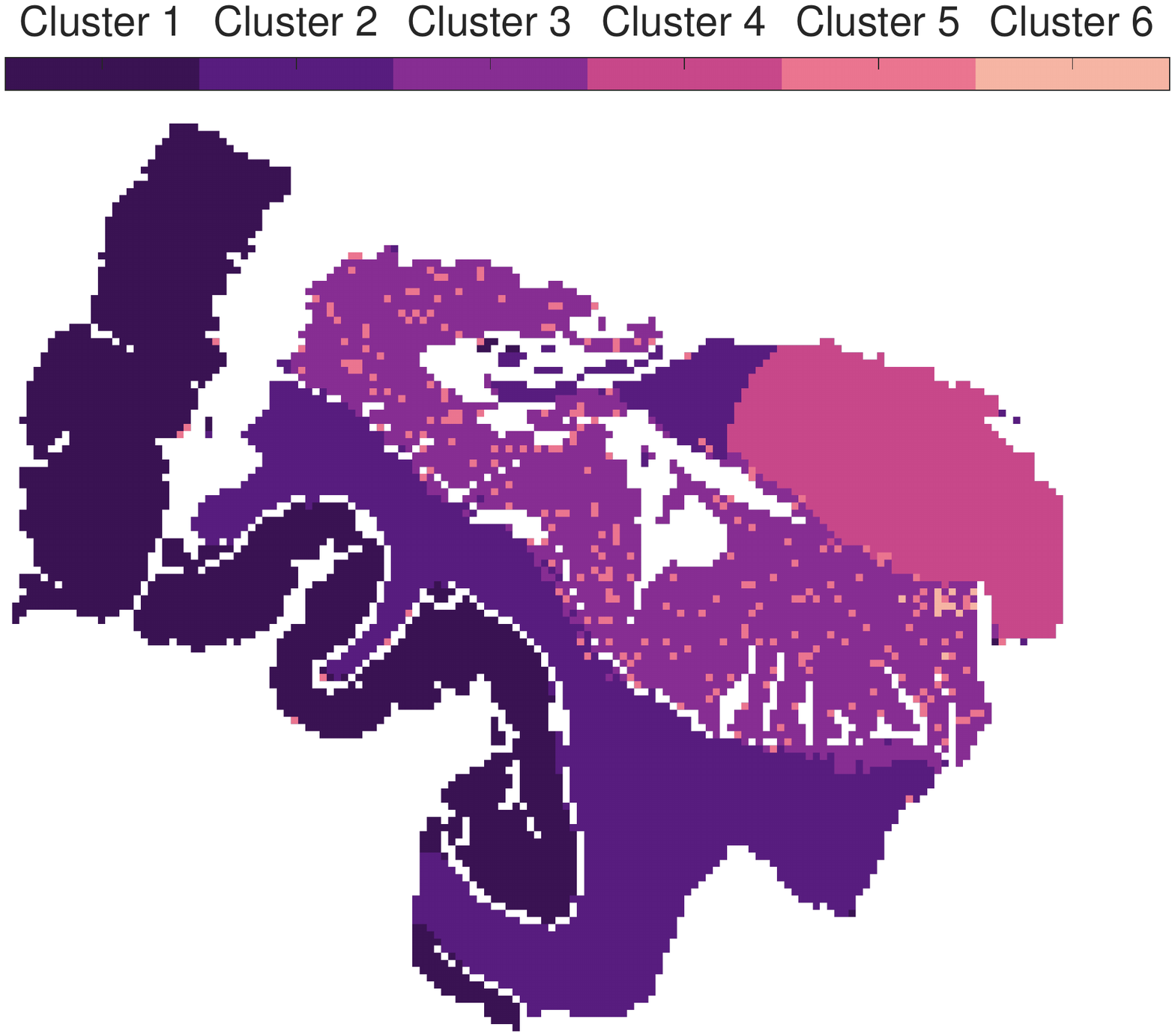}}\hfill
                \subfloat[\texttt{ONMFTV-SPRING}]{\label{fig:bestClusters-UTV:ONMFTV-SPRING}\includegraphics[trim=5.31cm 2cm 4.27cm 1.1cm, clip, width=0.29\textwidth]{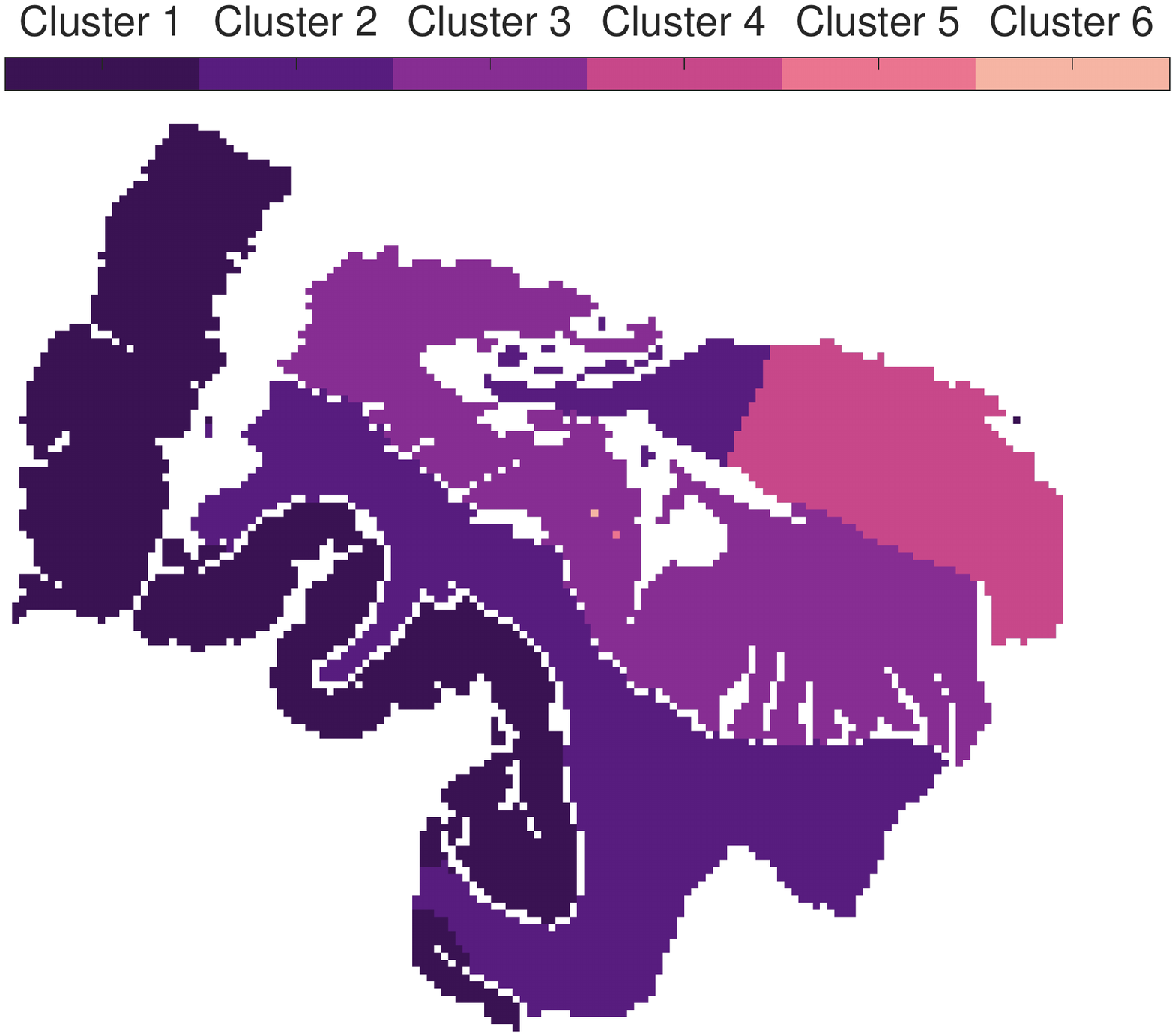}}
                \caption{Clusterings of all considered methods including the TV regularization. The best-performing replicate including the TV regularization is chosen based on the normalized van Dongen criterion.}
                \label{fig:bestClusters-UTV}
            \end{figure}

            In this section, we present and discuss the numerical results obtained in the evaluation discussed in the previous sections. 
            \begin{figure}[tbhp]
                \centering
                \subfloat[\texttt{K-means-TV}]{\label{fig:bestClusters-USep:k-Means-TV}\includegraphics[trim=5.31cm 2cm 4.27cm 1.1cm, clip, width=0.29\textwidth]{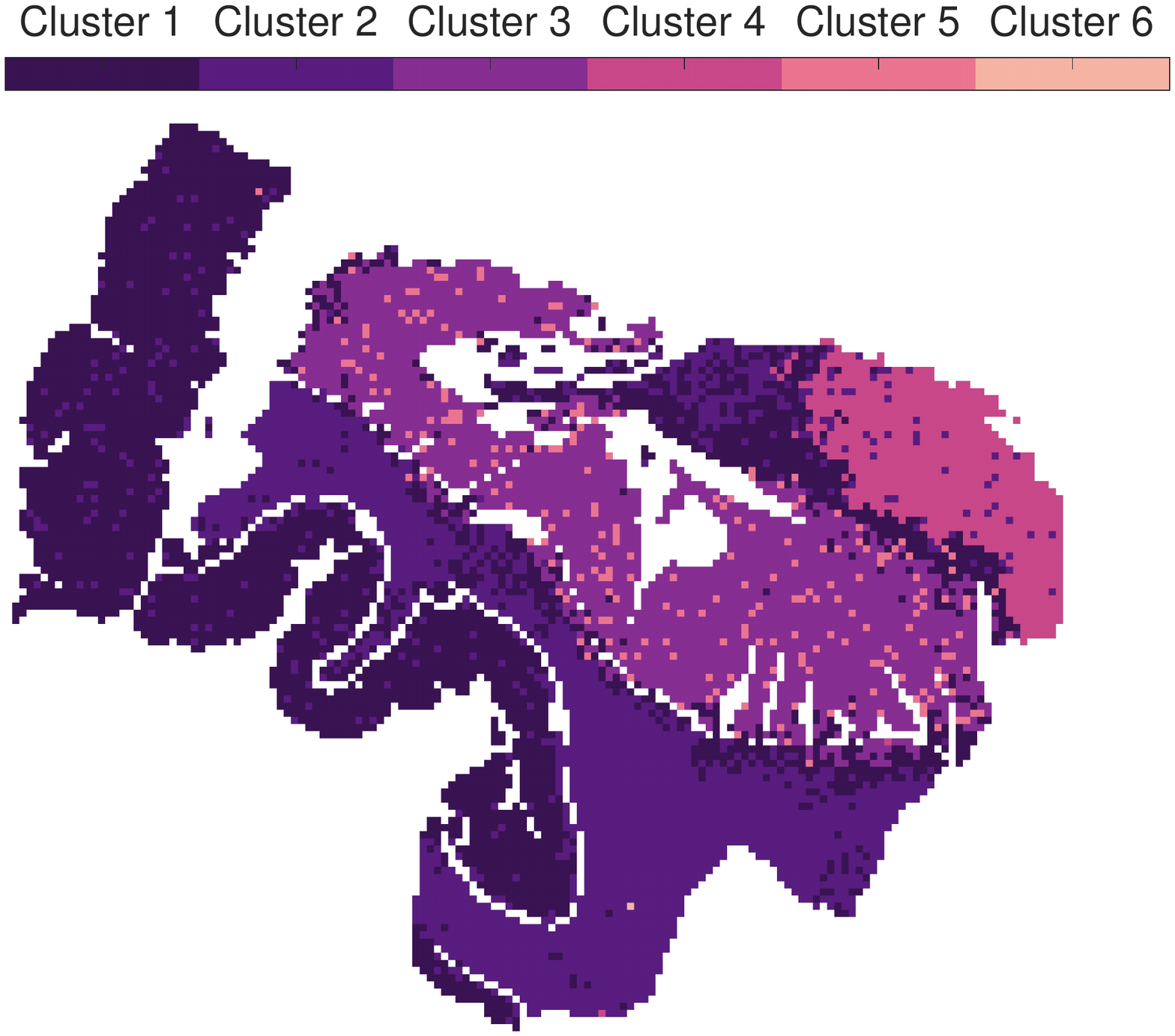}}\hfill
                \subfloat[\texttt{ONMF-TV-Choi}]{\label{fig:bestClusters-USep:ONMF-TV-Choi}\includegraphics[trim=5.31cm 2cm 4.27cm 1.1cm, clip, width=0.29\textwidth]{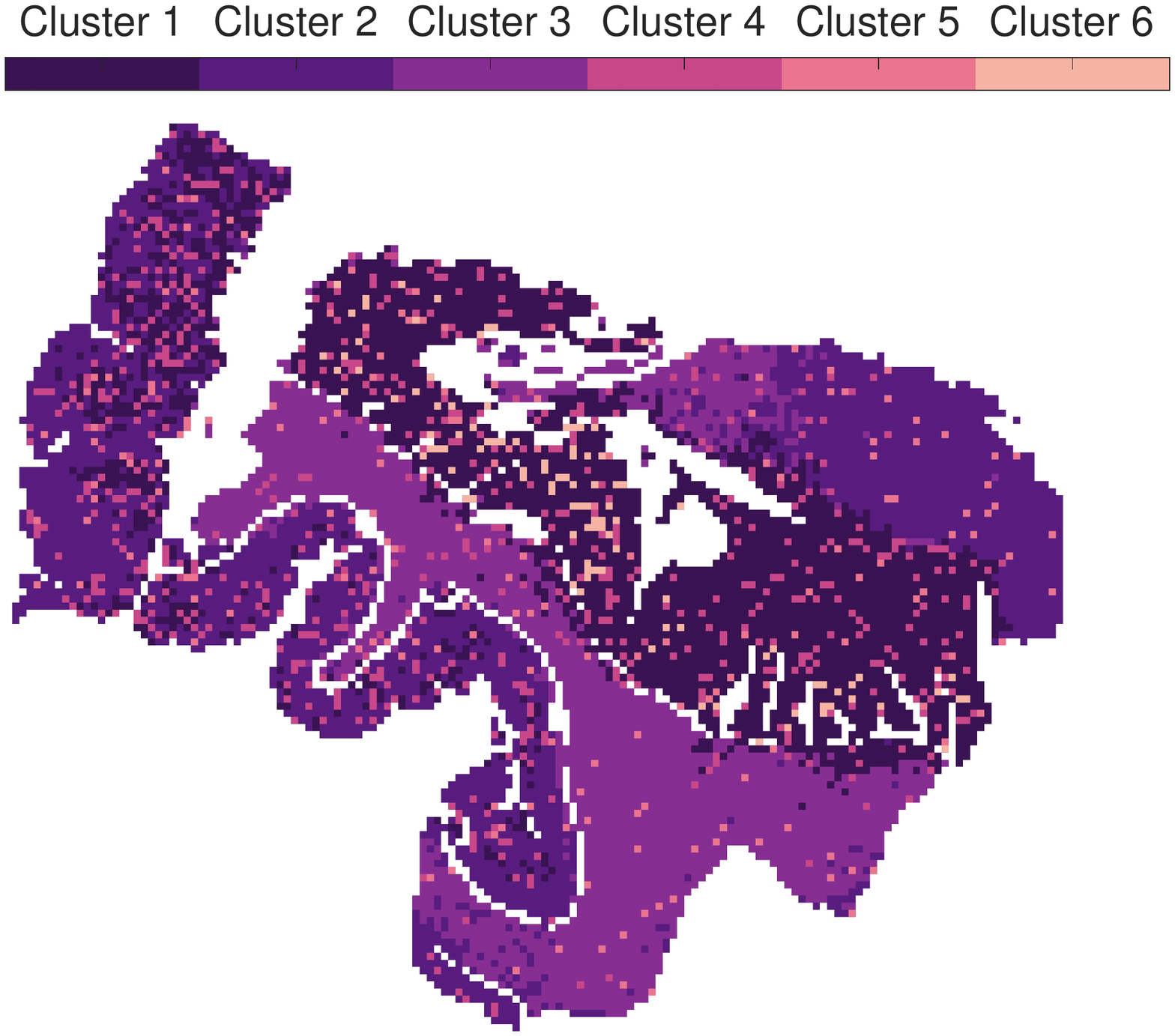}}\hfill
                \subfloat[\texttt{ONMF-TV-Ding}]{\label{fig:bestClusters-USep:ONMF-TV-Ding}\includegraphics[trim=5.31cm 2cm 4.27cm 1.1cm, clip, width=0.29\textwidth]{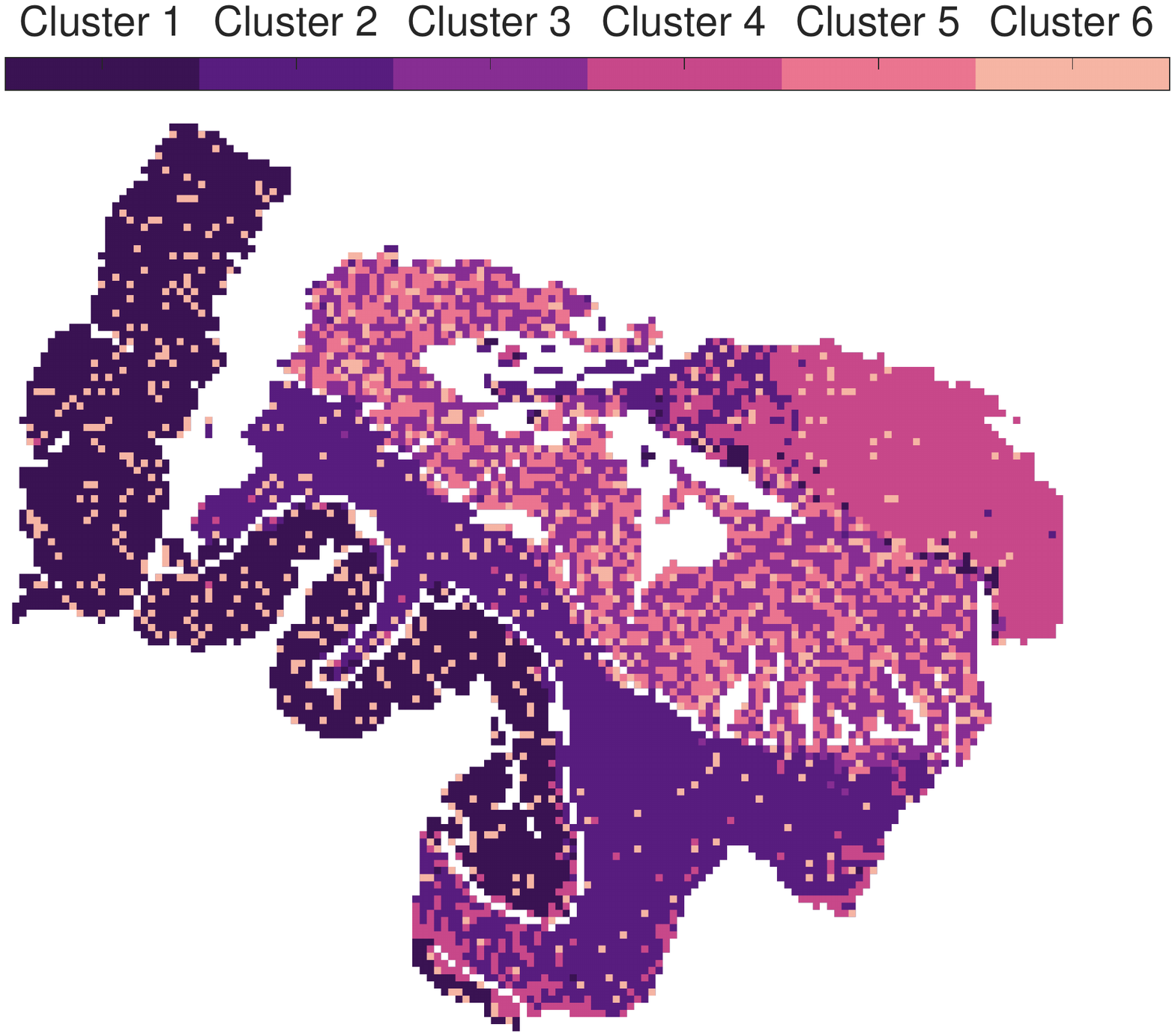}}\\
                \subfloat[\texttt{ONMF-TV-Pompili1}]{\label{fig:bestClusters-USep:ONMF-TV-Pompili1}\includegraphics[trim=5.31cm 2cm 4.27cm 1.1cm, clip, width=0.29\textwidth]{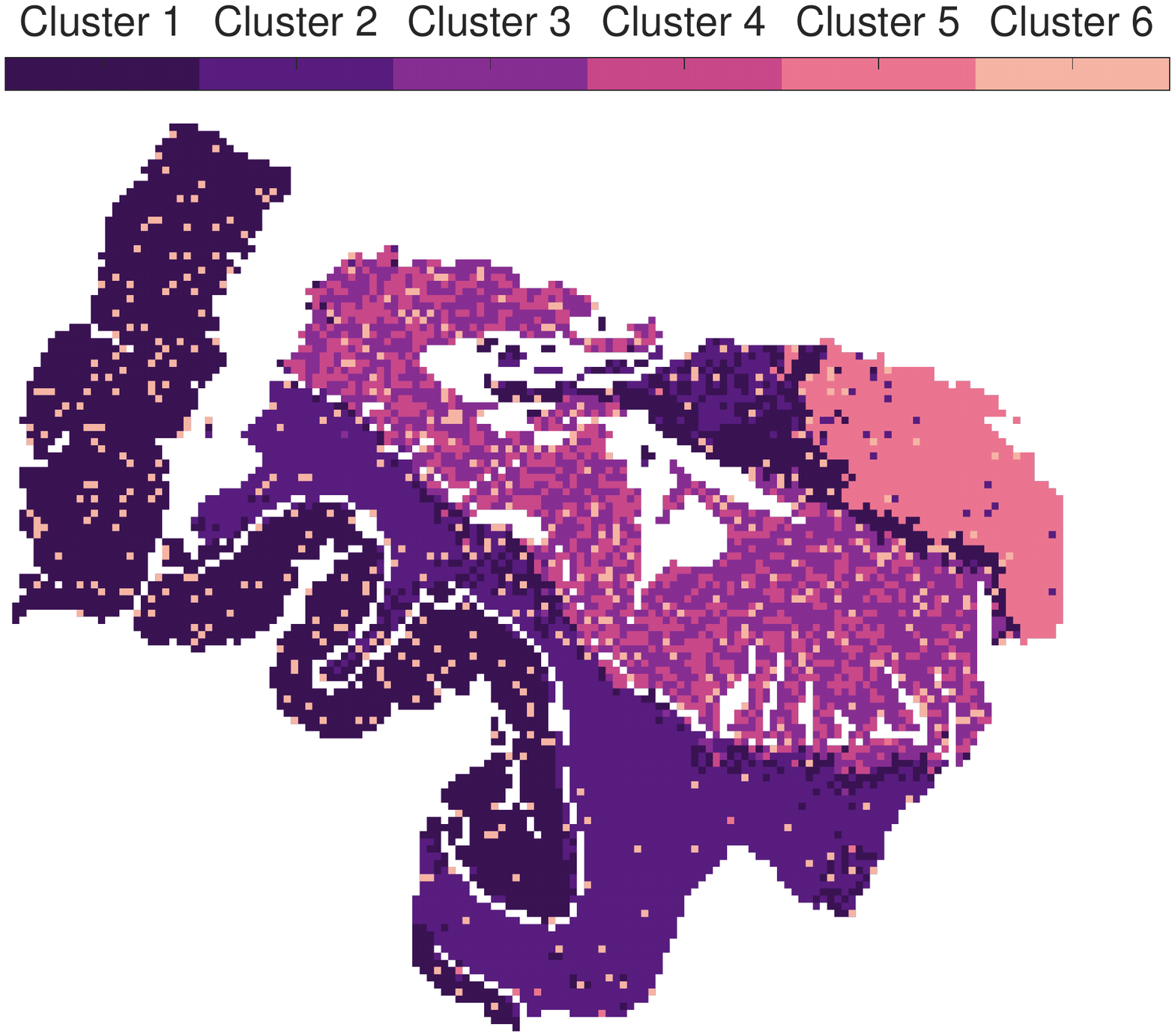}}\hfill
                \subfloat[\texttt{ONMF-TV-Pompili2}]{\label{fig:bestClusters-USep:ONMF-TV-Pompili2}\includegraphics[trim=5.31cm 2cm 4.27cm 1.1cm, clip, width=0.29\textwidth]{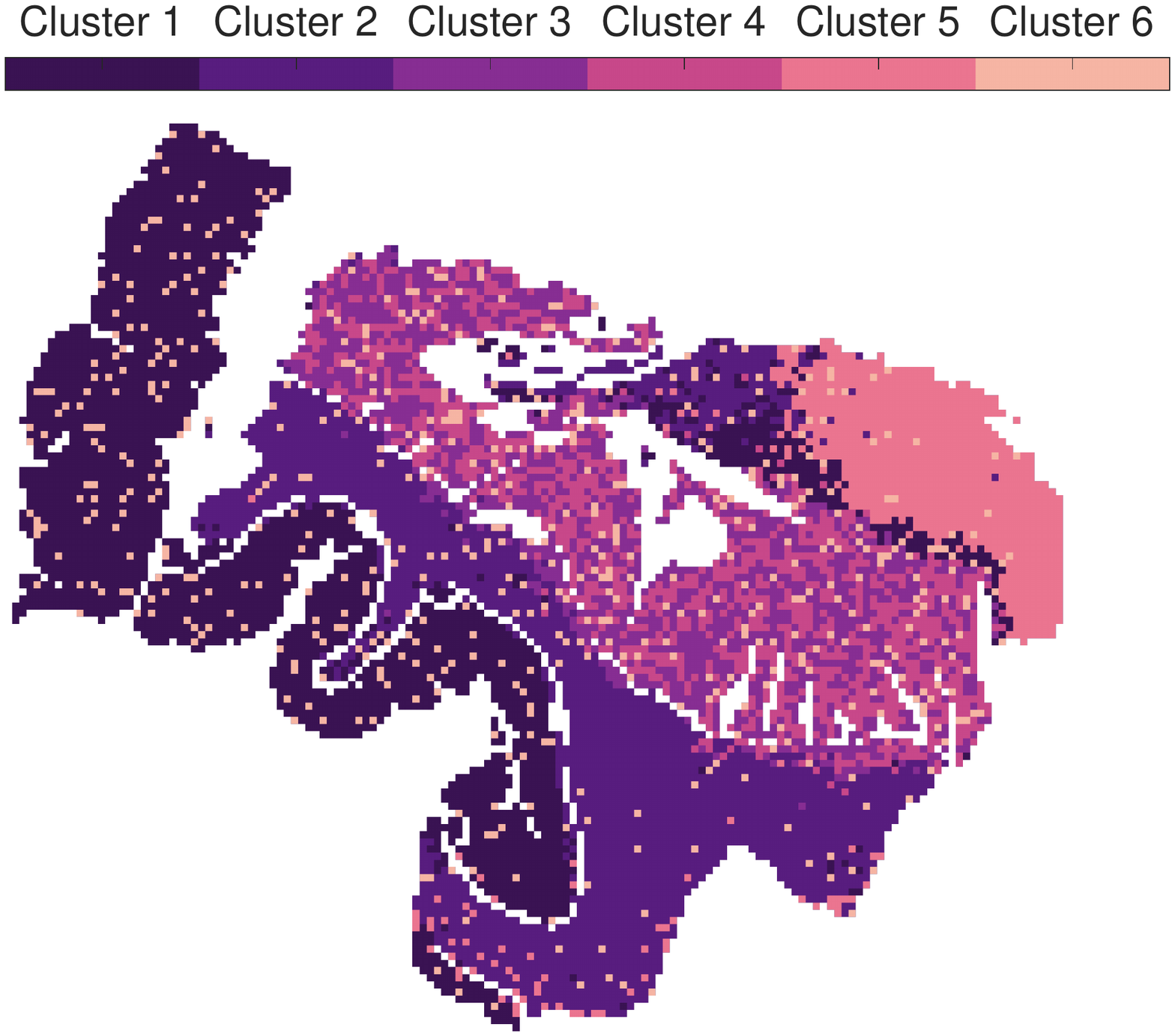}}\hfill
                \subfloat[\texttt{ONMF-TV-Kimura}]{\label{fig:bestClusters-USep:ONMF-TV-Kimura}\includegraphics[trim=5.31cm 2cm 4.27cm 1.1cm, clip, width=0.29\textwidth]{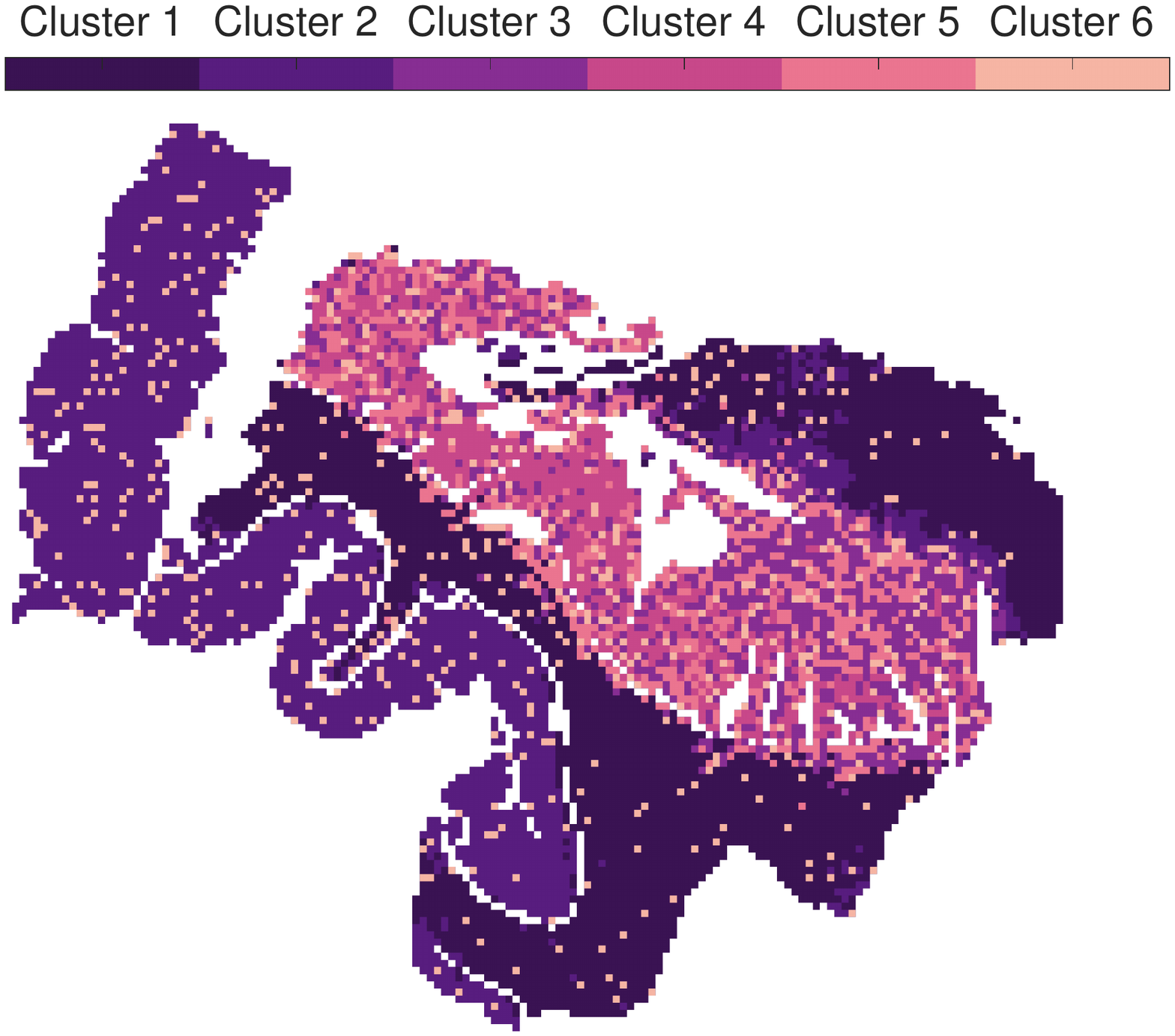}}\\
                \subfloat[\texttt{ONMF-TV-Li}]{\label{fig:bestClusters-USep:ONMF-TV-Li}\includegraphics[trim=5.31cm 2cm 4.27cm 1.1cm, clip, width=0.29\textwidth]{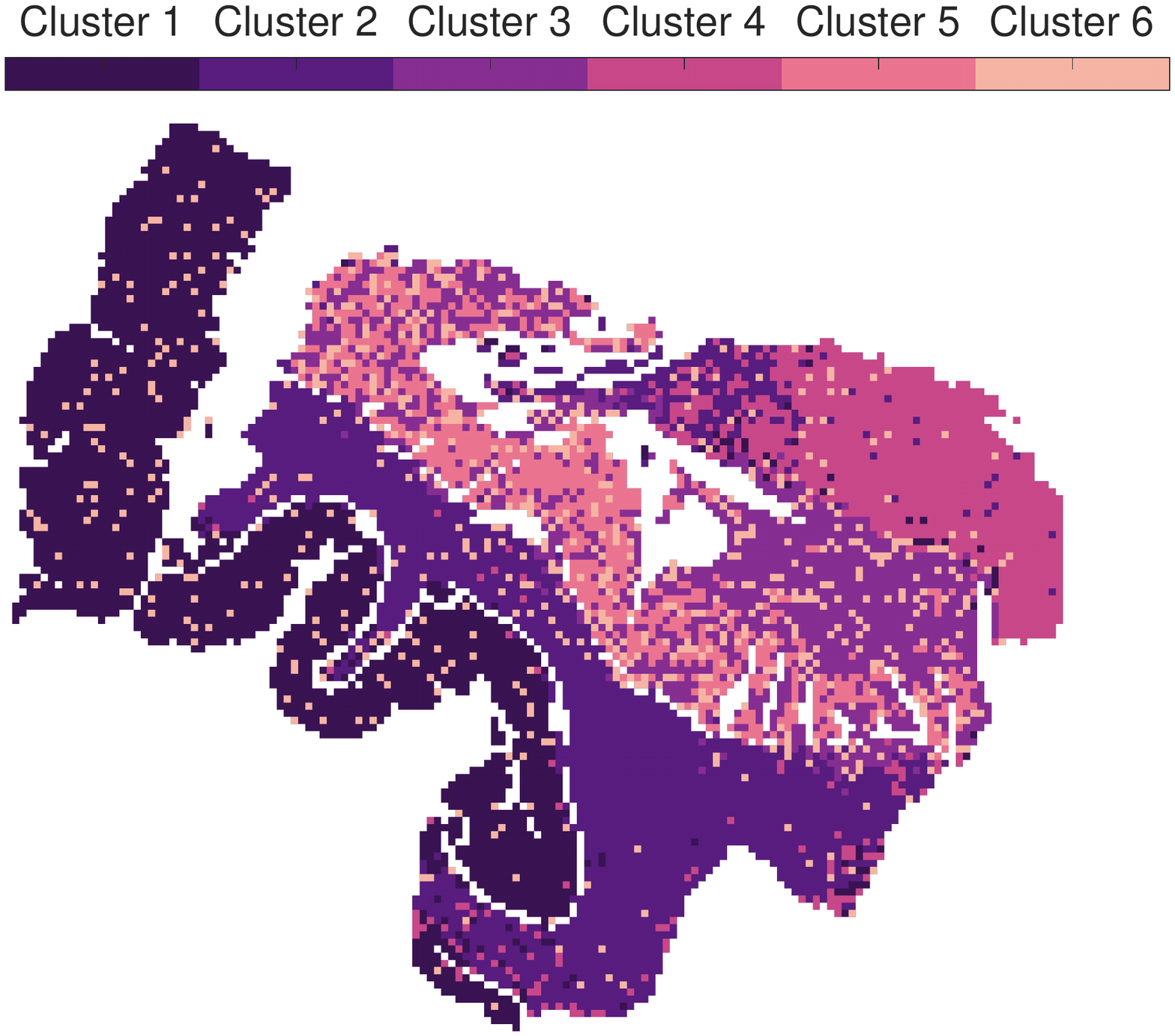}}\hfill
                \caption{Clusterings of all considered separated methods without TV regularization. The best-performing replicate including the TV regularization is chosen based on the normalized van Dongen criterion.}
                \label{fig:bestClusters-USep}
            \end{figure}
            \cref{fig:bestClusters-USep} shows the clusterings of the separated methods without the subsequent TV regularization step of the best performed replicate with respect to the VD\textsubscript{n}, which was obtained after the TV postprocessing. In general, every separated method is able to identify all classes shown in \cref{fig:dataSetAndGroundTruth:groundTruth} except the distinction between the tumor and mucosa. Since the combined methods are also not able to distinguish between both classes (see \cref{fig:bestClusters-UTV}), this is probably due to the fundamental, underlying NMF model with the orthogonality contraints, which are not able to identify the regions in an unsupervised workflow. In the case that annotated datasets are available for training, supervised machine learning methods could lead to different results. However, this work does not focus on these kind of approaches.
            
            Furthermore, every separated method is not able to guarantee a spatially coherent clustering in the region of the muscularis. Moreover, every method except \texttt{K-means-TV} does not provide any spatial coherence in any of the classes of the tissue slide. This is in contrast to the results in \cref{fig:bestClusters-UTV}, which shows the clusterings of all considered methods of the best performed replicate with respect to the VD\textsubscript{n} including the TV regularization. Every method is able to provide a spatially coherent clustering with some few exceptions in the region of the muscularis and hence leads to significantly improved clusterings in general. Comparing the results within \cref{fig:bestClusters-UTV}, the clustering of \texttt{ONMFTV-SPRING} seems to be the one, which best reproduces the given annotations. Furthermore, we note that some methods lead to clusterings with clusters, which do not contain any data points (see i.e.\ \texttt{K-means-TV}, \texttt{ONMF-TV-Pompili1}, \texttt{ONMFTV-PALM}). However, this is also owing to the fact that the class of the lymphocytes is significantly underrepresented compared to the other ones (see \cref{fig:dataSetAndGroundTruth:groundTruth}). Furthermore, this behaviour is also dependent on the choice of the TV regularization parameter $\tau.$
            
            \begin{figure}[tbhp]
                \centering
                \subfloat[Normalized van Dongen criterion (VD\textsubscript{n})]{\label{fig:boxPlots:VDN-VIN:VDN}\includegraphics[width=1.0\textwidth]{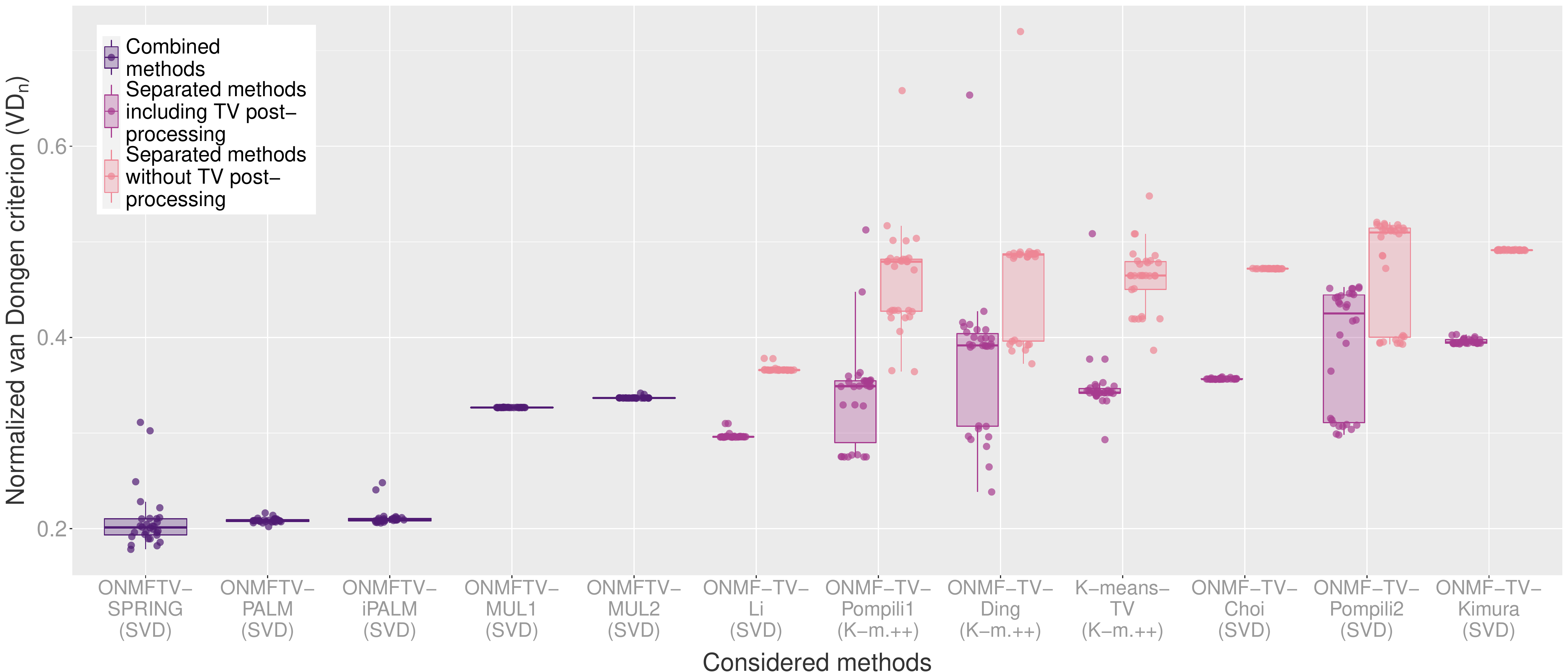}}\\
                \subfloat[Normalized variation of information (VI\textsubscript{n})]{\label{fig:boxPlots:VDN-VIN:VIN}\includegraphics[width=1.0\textwidth]{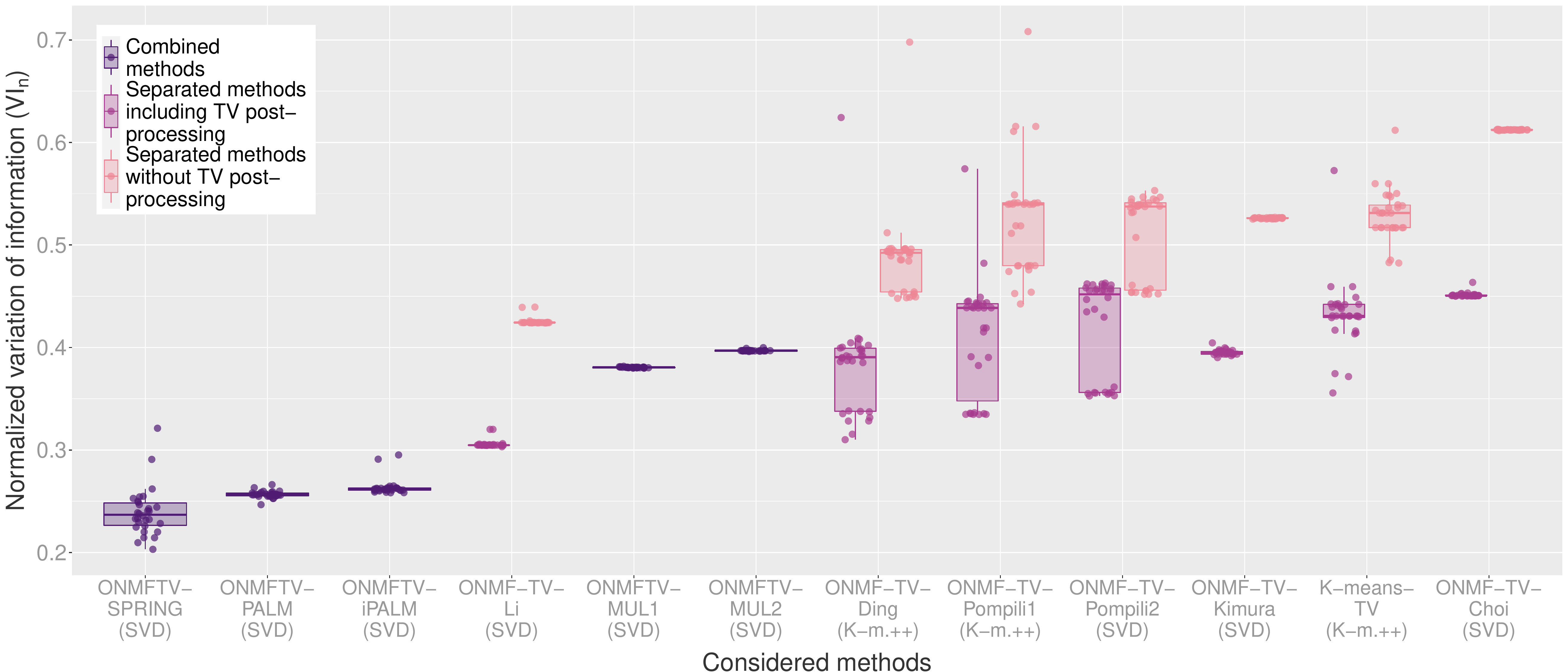}}
                \caption{Box plots of the normalized van Dongen criterion (VD\textsubscript{n}) and the normalized variation of information (VI\textsubscript{n}) of all performed experiments.}
                \label{fig:boxPlots:VDN-VIN}
            \end{figure}
            
            For a quantitative evaluation, we provide in \cref{fig:boxPlots:VDN-VIN} box plots of all replicates of every performed experiment and for all considered methods with respect to the VD\textsubscript{n} and VI\textsubscript{n}. Furthermore, regarding the separated methods, we also plot the validation measures for the clusterings with and without the TV postprocessing. For both measures, the observations of the qualitative evaluation above can be confirmed. First, we note that for all separated methods, the TV postprocessing does indeed lead to clusterings with better cluster validation measures. Moreover, the combined methods based on the PALM scheme achieve the best results with respect to both measures. While some experiments of \texttt{ONMFTV-SPRING} attain the best values compared to all other methods, this approach is less stable than the non stochastic approaches \texttt{ONMFTV-PALM} and \texttt{ONMFTV-iPALM}. Furthermore, we note that both combined methods \texttt{ONMFTV-MUL1} and \texttt{ONMFTV-MUL2} based on the multiplicative update rules do not perform better than some of the other separated methods. Comparing the separated methods with each other, we see that \texttt{ONMF-TV-Li} performs remarkably well with a high stability compared to the other approaches. Note that the stability also seems to depend on the initialization procedure. Regarding this, the SVD seems to favour more stable results than the K-means++ initialization.
            
            \begin{figure}[tbhp]
            	\centering
            	\includegraphics[width=1.0\textwidth]{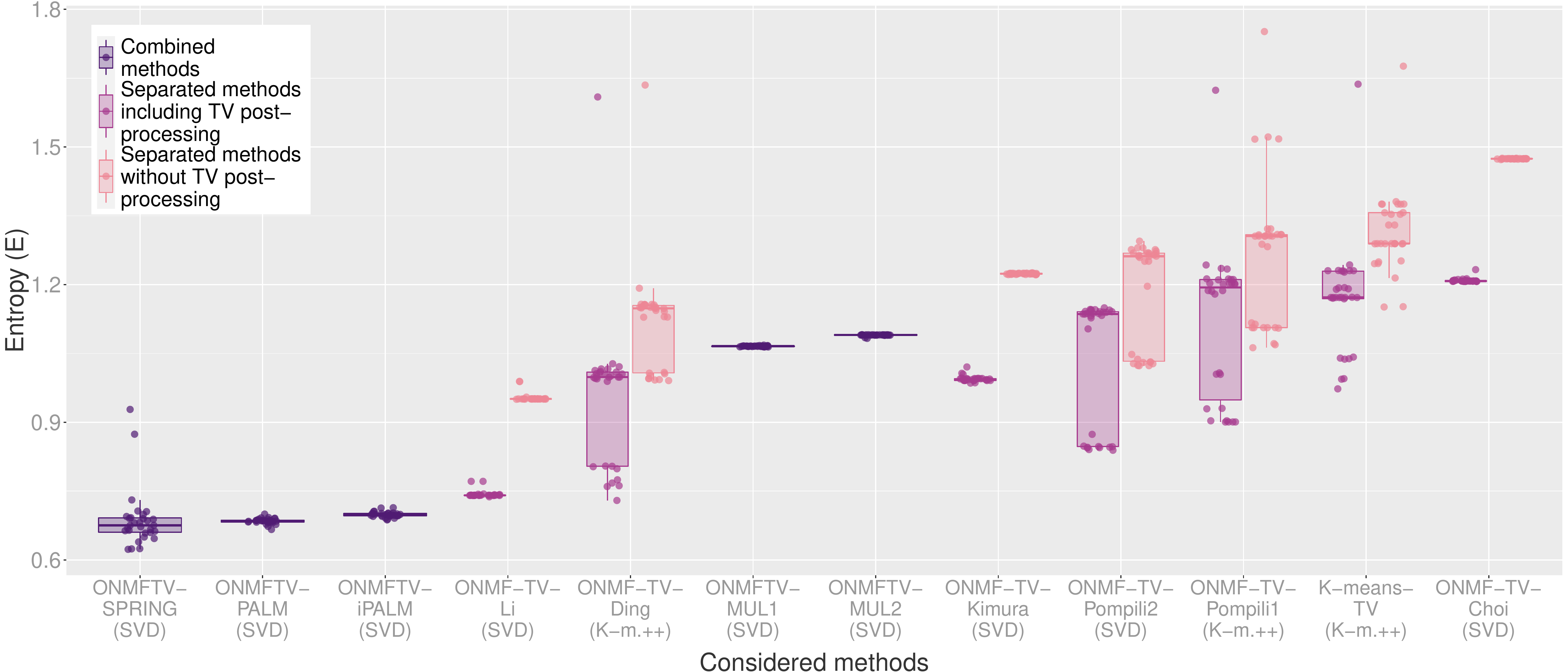}
            	\caption{Box plot of the Entropy (E) of all performed experiments.}
            	\label{fig:boxPlot:E}
            \end{figure}
            
            \cref{fig:boxPlot:E} shows the entropy of all performed experiments for all methods. This measure also confirms the observation, that the combined methods \texttt{ONMFTV-SPRING}, \texttt{ONMFTV-PALM} and \texttt{ONMFTV-iPALM} achieve the best results. Concerning the other methods, the outcomes are similar to the ones of the VD\textsubscript{n} and VI\textsubscript{n} and shall not be discussed in detail.
            
            \begin{figure}[tbhp]
            	\centering
            	\includegraphics[width=1.0\textwidth]{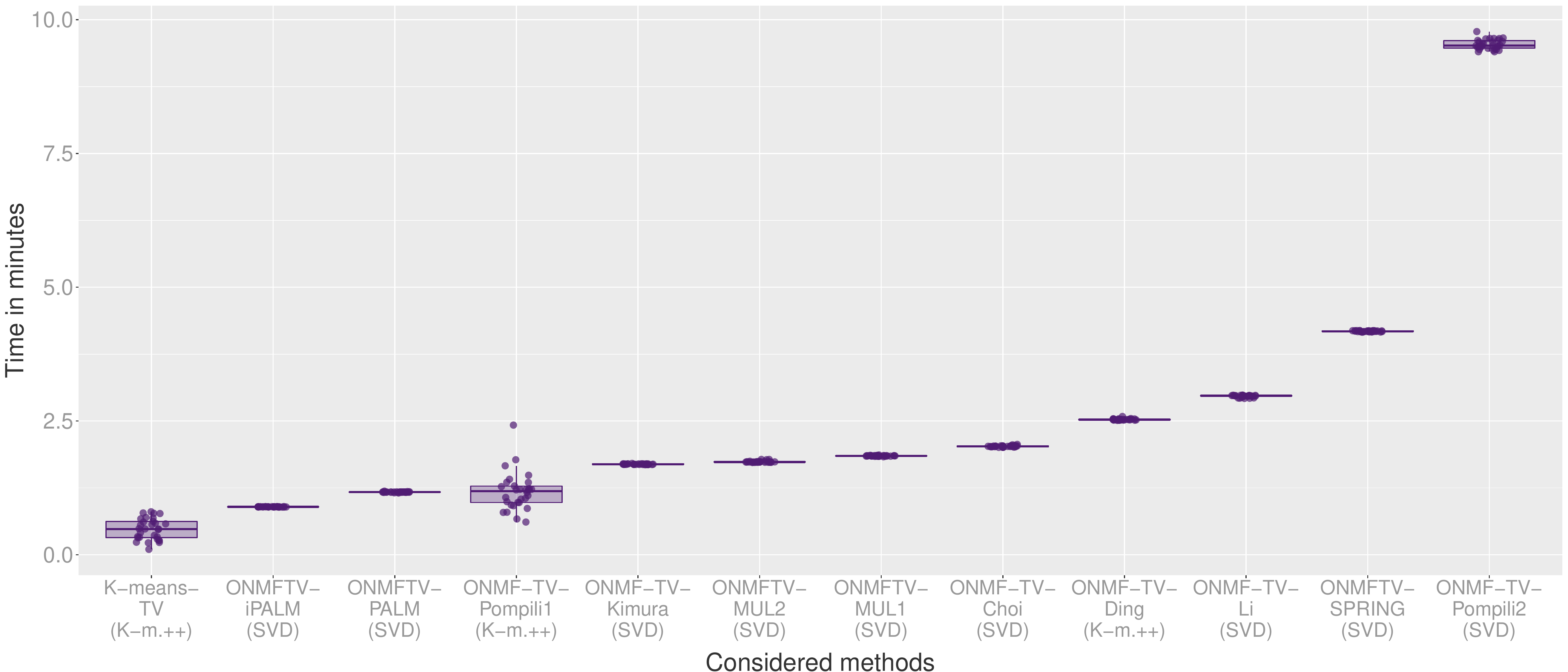}
            	\caption{Box plot of the computational times in minutes of all performed experiments.}
            	\label{fig:boxPlot:compTimes}
            \end{figure}
            
            \cref{fig:boxPlot:compTimes} shows the box plot of the computational times of all replicates for every method. The combined methods \texttt{ONMFTV-PALM} and \texttt{ONMFTV-iPALM} are one of the fastest methods together with \texttt{K-means-TV}, which requires the least time to compute the experiments. Furthermore, we note that \texttt{ONMFTV-SPRING} needs significantly more time compared to the other PALM algorithms. The other separated methods are faster than \texttt{ONMFTV-SPRING} except \texttt{ONMF-TV-Pompili2}, which needs significantly more time than every other considered approach.
            
            In \cref{app:sec:Further Experiments}, we show further quantitative results of the clusterings in terms of the VD and the VI. As expected, these results are very similar to the ones of the VD\textsubscript{n} and VI\textsubscript{n} and hence confirm our previous observations.
            
            All in all, based on the experiments performed on the MALDI dataset, we recommend the methods based on the PALM scheme and in particular \texttt{ONMFTV-PALM} as well as \texttt{ONMFTV-iPALM}, since they give the most stable results while achieving comparatively good results with low computational effort.
    
    \section{Conclusion} \label{sec: Conclusion}
        In this work, we have considered various orthogonal nonnegative matrix factorization (ONMF) models together with different optimization approaches for clustering hyperspectral data as the main application field. Furthermore, we have introduced total variation regularization in the proposed ONMF models to ensure spatial coherence in the obtained clusterings constituting the main innovation in this paper motivated by numerous spectral imaging applications, which naturally satisfy the spatial coherence in the data.
        
        After a brief description of the main principles of ONMF models, their relation to classical clustering methods and different optimization techniques, we have proposed so-called separated methods, which apply the TV denoising step after the computation of a classical ONMF algorithm. Furthermore, we have introduced more sophisticated combined methods with different optimization procedures, which include the TV regularization into the considered ONMF model.
        
        For the numerical evaluation, we have compared 12 different TV regularized ONMF methods on a MALDI-MSI human colon hyperspectral dataset with six different spatially coherent tissue regions, which constitute the ground truth for the clustering problem. The qualitative and quantitative results confirmed our expectation that the TV regularization significantly improves the clustering performance. Furthermore, the combined methods based on the proximal alternating linearized minimization have led to the best clustering outcomes and performance stability. Hence, based on the numerical evaluation of the MALDI dataset, we recommend the methods \texttt{ONMFTV-PALM}, \texttt{ONMFTV-iPALM} as well as \texttt{ONMFTV-SPRING}.
        
        Several further research directions could be of interest. One limitation of the presented approaches is the need of a manual a-priori choice of the needed hyperparameters. Hence, a useful extension of the proposed methods could be to introduce an automated way to choose suitable parameters. Another aspect is the analysis and the derivation of optimization algorithms for the case of discrepancy terms different from the Frobenius norm. Moreover, further gradient estimators different from the SGD could be examined for the method \texttt{ONMFTV-SPRING}. Furthermore, another major point is the consideration of more hyperspectral datasets from different application fields and a more thorough numerical evaluation of the different ONMF methods.
        
        Finally, a more theoretical aspect could be the investigation of spatially coherent clustering models in infinite dimension space leading to ``continuous'' factorization problems with gradient based penalty terms. In this setting, the analysis of first order conditions could lead to connections to corresponding K-means clustering models and partial differential equations, whose solutions give insight to the according distance measures and clusters. A first step for such an investigation could be to start with a finite dimensional space based on ONMF models.
    \appendix
    
    \section{Derivation of the Algorithm \texttt{ONMFTV-MUL1}} \label{app:sec:Derivation of the Algorithm for ONMFTV-MUL1}
        In this section, we give an outline of the proof of \cref{thm:ONMFTV-MUL1}. For more details on the derivation, we refer the reader to the work \cite{Fernsel:2018:Survey}.
        
        The update rules in \cref{alg:ONMFTV-MUL1} are based on the majorize-minimization (MM) principle \cite{Lange:13:Optimization}. The basic idea behind this concept is to replace the cost function $\mathcal{F}$ by a so-called surrogate function $\mathcal{Q}_\mathcal{F}$, whose minimization is simplified and leads to the desired multiplicative algorithms.
        \begin{definition}
            Let $\Omega \subset \R^{M\times K}$ be an open set and $\mathcal{F}: \Omega \to \R$ a given cost function. A function $\mathcal{Q}_\mathcal{F}: \Omega \times \Omega \to \R$ is called a surrogate function, if it satisfies the following properties:
            \begin{itemize}
                \item[(i)] $ \mathcal{Q}_\mathcal{F}(U, A) \geq \mathcal{F}(U)$ for all $ U, A \in \Omega, $ \label{itm:def:surrogate property 1}
        		\item[(ii)] $ \mathcal{Q}_\mathcal{F}(U, U) = \mathcal{F}(U) $ for all $ U\in \Omega. $ \label{itm:def:surrogate property 2}
            \end{itemize}
        \end{definition}
        The minimization procedure based on the MM approach is defined by
        \begin{equation*}
            U^{[i+1]} \coloneqq \argmin_{U\in \Omega} \mathcal{Q}_\mathcal{F}(U, U^{[i]}).
        \end{equation*}
        Together with the properties of $\mathcal{Q}_\mathcal{F},$ this leads directly to the monotone decrease of the cost function $\mathcal{F},$ since
        \begin{equation*}
            \mathcal{F}(U^{[i+1]}) \leq \mathcal{Q}_\mathcal{F}(U^{[i+1]}, U^{[i]}) \leq \mathcal{Q}_\mathcal{F}(U^{[i]}, U^{[i]}) = \mathcal{F}(U^{[i]}).
        \end{equation*}
        Tailored construction techniques for surrogate functions leads additionally to the desired multiplicative structure of the update rules and are typically based on Jensen's inequality or the so-called quadratic upper bound principle \cite{Lange:13:Optimization,Fernsel:2018:Survey}.
        
        We first focus on the minimization with respect to $U.$ It can be shown that for each term of the cost function $\mathcal{F}$ in \cref{eq:CombinedMethods:MUL:ONMFTV-MUL1:NMFProblem2}, a surrogate function can be constructed, which finally results to a suitable surrogate $\mathcal{Q}_\mathcal{F} \coloneqq \mathcal{Q}_{\mathcal{F}_1} + \mathcal{Q}_{\mathcal{F}_2} + \mathcal{Q}_{\mathcal{F}_3}$ defined by
        \begin{align*}
            \mathcal{Q}_{\mathcal{F}_1}(U, A) &\coloneqq \dfrac{1}{2} \sum_{m=1}^{M} \sum_{n=1}^{N} \dfrac{1}{(AV)_{mn}} \sum_{k=1}^{K} A_{ik}V_{kn} \left( X_{mn} - \dfrac{U_{mk}}{A_{mk}}(AV)_{mn}\right)^2,\\
            \mathcal{Q}_{\mathcal{F}_2}(U, A)  &\coloneqq\dfrac{\sigma_1}{2} \sum_{k=1}^{K} \sum_{\ell=1}^{K} \dfrac{1}{(W^\intercal A)_{k \ell}} \sum_{m=1}^{M} W_{m k} A_{m \ell} \left ( \delta_{k \ell} - \dfrac{U_{m \ell}}{A_{m \ell}} (W^\intercal A)_{k \ell} \right )^2\! + \dfrac{\sigma_2}{2} \Vert W - U\Vert_F^2, \\
            \mathcal{Q}_{\mathcal{F}_3}(U, A) &\coloneqq \frac{\tau}{2} \left( \sum_{k=1}^K \sum_{m=1}^M \left [ P(A)_{mk} (U_{mk} - Z(A)_{mk})^2  \right ] + C(A) \right),
        \end{align*}
        where $C(A)$ is some function depending only on $A$ and with $P(A)$ and $Z(A)$ given as in \cref{eq:ONMFTV-MUL1-P,eq:ONMFTV-MUL1-Z}. Computing the first order condition $\nabla_U \mathcal{Q}_\mathcal{F}(U,A) = 0$ and applying classical calculation rules for derivatives leads then to the desired update rule given in \cref{alg:ONMFTV-MUL1-Line4} of \cref{alg:ONMFTV-MUL1}.
        
        The update rules for $V$ and $W$ are treated similarly.
    \section{Details on the Algorithm \texttt{ONMFTV-MUL2}} \label{app:sec:Details on the Algorithm ONMFTV-MUL2}
        \subsection{Derivation of the Update Rules}\label{app:subsec:Derivation of the Update Rules}
            The \cref{alg:ONMFTV-MUL2} is based on a classical gradient descent approach with a suitably chosen step size to ensure a multiplicative structure of the update rules. We will discuss here only the derivation for the minimization with respect to $U.$ The update rules for $V$ are classical results and can be found in various works \cite{Lee:2001:NMFMultiplicative,Fernsel:2018:Survey}.
            
            Regarding the update rule for $U,$ we consider the classical gradient descent step
            \begin{equation*}
                U^{[i+1]} \coloneqq U^{[i]} - \Gamma^{[i]} \circ \nabla_U \mathcal{F}(U^{[i]},V)
            \end{equation*}
            given also in \cref{eq:ONMFTV-MUL2-GradientDescentStep}. Thus, we need an explicit expression for $\nabla_U \mathcal{F}(U,V).$ Classical calculation rules for derivatives yields
            \begin{equation*}
                \nabla_U \mathcal{F}(U,V) = UVV^\intercal - XV^\intercal + \sigma_1 (UU^\intercal U - U) + \nabla_U \TVeps(U).
            \end{equation*}
            The gradient of $\TVeps(U)$ can be acquired via the Euler-Lagrange equation and is a classical result. By considering the continuous case with the function $u:\Omega \to \R$ of \cref{par:ONMFTV-MUL2} and defining $\mathcal{L}(x_1, x_2, u, u_1, u_2) \coloneqq \Vert \nabla u \Vert_{\varepsilon_{\TV}}$ with $u_i \coloneqq \partial u / \partial x_i$ for $i\in \{ 1, 2\},$ the Euler Lagrange Equation
            \begin{equation*}
                \frac{\partial \mathcal{L}}{\partial u} - \sum_{i=1}^2 \frac{\partial}{\partial x_i} \left( \frac{\partial \mathcal{L}}{\partial u_i} \right) = 0
            \end{equation*}
            gives the formal relationship 
            \begin{equation*}
                \nabla_u \TVeps(u) \coloneqq -\divergence \left( \frac{\nabla u}{\Vert \nabla u \Vert_{\varepsilon_{\TV}}} \right).
            \end{equation*}
            Thus, the gradient descent step is given by
            \begin{equation*}
                U^{[i+1]} \coloneqq U^{[i]} - \Gamma^{[i]} \circ \left( U^{[i]} VV^\intercal - XV^\intercal - \tau \divergence \left( \frac{\nabla U^{[i]}}{\Vert \nabla U^{[i]} \Vert_{\varepsilon_{\TV}}} \right) + \sigma_1 U^{[i]} {U^{[i]}}^\intercal U^{[i]} - \sigma_1 U^{[i]} \right).
            \end{equation*}
            To ensure a multiplicative structure of the update rules, we set the step size to be
            \begin{equation*}
                \Gamma^{[i]} \coloneqq \frac{U^{[i]}}{U^{[i]}VV^\intercal + \sigma_1 U^{[i]}{U^{[i]}}^\intercal U^{[i]}},
            \end{equation*}
            which leads directly to the update rule in \cref{alg:ONMFTV-MUL2_Line4} of \cref{alg:ONMFTV-MUL2}.
        \subsection{Discretization of the TV Gradient}\label{app:subsec:Discretization of the TV Gradient}
            In this section, we describe the discretization procedure of the divergence term $\divergence\left( \nabla U^{[i]} / \Vert \nabla U^{[i]} \Vert_{\varepsilon_{\TV}} \right),$ which occurs in \cref{alg:ONMFTV-MUL2_Line4} of \cref{alg:ONMFTV-MUL2}.
            
            To perform such a discretization, it is needed for express the divergence term in terms of sums and products of first and second order derivatives. To simplify the notation, we stick in this first step to the continuous case and consider the function $u$ mentioned in \cref{eq:ONMFTV-MUL2-TV Penalty Term}. By considering the definition in \cref{eq:ONMFTV-MUL2-norm of nabla u} and applying classical calculation rules for derivatives, we get
            \begin{align*}
                 \divergence\left( \frac{\nabla u}{ \Vert \nabla u \Vert_{\varepsilon_{\TV}} } \right) = \frac{\varepsilon_{\TV}^2 \left( \frac{\partial^2 u}{\partial x^2} + \frac{\partial^2 u}{\partial y^2} \right) + \left( \frac{\partial u}{\partial x}\right)^2 \frac{\partial^2 u}{\partial y^2} + \left( \frac{\partial u}{\partial y}\right)^2 \frac{\partial^2 u}{\partial x^2} - 2 \frac{\partial u}{\partial x} \frac{\partial u}{\partial y} \frac{\partial^2 u}{\partial x \partial y} }{ \Vert \nabla u\Vert_{\varepsilon_{\TV}}^3}.
            \end{align*}
            For the discretization of this expression, we assume that $u$ is a given discretized image of size $d_1 \times d_2.$ Note that in the NMF models, which are considered in this work, this would correspond to one (reshaped) column of the matrix $U.$ To approximate the partial derivatives in the above expression, we use in the following a central differencing scheme and interpret the $x$ and $y$ directions as the vertical and horizontal axis in the image $u$ respectively. Thus, we define
            \begin{align*}
                (\Delta_x u)_{ij}       &\coloneqq \frac{u_{i+1,j} - u_{i-1,j}}{2},     & (\Delta_y u)_{ij}         &\coloneqq \frac{u_{i,j+1} - u_{i,j-1}}{2},\\
                (\Delta_{xx} u)_{ij}    &\coloneqq u_{i+1,j} - 2 u_{ij} + u_{i-1,j},    & (\Delta_{yy} u)_{ij}      &\coloneqq u_{i,j+1} - 2u_{ij} + u_{i,j-1},\\
                (\Delta_{xy} u)_{ij}    &\coloneqq \frac{u_{i+1,j+1} - u_{i+1,j-1} - u_{i-1,j+1} + u_{i-1,j-1}}{4}.
            \end{align*}
            The discretized gradients can also be interpreted as matrices of size $d_1 \times d_2,$ which finally leads to the discretization of the divergence term
            \begin{equation*}
                \scalebox{0.9}{\text{$\displaystyle \divergence\left( \frac{\nabla u}{ \Vert \nabla u \Vert_{\varepsilon_{\TV}} } \right) = \frac{\varepsilon_{\TV}^2 \left( \Delta_{xx} u + \Delta_{yy} u \right) +  (\Delta_x u)^2 \circ \Delta_{yy}u + (\Delta_y u)^2 \circ \Delta_{xx}u - 2 \Delta_x u \circ \Delta_y u \circ \Delta_{xy} u }{ \left( (\Delta_x u)^2 + (\Delta_y u)^2 + \varepsilon_{\TV}^2 1_{d_1\times d_2} \right)^{\nicefrac{3}{2}} }, $ } }
            \end{equation*}
            where $1_{d_1\times d_2}$ is a matrix of size $d_1\times d_2$ with ones in every entry.
    \section{Algorithmic Details on the Proximal Gradient Descent Approach} \label{app:sec:Algorithmic Details on the Proximal Gradient Descent Approach}
            In this section, we give some information about the involved gradients, Lipschitz constants and the computation of the step sizes of all algorithms based on the proximal gradient descent approach (see \cref{subsubsec: Promximal Alternating Linearized Minimization}).
            \subsection{\texttt{ONMFTV-PALM}} \label{app:subsec: ONMFTV-PALM}
                We start in this section with details on the gradients of the algorithm \texttt{ONMFTV-PALM}. Since the computations are straightforward, we will only state the final results of all three gradients:
                \begin{align*}
                    \nabla_U \mathcal{F}(U,V,W) &= UVV^\intercal - XV^\intercal + \sigma_1 (WW^\intercal U - W) + \sigma_2(U-W),\\
                    \nabla_V \mathcal{F}(U,V,W) &= U^\intercal UV - U^\intercal X,\\
                    \nabla_W \mathcal{F}(U,V,W) &= \sigma_1 (UU^\intercal W - U) + \sigma_2 (W-U).
                \end{align*}
                For the calculation of the Lipschitz constants, we compute for the minimization with respect to $U$ and arbitrary matrices $U_1, U_2 \in \mathbb{R}^{M\times K}$
                \begin{align*}
                    \Vert \nabla_U \mathcal{F}(U_1,V,W) - \nabla_U \mathcal{F}(U_2,V,W) \Vert &= \Vert (U_1 - U_2)VV^\intercal + \sigma_1 WW^\intercal (U_1 - U_2) + \sigma_2 (U_1 - U_2)\Vert \\
                    &\leq (\Vert VV^\intercal\Vert + \Vert \sigma_1 WW^\intercal + \sigma_2 I_{M\times M} \Vert)\cdot \Vert U_1 - U_2\Vert\\
                    &=(\lambda_{1,U} + \lambda_{2,U})\cdot \Vert U_1 - U_2\Vert,
                \end{align*}
                where $\Vert \cdot \Vert$ is the spectral norm and $\lambda_{1,U}, \lambda_{2,U}$ the maximal absolute eigenvalue of the symmetric matrices $VV^\intercal$ and $\sigma_1 WW^\intercal + \sigma_2 I_{M\times M}$ respectively, with $I_{M\times M}$ being the identity matrix of size $M\times M.$ The other cases are treated analogously and results in
                \begin{equation*}
                    \Vert \nabla_V \mathcal{F}(U,V_1,W) - \nabla_V \mathcal{F}(U,V_2,W) \Vert \leq \Vert U^\intercal U\Vert \cdot \Vert V_1 - V_2\Vert = \lambda_{1,V} \Vert V_1 - V_2\Vert
                \end{equation*}
                as well as
                \begin{align*}
                    \Vert \nabla_W \mathcal{F}(U,V,W_1) - \nabla_W \mathcal{F}(U,V,W_2) \Vert &\leq \Vert \sigma_1 UU^\intercal + \sigma_2 I_{M\times M} \Vert \cdot \Vert W_1 - W_2\Vert\\
                    &= \lambda_{1,W} \Vert W_1 - W_2\Vert.
                \end{align*}
                Numerically, the eigenvalues of the matrices are approximated via the power iteration with 5 iterations. The step sizes $\eta_U, \eta_V$ and $\eta_W$ described in \cref{alg:ONMFTV-PALM} are computed based on the Lipschitz constants, which are given by $L_U \coloneqq \lambda_{1,U} + \lambda_{2,U},$ $L_V \coloneqq \lambda_{1,V}$ and $L_W \coloneqq \lambda_{1,W}$ with the computation above. The standard choice for the step size of the \texttt{ONMFTV-PALM} algorithm are simply the inverse Lipschitz constants, i.e.\ $\eta_U = 1/L_U, \eta_V = 1/L_V$ and $\eta_W = 1/L_W$ \cite{BolteEtal:2014:PALM,DriggsEtal:2020:SPRING}.
            \subsection{\texttt{ONMFTV-iPALM}} \label{app:subsec: ONMFTV-iPALM}
                The calculation of the gradients and Lipschitz constants for \texttt{ONMFTV-\\iPALM} are the same as in the case of \texttt{ONMFTV-PALM} with the slight difference that the gradients are evaluated at other points. Furthermore, the step sizes are set by $\eta_U = 0.9/L_U, \eta_V = 0.9/L_V$ and $\eta_W = 0.9/L_W$ according to \cite{DriggsEtal:2020:SPRING}. For details, we refer the reader to \cite{PockSabach:2016:iPALM,DriggsEtal:2020:SPRING}.
            \subsection{\texttt{ONMFTV-SPRING}} \label{app:subsec:ONMFTV-SPRING}
                The computation of the gradients for \texttt{ONMFTV-SPRING} are based on the SGD estimator given in \cref{eq:ONMFTV-SPRING-SGDEstimator}. As in the case of \texttt{ONMFTV-PALM} the computations are straightforward and we get by defining the mini-batch $\mathcal{B}_{i,j}^V \subset \{ 1,\dots,M \}$ for the minimization with respect to $V$ and with $\mathcal{F}_m(U,V,W) \coloneqq \nicefrac{1}{2}\Vert X_{m,\bullet} - U_{m,\bullet}V \Vert_2^2, $ we compute the gradients
                \begin{align*}
                    \tilde{\nabla}_U^{i,j} \mathcal{F}(U,V,W) &= \sum_{n\in \mathcal{B}_{i,j}^U} \left[ (UV)_{\bullet, n} V^\intercal_{n,\bullet} - X_{\bullet,n} V^\intercal_{n,\bullet} + \nicefrac{1}{N} \left(\sigma_1 (WW^\intercal U - W) + \sigma_2 (U-W) \right) \right]\\
                    &=UV_{\bullet, \mathcal{B}_{i,j}^U}V^\intercal_{\mathcal{B}_{i,j}^U, \bullet} - X_{\bullet, \mathcal{B}_{i,j}^U} V^\intercal_{\mathcal{B}_{i,j}^U, \bullet} + \frac{\vert \mathcal{B}_{i,j}^U\vert}{N} \left(\sigma_1 (WW^\intercal U - W) + \sigma_2 (U-W) \right),\\
                    \tilde{\nabla}_V^{i,j} \mathcal{F}(U,V,W) &= \sum_{m\in \mathcal{B}_{i,j}^V} \nabla_V \mathcal{F}_m(U,V,W) = \sum_{m\in \mathcal{B}_{i,j}^V} U^\intercal_{\bullet, m} (UV)_{m,\bullet} - U^\intercal_{\bullet, m}X_{m,\bullet} \\
                    &\hspace{25.9ex}=U^\intercal_{\bullet, \mathcal{B}_{i,j}^V} U_{\mathcal{B}_{i,j}^V, \bullet} V - U^\intercal_{\bullet, \mathcal{B}_{i,j}^V} X_{\mathcal{B}_{i,j}^V, \bullet},
                \end{align*}
                where $V_{\bullet, \mathcal{B}_{i,j}^U}, X_{\bullet, \mathcal{B}_{i,j}^U}$ and $U_{\mathcal{B}_{i,j}^V, \bullet}, X_{\mathcal{B}_{i,j}^V, \bullet}$ are the submatrices of $U, V$ and $X,$ which are constrained column-wise and row-wise based on the index sets $\mathcal{B}_{i,j}^U$ and $\mathcal{B}_{i,j}^V$ respectively. For the minimization with respect to $W,$ the full partial gradient of $\mathcal{F}$ is used which was already computed in \cref{app:subsec: ONMFTV-PALM}.
                
                The computation of the Lipschitz constants of the partial gradients of $\mathcal{F}$ with respect to $U$ and $V$ is based on the SGD estimator and goes analogously to the steps in \cref{app:subsec: ONMFTV-PALM}. Hence, we have
                
                \begin{align*}
                    \Vert \tilde{\nabla}_U^{i,j} \mathcal{F}(U_1,V,W) - \tilde{\nabla}_U^{i,j} \mathcal{F}(U_2,V,W) \Vert &\leq \scalebox{0.87}{\text{$\displaystyle \left( \Vert V_{\bullet, \mathcal{B}_{i,j}^U}V^\intercal_{\mathcal{B}_{i,j}^U, \bullet} \Vert + \frac{\vert \mathcal{B}_{i,j}^U\vert}{N} \Vert \sigma_1 WW^\intercal + \sigma_2 I_{M\times M} \Vert \right) \Vert U_1 - U_2$ }} \Vert \\
                    &= (\lambda_{1,U} + \lambda_{2,U})\cdot \Vert U_1 - U_2\Vert,\\
                    \Vert \tilde{\nabla}_V^{i,j} \mathcal{F}(U,V_1,W) - \tilde{\nabla}_V^{i,j} \mathcal{F}(U,V_2,W) \Vert &\leq \Vert U^\intercal_{\bullet, \mathcal{B}_{i,j}^V} U_{\mathcal{B}_{i,j}^V, \bullet} \Vert \cdot \Vert V_1 - V_2\Vert = \lambda_{1,V} \Vert V_1 - V_2\Vert
                \end{align*}
                together with the Lipschitz constants $L_U, L_V$ analogously to \cref{app:subsec: ONMFTV-PALM}. Since we consider the full partial gradient of $\mathcal{F}$ with respect to $W,$ we also get the same Lipschitz constant $L_W$ as in \cref{app:subsec: ONMFTV-PALM}.
                
                The choice of of the step sizes $\eta_{U^{[i,j]}}$ and $\eta_{V^{[i,j]}}$ are chosen according to the work \cite{DriggsEtal:2020:SPRING} by defining
                \begin{align*}
                    \eta_{U^{[i,j]}} \coloneqq \min \left\{ \frac{1}{ \sqrt{ \left\lceil i\cdot \nicefrac{\vert \mathcal{B}_{i,j}^U \vert}{N} \right\rceil }\cdot L_U },\ \frac{1}{L_U} \right\}, && \eta_{V^{[i,j]}} \coloneqq \min \left\{ \frac{1}{ \sqrt{ \left\lceil i\cdot \nicefrac{\vert \mathcal{B}_{i,j}^V \vert}{M} \right\rceil }\cdot L_V },\ \frac{1}{L_V} \right\}.
                \end{align*}
                Furthermore, as described in \cref{subsec: Setup}, we perform an additional projection step of the parameter $\tau \eta_{U^{[i]}}$ (see \cref{alg:ONMFTV-SPRING_Line6} in \cref{alg:ONMFTV-SPRING} and the application of the $\prox_{\tau \eta_{U^{[i,j]}} \mathcal{J}}$ operator) to avoid too large regularization parameters of the TV penalty term by applying $\prox_{\tau_{i,j} \mathcal{J}}$ with $\tau_{i,j} \coloneqq \min \{ \tau \eta_{U^{[i,j]}},\ 1\cdot 10^{-3} \}.$
    \section{Parameter Choice} \label{app:sec:Parameter Choice}
        As described in \cref{subsec: Setup}, we choose the regularization parameters for the numerical experiments of every considered method empirically by performing multiple experiments for different parameters and choosing the ones which leads to stable experiments and the best VD\textsubscript{n} and VI\textsubscript{n}. For the separated methods, we also follow partially the recommendations of the corresponding works. In the following sections, we describe the choice of these and other hyperparameters in more detail.
        
        \subsection{Separated Methods} \label{app:subsec:Parameter Choice of Separated Methods}
            In this Section, we describe the choice of the main hyperparameters of all separated methods. \cref{tab:parameterChoice:SeparatedMethods} shows the selected parameters for the numerical experiments in \cref{sec: Numerical Experiments} along with the used stopping criteria (stopCrit) and initialization methods (initMethod).
        
            \begin{table}[tbhp]
                    \caption{Parameter choice of the separated methods for the numerical experiments with the MALDI dataset.}
                    \ra{1.1}
                    \centering
                    \resizebox{0.85\columnwidth}{!}{%
                    \pgfplotstabletypeset[
                        col sep=&, row sep=\\,
                        string type,
                        every head row/.style={
                            before row={
                                \toprule
                            },
                            after row=\midrule,
                        },
                        every last row/.style={
                            after row=\bottomrule},
                        every even row/.style={
                            before row={\rowcolor[gray]{0.9}}
                        },
                        columns/method/.style       ={column name=\textbf{Method}, column type = l},
                        columns/tau/.style          ={column name=$\tau$},
                        columns/stopCrit/.style     ={column name=stopCrit, column type = l},
                        columns/iMax/.style         ={column name=$i_{\text{max}}$},
                        columns/initMethod/.style   ={column name=initMethod, column type = l},
                        ]{
                            method                      & tau                   & stopCrit          & iMax                  & initMethod \\
                			\texttt{K-means-TV}            & $1$                   & Cluster assignment     & --               & K-means++\\
                			\texttt{ONMF-TV-Choi}       & $2\cdot 10^{-2}$      & maxIt             & $6\cdot 10^{2}$       & SVD\\
                			\texttt{ONMF-TV-Ding}       & $2\cdot 10^{-2}$      & maxIt             & $8\cdot 10^{2}$       & K-means++\\
                			\texttt{ONMF-TV-Pompili1}   & $1$                   & Cluster assignment     & --               & K-means++\\
                			\texttt{ONMF-TV-Pompili2}   & $4\cdot 10^{-2}$      & Pompili Internal  & --                    & SVD\\
                			\texttt{ONMF-TV-Kimura}     & $2\cdot 10^{-2}$      & maxIt             & $7\cdot 10^{2}$       & SVD\\
                			\texttt{ONMF-TV-Li}         & $3\cdot 10^{-2}$      & maxIt             & $2\cdot 10^{2}$       & SVD\\
                        }
                    }
            		\label{tab:parameterChoice:SeparatedMethods}
                \end{table}
                All separated methods are initialized either based on the SVD approach described in \cref{subsec: Setup} or via K-means++. Both methods were tried out for every separated method in the numerical experiments and one of them was chosen based on the stability and quality measures of the results. For the specific case of \texttt{ONMF-TV-Pompili1}, we use a mixture of the K-means++ method and the internal initialization procedure of \cite{Pompili:2014:CompMethod}.
                
                Furthermore, different stopping criteria were used. Besides the classical stopping criterion based on a maximal iteration number $i_{\text{max}}$ (maxIt), the clustering algorithms \texttt{K-means-TV} and \texttt{ONMF-TV-Pompili1} (i.e.\ \cref{alg:ONMF-TV_Line3} in \cref{alg:ONMF-TV}) are stopped until the cluster assignments in the cluster membership matrix $U$ do not change anymore. For the special case of \texttt{ONMF-TV-Pompili2}, the algorithm stops until the current iterates are ``sufficiently'' nonnegative (Pompili Internal, see \cite{Pompili:2014:CompMethod}).
                
                Finally, the TV denoising algorithm in \cref{alg:ONMF-TV_Line4} of all separated methods is based on a fast iterative shrinkage-thresholding Algorihtm (FISTA, see \cite{BeckTeboulle:2009:FISTA}) with a maximal iteration number of 100.
            
        \subsection{Combined Methods} \label{app:subsec:Parameter Choice of Combined Methods}
            In this Section, we describe the choice of the main hyperparameters of all combined methods. As in the previous section, \cref{tab:parameterChoice:CombinedMethods} shows the selected parameters for the numerical experiments in \cref{sec: Numerical Experiments} along with the used stopping criteria (stopCrit) and initialization methods (initMethod).
            
            \begin{table}[tbhp]
                \caption{Parameter choice of the combined methods for the numerical experiments with the MALDI dataset.}
                \ra{1.1}
                \resizebox{\columnwidth}{!}{%
                \pgfplotstabletypeset[
                    col sep=&, row sep=\\,
                    string type,
                    every head row/.style={
                        before row={
                            \toprule
                        },
                        after row=\midrule,
                    },
                    every last row/.style={
                        after row=\bottomrule},
                    every even row/.style={
                        before row={\rowcolor[gray]{0.9}}
                    },
                    columns/method/.style       ={column name=\textbf{Method}, column type = l},
                    columns/sigma1/.style       ={column name=$\sigma_1$},
                    columns/sigma2/.style       ={column name=$\sigma_2$},
                    columns/tau/.style          ={column name=$\tau$},
                    columns/epsTV/.style        ={column name=$\varepsilon_{\TV}$},
                    columns/subRat/.style       ={column name=$s_r$},
                    columns/stopCrit/.style     ={column name=stopCrit, column type = l},
                    columns/iMax/.style         ={column name=$i_{\text{max}}$},
                    columns/initMethod/.style   ={column name=initMethod, column type = l},
                    ]{
                        method                  & sigma1    & sigma2    & tau                   & epsTV             & subRat  & stopCrit   & iMax                 & initMethod \\
            			\texttt{ONMFTV-MUL1}    & $0.5$       & $0.5$       & $5\cdot 10^{-3}$      & $\sqrt{1\cdot 10^{-5}}$  & --      & maxIt      & $8\cdot 10^{2}$      & SVD\\
            			\texttt{ONMFTV-MUL2}    & $1$       & --        & $1\cdot 10^{-3}$      & $\sqrt{1\cdot 10^{-5}}$  & --      & maxIt      & $7\cdot 10^{2}$      & SVD\\
            			\texttt{ONMFTV-PALM}    & $0.1$     & $0.1$     & $1\cdot 10^{-1}$      & --                & --      & maxIt      & $4\cdot 10^{2}$      & SVD\\
            			\texttt{ONMFTV-iPALM}   & $0.1$     & $0.1$     & $1\cdot 10^{-1}$      & --                & --      & maxIt      & $3\cdot 10^{2}$      & SVD\\
            			\texttt{ONMFTV-SPRING}  & $0.1$     & $0.1$     & $1\cdot 10^{-4}$      & --                & $40$    & maxIt      & $1\cdot 10^{2}$      & SVD\\
                    }
                }
        		\label{tab:parameterChoice:CombinedMethods}
            \end{table}
            All combined methods are initialized based on the SVD approach described in \cref{subsec: Setup} and are stopped until $i_{\text{max}}$ iterations are reached (maxIt). Regarding the initialization methods, kmeans++ and the SVD method were tried out and one of them was finally chosen based on the quality of the results as it is the case for the separated methods. For \texttt{ONMFTV-SPRING}, we choose $s_r = 40$ according to \cite{DriggsEtal:2020:SPRING}.
            
            Furthermore, the computation of the proximal mapping for the combined methods in \cref{subsubsec: Promximal Alternating Linearized Minimization} is also based on the FISTA algorithm \cite{BeckTeboulle:2009:FISTA}, where the maximal iteration number is reduced to 5.
    
    
    
    \section{Further Experiments}\label{app:sec:Further Experiments}
\begin{figure}[tbhp]
            \centering
            \subfloat[Van Dongen criterion (VD)]{\label{fig:boxPlots:VD-VI:VD}\includegraphics[width=1.0\textwidth]{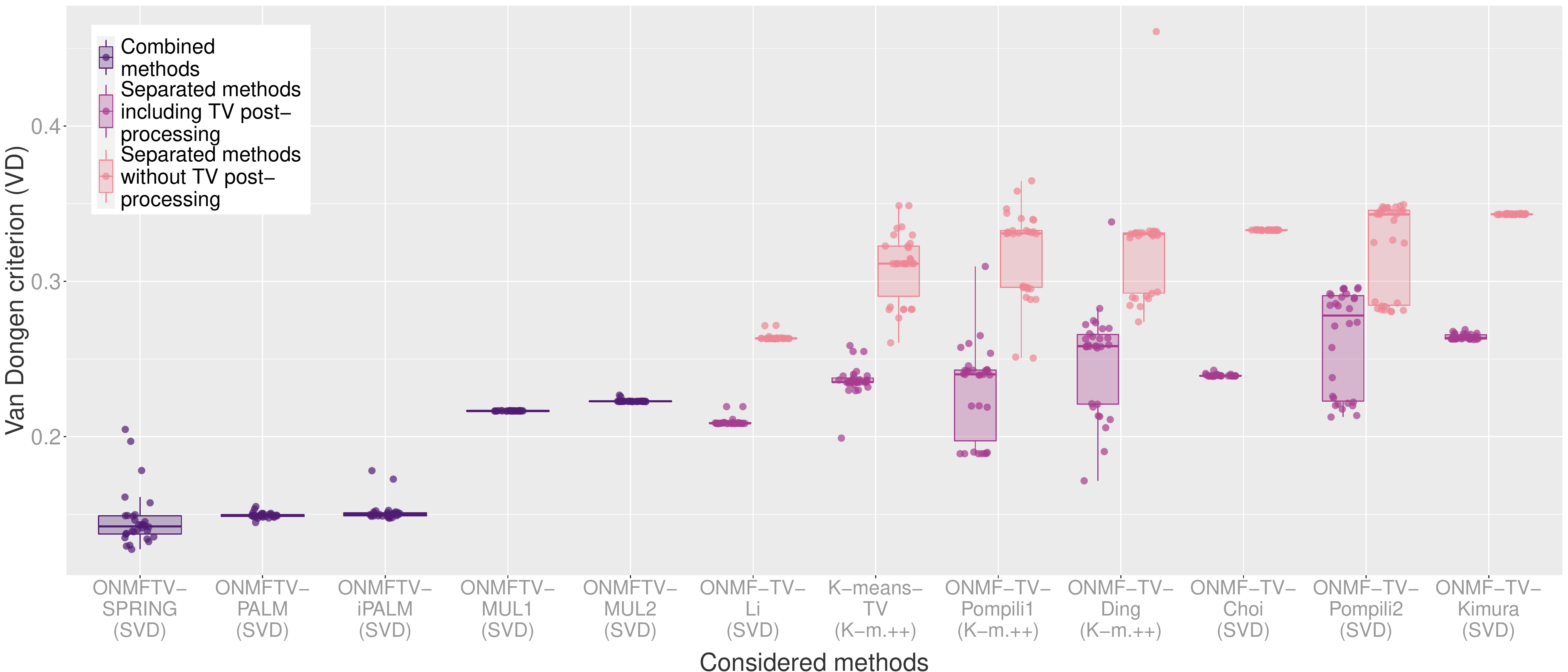}}\\
            \subfloat[Variation of information (VI)]{\label{fig:boxPlots:VD-VI:VI}\includegraphics[width=1.0\textwidth]{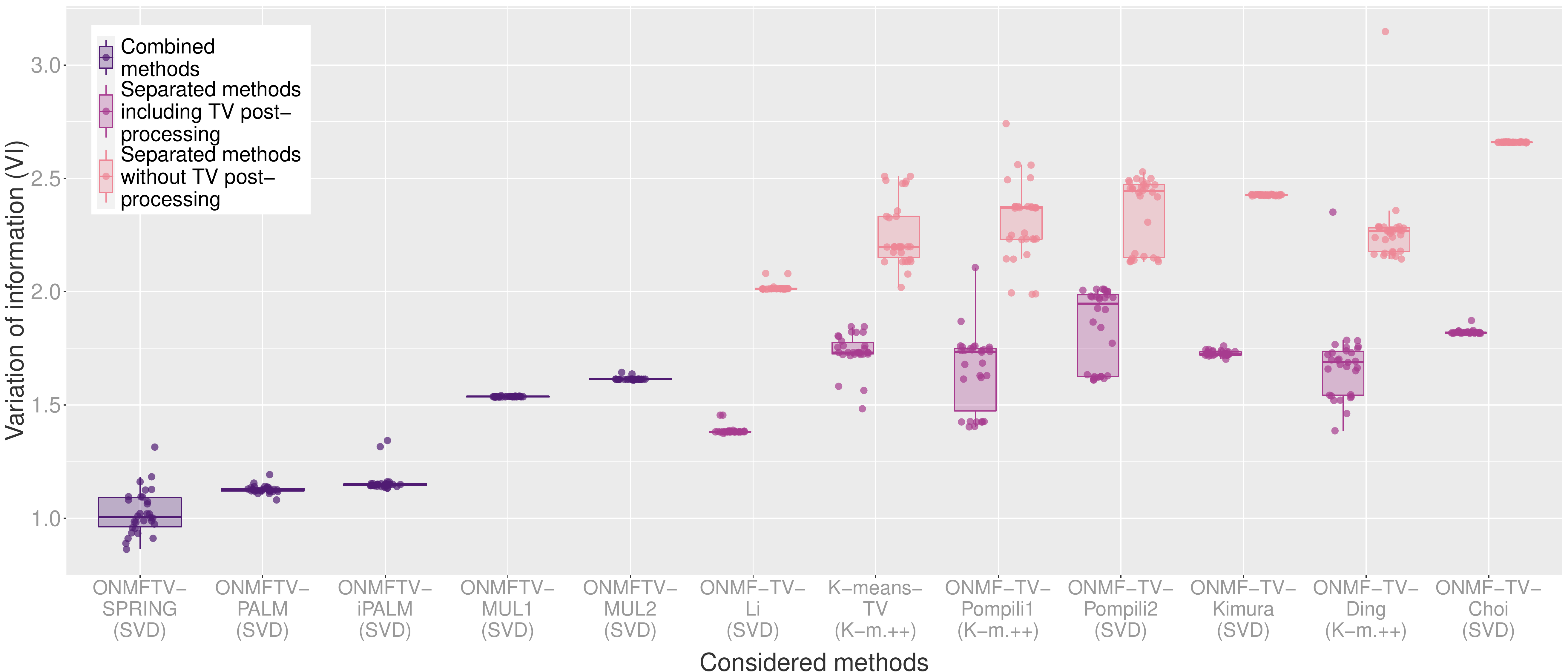}}
            \caption{Box plots of the van Dongen criterion (VD) and the variation of information (VI) of all performed experiments.}
            \label{fig:boxPlots:VD-VI}
        \end{figure}

        In this section, we present further experiments, which were done on the considered MALDI dataset. \cref{fig:boxPlots:VD-VI} shows the results of all performed experiments in terms of the van Dongen criterion (VD) and the variation of information (VI), which confirm the previous observations in \cref{subsec: Results and Discussion}.
    
    \section*{Acknowledgments}
    This project was funded by the Deutsche Forschungsgemeinschaft (DFG, German Research Foundation) within the framework of RTG ``$\pi^3$: Parameter Identification -- Analysis, Algorithms, Applications'' -- Project number 281474342/GRK2224/1. The used tissue sample is a SCiLS example dataset provided by Axel Walch (German Research Center for Environmental Health, Helmholtz Zentrum München). Furthermore, the author would like to thank Sibylle Hess (Eindhoven University of Technology) for fruitful discussions.
    
    \bibliographystyle{siam}
    \bibliography{main}

\end{document}